\documentclass[reqno, 11pt, a4paper]{amsart}
\usepackage[english]{babel}
\usepackage[letterpaper,top=2.8cm,bottom=2.8cm,left=3.5cm,right=3.5cm,marginparwidth=1.75cm]{geometry}

\DeclareFontFamily{U}{BOONDOX-calo}{\skewchar\font=45 }
\DeclareFontShape{U}{BOONDOX-calo}{m}{n}{
  <-> s*[1.05] BOONDOX-r-calo}{}
\DeclareFontShape{U}{BOONDOX-calo}{b}{n}{
  <-> s*[1.05] BOONDOX-b-calo}{}
\DeclareMathAlphabet{\mathcalboondox}{U}{BOONDOX-calo}{m}{n}
\SetMathAlphabet{\mathcalboondox}{bold}{U}{BOONDOX-calo}{b}{n}
\DeclareMathAlphabet{\mathbcalboondox}{U}{BOONDOX-calo}{b}{n}
\newcommand{\mcb}[1]{{\mathcalboondox #1}}

\usepackage{mathrsfs}

\usepackage{amsmath}
\usepackage{mathtools}
\usepackage{amsfonts}
\usepackage{graphicx}
\usepackage[colorlinks=true, allcolors=blue]{hyperref}
\usepackage{accents}
\usepackage{amssymb}

\usepackage{dsfont}
\usepackage{bbm}
\usepackage{indentfirst}

\newcommand{\lsim}{\lesssim}

\usepackage{amsthm}

\usepackage{tikz}
\usetikzlibrary{shapes.misc}
\usetikzlibrary{arrows}
\tikzset{cross/.style={cross out, draw=black, minimum size=2*(#1-\pgflinewidth), inner sep=0pt, outer sep=0pt},
cross/.default={0.2cm}}
\definecolor{blue-violet}{rgb}{0.54, 0.17, 0.89}
\usetikzlibrary{decorations.pathreplacing}

\begin{document}
\title[Exclusion processes with non-reversible boundary]{Exclusion processes with non-reversible boundary: hydrodynamics and large deviations}

\author[C. Landim]{Claudio Landim}
\address{Instituto de matemática pura e aplicada, Estrada Dona Castorina 110, J. Botanico, 22460 Rio de Janeiro, Brazil}
\email{landim@impa.br}

\author[J. Mangi]{João Pedro Mangi}
\address{Instituto de matemática pura e aplicada, Estrada Dona Castorina 110, J. Botanico, 22460 Rio de Janeiro, Brazil}
\email{mangi.joao@impa.br}

\author[B. Salvador]{Beatriz Salvador}
\address{Center for Mathematical Analysis, Geometry and Dynamical Systems, Instituto Superior Técnico, Universidade de Lisboa, Av. Rovisco Pais, 1049-001 Lisboa, Portugal.}
\email{beatriz.salvador@tecnico.ulisboa.pt}

\maketitle

\theoremstyle{plain}
\newtheorem{theorem}{Theorem}[section] 
\newtheorem{lemma}[theorem]{Lemma}
\newtheorem{corollary}[theorem]{Corollary}
\newtheorem{proposition}[theorem]{Proposition}

\theoremstyle{definition}
\newtheorem{definition}[theorem]{Definition}
\newtheorem{remark}[theorem]{Remark}
\newtheorem{assumption}[theorem]{Assumption}
\newtheorem{example}[theorem]{Example}

\makeatletter
\renewcommand{\theequation}{%
\thesection.\arabic{equation}}
\@addtoreset{equation}{section}
\makeatother

\makeatletter
\renewcommand{\p@enumi}{A}
\makeatother


\providecommand{\keywords}[1]
{
  \small	
  \textbf{\textbf{Keywords:  }} #1
}

\newcommand{\acks}{\textbf{Acknowledgements.}}

\newenvironment{acknowledgements}{%
  \renewcommand{\abstractname}{Acknowledgements}
  \begin{abstract}
}{%
  \end{abstract}
}

\allowdisplaybreaks

\begin{abstract} We consider a one-dimensional
exclusion dynamics in mild contact with boundary reservoirs. In the
diffusive scale, the particles' density evolves as the solution of the
heat equation with non-linear Robin boundary conditions. For
appropriate choices of the boundary rates, these partial differential
equations have more than one stationary solution. We prove the
dynamical large deviations principle.
\end{abstract}
\vspace{0.5cm} \keywords{Hydrodynamic Limit, Large Deviations,
Boundary-driven particle systems, Non-reversible
boundary, Non-linear Robin boundary conditions}

\vspace{1cm}

\section{Introduction}
Interacting particle systems in contact with boundary
reservoirs have attracted a lot of attention these last decades as a
model for out of equilibrium systems which exhibit long-range
correlations \cite{spohn1983long}, non-local thermodynamical
variables, and dynamical phase transition
\cite{bertini2015macroscopic}.

In this article, we address a long-standing open problem which
consists in deriving the hydrodynamic limit and the dynamical large
deviations principles of systems in which the boundary and the bulk
stationary states are different.

It is well known that the canonical steady states of {the} simple exclusion
dynamics are the uniform measures concentrated in configurations with
a fixed number of particles, while the grand canonical ones are the
Bernoulli product measures. By superposing a simple symmetric
exclusion dynamic, say a one-dimensional grid to simplify the
exposition (although this hypothesis is not necessary), in contact at
the boundary with a non-conservative dynamics whose stationary states
are product Bernoulli measures with distinct densities at each end of
the interval, we obtain a dynamics where a particle flow is observed
from the end  {point} with larger density to the other.

Since the boundary and the bulk dynamics have the Bernoulli product
measures as steady states, it is expected that in the steady state of
the full dynamics, near the boundary, the distribution of particles is
close to a Bernoulli measure with density determined by the
non-conservative boundary dynamics.
Therefore, even if the steady state of the global dynamics is not
explicitly known, it is still possible to apply the entropy or the
relative entropy methods \cite{KL99} to derive the hydrodynamic limit.  {Thus, for its derivation, one can
use} as a reference measure a Bernoulli product measure that smoothly
interpolates between the two densities determined by the boundary
dynamics, see \cite[Proposition 2.3]{landim1998driven}.

When the steady state of the boundary dynamics is not a Bernoulli
measure, near the boundary, the distribution of particles in the
steady state is not known. Moreover, a good reference measure is not available
and the proof of the hydrodynamic behavior of the particle system is
still an open problem, with few exceptions for very particular choices of boundary dynamics \cite{Nahum2020,de2012truncated} or  dynamics in which
duality techniques can be applied \cite{erignoux2018stationary,
erignoux2018hydrodynamic, erignoux2019hydrodynamics}.

By decreasing the intensity of the system's interaction with the
reservoirs, the bulk dynamics prevail and it is expected that the
distribution of particles near the boundary will be close to a
Bernoulli measure with a density determined by the boundary dynamics. In this case, it is possible using
as a reference measure a Bernoulli product measure that smoothly
interpolates between the two densities determined by the boundary
dynamics to derive the hydrodynamic limit. The
boundary conditions of the hydrodynamic equations, which were of
Dirichlet type in the original model, become of mixed type.

By taking a double limit, letting first the mesh of the grid decrease
to $0$ and then the intensity of the interaction of the system with
the boundary reservoirs to increase to $+\infty$, one could recover
the hydrodynamic limit and the dynamic and static large deviations
rate functional (by taking a $\Gamma$-convergence) of the original
dynamics.  {This is the outline of a project in progress for which this article can be though as the first step of the program.}

In this article, we investigate the hydrodynamic limit and the dynamic
large deviations of a class of exclusion processes evolving on the
one-dimensional lattice with mild interactions with two stochastic
boundary reservoirs. The model is a result of a superposition of two dynamics: in the interior points, the usual exclusion rule, i.e. only one particle is allowed per site with the nearest neighbor interaction, while the boundary dynamics is a creation and annihilation of particles in a window of fixed size $l$ attached to each of the extreme points of the boundary. This choice of boundary dynamics is non-conservative and considered under very general rates. This last property is the main difference between the class of models we consider here and the exclusion process with non-reversible boundary studied in \cite{erignoux2018hydrodynamic,Nahum2020, erignoux2018stationary, de2012truncated}, where very specific choices of boundary dynamics were made.

One of the main features of the model considered
here is that the hydrodynamic equation may have, according to the
choice of the rates, multiple stationary solutions. This property
raises a series of questions, such as the identification of stable and
unstable density profile equilibria, the transition times between
stable equilibria, and the metastable behavior of the system, similar
problems, and probably simpler, to those examined for the
superposition of the Glauber and Kawasaki dynamics \cite{faris1982large, landim2018hydrostatics, farfan2019static}.
The main novelty of this work is the proof of Theorem \ref{uniqueness_PDE_hydro}, which allowed us to carry on the usual techniques to prove both Theorems \ref{th_hydrodynamics} and \ref{main_result_large_dev}.

The paper is organized as follows: 
in Section \ref{intro} we introduce the model and state the main results, Theorems \ref{th_hydrodynamics} and \ref{main_result_large_dev}. Section \ref{section_hydrodynamics} is dedicated to the proof of Theorem \ref{th_hydrodynamics} and is divided into three parts: tightness, characterization of the limit points, and a small comment on uniqueness, which is later proved in Appendix \ref{appB}. Next, Section \ref{large_dev_proof} contains the main ingredients to prove Theorem \ref{main_result_large_dev} and also its own proof. We prove the upper bound and lower bounds in Section \ref{upperbound} and \ref{lowerbound}, respectively, and the I-density and properties of the large deviations rate functional in Section \ref{I-density}. Finally, in Appendix \ref{appA} we estimate the Dirichlet form of the process and prove a replacement lemma that is needed for the characterization of the limit point in the hydrodynamic limit result. In Appendix \ref{model_1.5} we provide an example of a choice of boundary rates for which we construct a model, which we call Exclusion $l3$, where we have more than one stationary profile but the hydrodynamic equation has a unique weak solution. We conclude with Appendix \ref{appB} where we not only show the uniqueness and properties of the solution to the initial value problem that arises from the hydrodynamic limit but also prove the uniqueness of the solution to the initial value problem associated with the perturbed process that arises on the proof of the lower bound of the dynamical large deviations principle.

\section{Model and main results}\label{intro}
\subsection{The model.} \label{model}
Given $N \in \mathbb{N}$, let $\epsilon_N = 1/N$ and $\tau_N = (N-1)/N$. We define an exclusion process in $\Lambda_N:= \{\epsilon_N, \dots, \tau_N\}$, which we call bulk, with state space $\Omega_N:= \{0,1\}^{\Lambda_N}$. We denote by $\Lambda_N^\circ$ the set $\Lambda_N \setminus \tau_N$. For every $x \in \Lambda_N$ and $\eta \in \Omega_N$, $\eta(x)$ represents the number of particles of the configuration $\eta$ at the site $x$ and for simplicity we will call it occupation variable. Also, fix $l \in \{1, \dots, N-1\}$, which represents the size of the left and right windows of the interaction of the boundary dynamics, that we will prescribe later, with the bulk. Let $\Sigma^-_l:=\{\epsilon_N,...,l\epsilon_N\}$ and $\Sigma^+_l:=\{(N-l)\epsilon_N,...,\tau_N\}$. Consider $\{0,1\}^{\Sigma^-_l}$, respectively $\{0,1\}^{\Sigma^+_l}$, which corresponds to the space of configuration of particles of size $l$ close to the boundary point $0$, respectively $1$.

For $1 \leq j < k \leq N-1$ and $\eta \in \Omega_N$, we define the projection $\Pi_{j,k}$ as $$(\Pi_{j,k} \eta) (x) = \eta(x),\;\text{ for } x \in \{j\epsilon_N, \dots, k \epsilon_N\}.$$The configuration $\Pi_{j,k} \eta$ is the one obtained from $\eta$ when restricting $\eta$ to the values of $\eta(x)$ for $x \in \{j\epsilon_N, \dots, k \epsilon_N\}$. Given $x \in$ $\Lambda_N^\circ$ and $\eta \in \Omega_N$, we define $\sigma^{x,x+\epsilon_N} \eta \in \Omega_N$ as, for every $z \in \Lambda_N$, $$[\sigma^{x,x+\epsilon_N} \eta](z) = \eta(z) \mathds{1}_{\{z \neq x, x+\epsilon_N\}} + \eta(x+\epsilon_N) \mathds{1}_{\{z = x\}} + \eta(x) \mathds{1}_{\{z = x + \epsilon_N \}}.$$ The configuration $\sigma^{x,x+\epsilon_N} \eta$ is the new configuration obtained from $\eta$ when swapping the number of particles in $x$ with the number of particles in $x + \epsilon_N$. Given $\xi \in \{0,1\}^{\Sigma^-_l}$, we also define
$\xi||\Pi_{l+1,N-1} \eta \in \Omega_N$ where 
$$
(\xi||\Pi_{l+1,N-1} \eta)(x) := \xi(x) \mathds{1}_{\{x = \epsilon_N,\dots,l\epsilon_N\}} + \eta(x)\mathds{1}_{\{x =(l+1)\epsilon_N,\dots,\tau_N\}},
$$
for every $x \in \Lambda_N$. Given $\xi \in \{0,1\}^{\Sigma^+_l}$, we define
$\Pi_{1,N-l-1} \eta || \xi \in \Omega_N$ where 
$$
(\Pi_{1,N-l-1} \eta || \xi)(x) := \eta(x)\mathds{1}_{\{x =\epsilon_N,\dots,(N-1-l)\epsilon_N\}}+ \xi(x) \mathds{1}_{\{x = (N -l)\epsilon_N,\dots,\tau_N\}},
$$
for every $x \in \Lambda_N$.
The configurations $\xi||\Pi_{l+1,N-1} \eta$ and
$\Pi_{1,N-l-1} \eta || \xi$ represent the concatenation of the configuration $\xi$ with $\Pi_{l+1,N-1} \eta$ and $\Pi_{1,N-l-1} \eta$ with $\xi$, respectively.
The particle system we are interested in, and that we denote by $\{\eta_{t};{t\geq 0}\}$ has a Markov generator $\mathscr{L}_N$ given by
\begin{equation}
    \mathscr{L}_N = \mathscr{L}^{-}_N + \mathscr{L}^{bulk}_N + \mathscr{L}^{+}_N,
\end{equation}
where, for every $f: \Omega_N \to \mathbb{R}$ and $\eta \in \Omega_N$,
\begin{align}
    (\mathscr{L}^{bulk}_N f)(\eta) = N^2 \sum_{x \in \Lambda_N^\circ} [f(\sigma^{x,x+\epsilon_N} \eta) - f(\eta)],
\end{align}
\begin{align}
    (\mathscr{L}^{-}_N f)(\eta) = N \sum_{\xi \in \{0,1\}^{\Sigma^-_l}} R^-\left(\Pi_{1,l} \eta, \xi \right) [f(\xi||\Pi_{l+1,N-1} \eta) - f(\eta)],
\end{align}
and
\begin{align}
    (\mathscr{L}^{+}_N f)(\eta) = N \sum_{\xi \in \{0,1\}^{{\Sigma^+_l}}} R^+\left(\Pi_{N-l,N-1} \eta, \xi \right) [f(\Pi_{1,N-1-l} \eta || \xi) - f(\eta)],
\end{align} where $R^-$ and $R^+$ are non-negative functions that will determine the boundary rates.


This model is the result of a superposition of two dynamics: the symmetric simple exclusion (Kawasaki's dynamics) on the bulk and a creation and annihilation dynamics at the points of a neighborhood of size $l$ of each of the boundary points, $0$ and $1$, with transition rates given by $R^-$ and $R^+$, respectively. It is important to remark that, since these rates are general, we can inject and extract more than one particle at the same time on the system provided that we obtain a configuration with no more than one particle at each site of $\Lambda_N$. In \cite{Nahum2020}, this model was considered with a specific choice of the rates $R^-$ and $R^+$. Here, we will consider the general case with the only assumption that the rates $R^{\pm}$ must be chosen in such a way that the dynamics in each of the windows, i.e., the superposition of the boundary dynamics with the Kawasaki dynamics restricted to $\{0,1\}^{\Sigma^-_l}$, respectively $\{0,1\}^{\Sigma^+_l}$, must be an irreducible dynamics. This assumption implies the irreducibility of the process.



\subsection{Notation of function spaces}
Fix $T>0$ and a finite interval of time $[0,T]$. Let $\Omega_T := [0,T] \times [0,1]$. We will now introduce some function spaces that we will use throughout the article.
\begin{itemize}
    \item For every $m \in \mathbb{N} \cup \{0\}$, we denote by $C^{m}([0,1])$ the set of continuous functions defined on $[0,1]$ that are $m$ times continuously differentiable. For $m=0$, instead of writing $C^0([0,1])$, we will simply write $C[0,1]$;
    \item For every $m,n \in \mathbb{N} \cup \{0\}$, we denote by $C^{m,n}([0,T]\times(0,1))$ the set of real-valued functions defined on $[0, T]\times(0,1)$ that are $m$ times continuously differentiable in time, and $n$ times continuously differentiable on space. Moreover, for every $m,n \in \mathbb{N} \cup \{0\}$, we denote by $C^{m,n}(\Omega_T)$ the set of functions defined on $[0, T]\times[0,1]$ that belong to $C^{m,n}([0,T]\times(0,1))$ and are continuously extended to the boundary;
    \item For every $m,n \in \mathbb{N} \cup \{0\}$, we denote by $C^{n,m}_{0}(\Omega_T)$ be the set of all functions $H:[0,T]\times[0,1]$ that below to $C^{n,m}([0,T]\times(0,1))$ and which are continuously extended to the boundary as $H_t(0) = H_t(1) = 0$ for any $t \in [0,T]$;
    \item We denote by $H^1([0,1])$ the Sobolev space of measureable functions $u:[0,1]\to\mathbb{R}$ with weak derivative $\nabla u$ in $L^2([0,1])$. The space $H^1$ endowed with the scalar product $\langle\cdot,\cdot\rangle_1$ given, for every $f,g \in H^1$, by 
\begin{equation*}
\langle f,g\rangle_1 = \int_0^1 f(x)g(x)\;dx + \int_0^1 \nabla f(x) \nabla g(x)\;dx,
\end{equation*} is an Hilbert space and we denote by $||.||_1$ the norm induced by the inner product in $H^1$.
    \item We introduce the Sobolev space $L^2([0, T], H^1([0,1]))$ which is the Banach space of measurable functions $u: [0, T] \to H^1([0,1])$ endowed with the norm $||\cdot||_{2,1}$ which is given by $||u||_{2,1}:= \left(\int_0^t ||u(s)||^2_1 ds\right)^{1/2}$.
\end{itemize}

\subsection{Hydrodynamic Limit} \label{hydrodynamic}
We observe that given an initial configuration $\eta_0 \in \Omega_N$, the evolution of the Markov process $\{\eta_{t}; t \geq 0\}$ is a trajectory in $\Omega_N$. For a metric space $\mathbb{X}$, we denote by $\mathcal{D}([0,T],\mathbb{X})$ the space of càdlàg trajectories endowed with the Skorohod topology with values in $\mathbb{X}$. Let $\mcb {M}$ be the space of non-negative measures on $[0,1]$ with total mass bounded by one, equipped with the weak topology. For every configuration  $\eta \in \Omega_N$, we define the empirical measure $\pi^{N}(\eta,du)$ by 
\begin{equation*}\label{MedEmp}
\pi^{N}(\eta, du):=\dfrac{1}{N-1}\sum _{x\in\Lambda_N }\eta(x)\delta_{x}\left( du\right) \in \mcb{M},
 \end{equation*}
where $\delta_{b}$ is a Dirac measure of $b \in \mathbb{R}$. To simplify notation, we will use $\pi^N$ to denote $\pi^N(\eta,du)$. For every $H: [0,1] \rightarrow \mathbb{R}$, we denote the integral of $H$ with respect to $\pi^N$ by $\langle \pi^N, H \rangle$
\begin{align}~\label{eq:integral_G_emp}
\langle \pi^N, H \rangle := \frac{1}{N-1} \sum_{x \in \Lambda_N } \eta(x) H \left( x  \right).
\end{align}

Let $\gamma: [0,1] \rightarrow [0,1]$ be a measurable function. We say that a sequence $(\nu_N)_{N \geq 1}$ of probability measures on $\Omega_N$ is associated to the profile $\gamma$ if for every function $H \in C([0,1])$ and every $\delta >0$, it holds
\begin{equation*}
\lim_{N \rightarrow \infty} \nu_N \Big(\; \Big| \langle \pi^N, H \rangle - \int_0^1 H(x) \gamma (x) dx \Big| > \delta \Big) =0. 
\end{equation*}
For every $N \in \mathbb{N}$, let $\mathbb{P} _{\nu_{N}}$ be the probability measure on $\mathcal{D}([0,T],\Omega_N)$ induced by the Markov process $\{\eta_{t};{t\geq 0}\}$ and by the initial measure $\nu_{N}$. The expectation with respect to $\mathbb{P}_{\nu_{N}}$ is denoted  by $\mathbb{E}_{\nu_{N}}$. We also define $\pi^{N}_{t}(\eta, du):=\pi^{N}(\eta_{t}, du)$ a trajectory of measures in $\mcb{M}$. Let $(\mathbb{Q}_{\nu_N})_{N \in \mathbb{N}}$ be the sequence of probability measures on $\mathcal{D}([0,T],\mcb M)$ induced by the Markov process $(\pi_{t}^{N})_{0 \leq t \leq T}$ and $\mathbb{P} _{\nu_{N}}$, i.e. $\mathbb{Q}_{\nu_N}(\cdot) := \mathbb{P}_{\nu_{N}}( \pi^{N} \in \cdot)$. 

Denote by $h_{\pm}(\eta)$ the rate at which the boundary density changes, which is defined as
\begin{equation} \label{def_h_plus_minus}
h_{\pm}(\eta):=\mathscr{L}^{\pm}\left(\sum_{x\in\Sigma^\pm_l}\eta(x)\right)=\sum_{\xi\in\{0,1\}^{\Sigma^\pm_l}}R^{\pm}(\Pi\eta,\xi)\sum_{x\in\Sigma^\pm_l}(\xi(x)-\eta(x)),
\end{equation} where $\Pi\eta$ should be interpret as $\Pi_{1,l}\eta$, respectively $\Pi_{N-l,N-1}\eta$, in the case of the $-$ (minus), respectively $+$(plus), sign.
Denoting by $\nu_\alpha$ the Bernoulli product measure with density $\alpha$, set 
\begin{equation} \label{def_F_pm}
    F_{\pm}(\alpha): =\mathbb{E}_{\nu_\alpha}(h_{\pm}).
\end{equation}
We will denote by $B_{\pm}(\alpha)$ and $D_{\pm}(\alpha)$ the creation and destruction rates, respectively, given by 
\begin{equation}\label{creation_rate}
    B_{\pm}(\alpha)=\mathbb{E}_{\nu_\alpha}\Big[\sum_{\xi\in\{0,1\}^{\Sigma_l^{\pm}}}R^{\pm}(\Pi\eta,\xi)\Big(\sum_{x\in\Sigma_l^{\pm}}\xi(x)-\eta(x)\Big)^+\Big]
\end{equation}
\begin{equation}\label{destruction_rate}
    D_{\pm}(\alpha)=\mathbb{E}_{\nu_\alpha}\Big[\sum_{\xi\in\{0,1\}^{\Sigma_l^{\pm}}}R^{\pm}(\Pi\eta,\xi)\Big(\sum_{x\in\Sigma_l^{\pm}}\xi(x)-\eta(x)\Big)^-\Big],
\end{equation} where $\Pi\eta$ should again be interpret as $\Pi_{1,l}\eta$, respectively $\Pi_{N-l,N-1}\eta$, in the case of the $-$ (minus), respectively $+$(plus), sign.
Then, $F_{\pm}(\alpha)=B_{\pm}(\alpha)-D_{\pm}(\alpha)$.

\begin{definition}[Weak solution] \label{def_weak_solution_complete} Given $\gamma:[0,1] \to [0,1]$ a measurable function and $F_\pm$ the functions defined in \eqref{def_F_pm}, we say that $\rho:\Omega_T \to \mathbb{R}$ is a weak solution of \begin{align} \label{PDEproblem}
    \begin{cases}
    \partial_t \rho_t(x) = \Delta \rho_t(x),  &(t,x) \in (0,T)\times (0,1), \\
    \nabla \rho_t(0) = - F_-(\rho(t,0)), &t \in (0,T),\\
    \nabla \rho_t(1) = F_+(\rho(t,1)), &t \in (0,T),\\
    \rho_0(x) = \gamma(x), &x \in [0,1],
    \end{cases} 
\end{align}
if the following two conditions are satisfied:
\begin{enumerate}
    \item $\rho \in L^2([0,T];H^1([0,1]))$ ;
    \item for every $H \in C^{1,2}( {\Omega_T})$,
     {\begin{align} \nonumber
        &\int_0^1 \rho_T(x) H_T(x)dx - \int_0^1 \gamma(x) H_0(x) dx  =\int_0^T \int_0^1 \rho_s(x) \partial_s H_s(x) dx ds \\ \label{def_weak_sol_eq_hydrodynamcis}
        &- \int_0^T \int_0^1 \nabla \rho_s(x) \nabla H_s(x) dx ds + \int_0^T \left\{ F_+[ \rho_s(1)] H_s(1) + F_-[ \rho_s(0)] H_s(0) \right\} ds.
    \end{align}}
\end{enumerate}
\end{definition}

 {The initial value problem \eqref{PDEproblem} has one weak solution which is also unique. Its existence is a consequence of the tightness of $(\mathbb{Q}_{\nu_N})_{N \in \mathbb{N}}$ which is shown in Section \ref{tightness_section} and its uniqueness is shown in Appendix \ref{appB}.}

\begin{theorem}{(Hydrodynamic Limit)}\label{th_hydrodynamics}
Let $\gamma: [0,1] \rightarrow [0,1]$  {be a measurable function and $(\nu_N)_{N \geq 1}$ a }sequence of probability measures on $\Omega_N$ which are associated to the initial profile $\gamma$. Then,  {for every $G \in C([0,1])$, every $t \in [0,T]$ and every $\delta > 0$,}
\begin{equation}\label{limHidreform}
 \lim _{N\to\infty } \mathbb P_{\nu_{N}}\Big(\;\Big| \langle \pi^N_t,G\rangle - \int_0^1 G(x) \rho_t(x) dx\Big| > \delta \Big)= 0,
\end{equation}
where  $\rho_{t}(\cdot)$ is the unique weak solution of \eqref{PDEproblem}.
\end{theorem}

\begin{remark}{(Stationary Solutions)} \label{stationary_rmk}
A density profile $\rho:[0,1]\to[0,1]$ is a stationary solution of equation \eqref{PDEproblem} if and only if there exists $\alpha,\beta\in[0,1]$ such that
\begin{equation}\label{stationary_rmk_equation}
     {- F_{-}(\alpha) = F_{+}(\beta) = \beta - \alpha}
\end{equation}
and,  {in this case,}  {$\rho(x)= (\beta - \alpha)x+\alpha$}. A particular case occurs when $F_{-}(\alpha)=F_{+}(\alpha)=0$
for some $0\leq \alpha\leq1$, in which case the constant profile $\rho(x)=\alpha$ is a stationary solution of  {\eqref{PDEproblem}}.

If equation \eqref{stationary_rmk_equation} admits more than one solution, then there are  {more than one stationary solution} and the quasi-potential (that is, the large deviations rate function
of the stationary problem) exhibits at least two critical points. In Appendix \ref{model_1.5}, we present a choice of rates for which this behavior can be observed. This is one of the main features of the model considered in this paper.
\end{remark}

\subsection{Large Deviations} \label{largedeviations}
Let $\mcb{M}_{ac} \subset \mcb{M}$  {be} the set of all measures which are absolutely continuous with respect to the Lebesgue measure and  {whose densities take} values in $[0,1]$.

For $T>0$, define the energy $Q_{[0,T]}(\pi)$, for every $\pi \in \mathcal{D}([0,T],\mcb M_{ac})$ where $\pi(t,dx) = \rho_t(x)dx$ for some density $\rho$, by
\begin{equation}
    Q_{[0,T]}(\pi)=\sup_{H\in  {C^{1,2}_{0}(\Omega_T)}}\left\{\int_{0}^{T}\int_{0}^{1}\rho_t(x)\nabla H(t,x)\;dxdt-\frac{1}{2}\int_{0}^{T}\int_{0}^{1}\sigma( {\rho_t(x)})H(t,x)^2\;dxdt\right\}
\end{equation}
where $\sigma:[0,1]\to\mathbb{R}$ is the mobility of the exclusion process defined by $\sigma(a)=a(1-a)$. Note that $Q_{[0,T]}$ is convex and lower semi-continuous. 

\begin{remark}
To simplify notation, when $\pi(t,dx)=\rho_t(x)\;dx$ for some density $\rho$, we will write $Q_{[0,T]}(\rho)$ instead of $ Q_{[0,T]}(\pi)$, and  {we will} follow this convention for every functional of $A \subset \mathcal{D}([0,T],\mcb M_{ac})$, with $A$ closed.
\end{remark}

Denote by $\mathcal{D}_{\mathcal{E}}([0,T],\mcb{M}_{ac})$ the set of trajectories in $\mathcal{D}([0,T],\mcb M_{ac})$ with finite energy, and, for a fixed initial profile $\gamma$, denote by $\mcb{D}_{\gamma,\mathcal{E}}([0,T],\mcb{M}_{ac})$ the set of trajectories with finite energy and starting from $\gamma$. Define, for each $H\in C^{1,2}( {\Omega_T})$, the functional $J_{T,H}:\mathcal{D}_{\mathcal{E}}([0,T],\mcb{M}_{ac})\to\mathbb{R}$ given by
 {\begin{equation}\label{def_J}
\begin{split}
    &J_{T,H}(\rho)=\langle \rho_T,H_T\rangle-\langle \gamma,H_0\rangle-\int_0^T\langle \rho_t,\partial_tH_t\rangle\;dt-\int_0^T\langle \nabla \rho_t,\nabla H_t\rangle\;dt\\
    &-\int_0^T\langle \sigma(\rho),(\nabla H_t)^2\rangle\;dt- \int_0^T\{\mathfrak{b}^-(\rho_t(0),H_t(0))+\mathfrak{b}^+(\rho_t(1), H_t(1))\}\;dt,
    \end{split}
\end{equation}}
where $\langle\cdot,\cdot\rangle$ is the inner product in $L^2[0,1]$,  { and the functions $\mathfrak{b}^{\pm}(\cdot,\cdot)$ are given by
\begin{equation} \label{bfrakplus}
    \mathfrak{b}^{-}(\alpha,M):=\mathbb{E}_{\nu_\alpha}\left[\sum_{\xi\in {\{0,1\}^{\Sigma^-_l}}}R^{-}(\Pi_{1,l}\eta,\xi)\{e^{M\sum_{x\in {\Sigma^-_l}}[\xi(x)-\eta(x)]}-1\}\right],
\end{equation} and
\begin{equation} \label{bfrakminus}
    \mathfrak{b}^{+}(\alpha,M):=\mathbb{E}_{\nu_\alpha}\left[\sum_{\xi\in {\{0,1\}^{\Sigma^+_l}}}R^{+}(\Pi_{N-l,N-1}\eta,\xi)\{e^{M\sum_{x\in  {\Sigma^+_l}}[ {\xi(x)}-\eta(x)]}-1\}\right].
\end{equation} 
Note that we can also write $\mathfrak{b}^{\pm}(\alpha,M)$ as 
 {\begin{equation*}
    \mathfrak{b}^{ {\pm}}(\alpha,M):=\sum_{k=1}^{l}(e^{Mk}-1)B_k(\alpha) + \sum_{k=-l}^{-1}(e^{Mk}-1)D_{k}(\alpha),
\end{equation*}
where
\begin{equation} \label{B_D}
\textrm{ for every $k \in \{1,\dots,l\}$, } \quad B_k(\alpha) = \mathbb{E}_{\nu_\alpha}\left[R^{ {\pm}}_k(\eta)\right] \quad \textrm{ and } \quad D_{-k}(\alpha) = B_{k}(\alpha),
\end{equation} with $R^{ {\pm}}_k(\eta)=\sum_{\xi\in {\{0,1\}^{\Sigma^\pm_l}}}R^{\pm}(\Pi\eta,\xi)\boldsymbol{1}\{\sum_{x\in {\Sigma^\pm_l}} {\xi(x)}=\sum_{x\in {\Sigma^\pm_l}} {\eta(x)}+k\}$  {represents the rate at which  {$k$} particles are created},} and $\Pi\eta$ should be understood  {as usual.}
}

\begin{remark}
     {As} in the case of the exclusion process with linear Robin boundary conditions, we need to define the functional $J_{T,H}$ in $\mathcal{D}_{\mathcal{E}}([0,T],\mcb M_{ac})$ instead of $\mathcal{D}([0,T],\mcb M_{ac})$, as it is the case of Dirichlet boundary conditions. Since the boundary densities are not fixed by the dynamics, we need to restrict ourselves to trajectories with finite energy so that the density profile is well defined  {at} the boundary. 
\end{remark}
 {We extend the definition of $J_{T,H}$ to $\mathcal{D}([0,T],\mcb{M}_{ac})\to[0,\infty]$ by setting
    \begin{equation*}
        J_{T,H}(\pi) = + \infty \quad \textrm{ if } \quad \pi \notin \mathcal{D}_{\mathcal{E}}([0,T],\mcb M_{ac}).
    \end{equation*}}
    
Now define the functional $I_{[0,T]}: {\mathcal{D}}([0,T],\mcb{M}_{ac})\to[0,\infty]$ by
\begin{equation*}
I_{[0,T]}(\pi)=\sup_{H\in C^{1,2}(\Omega_T)}J_{T,H}(\pi).
\end{equation*}
For a fixed density profile $\gamma\in\mcb{M}_{ac}$ we define the \textbf{\textit{rate functional}} as $I_{[0,T]}(\cdot|\gamma): {\mathcal{D}}([0,T],\mcb{M}_{ac})\to\mathbb{R}$, given by
\begin{equation}
I_{[0,T]}(\cdot|\gamma) =
\begin{cases}
I_{[0,T]}(\pi),\;\;\text{if}\;\pi\in \mcb{D}_{\gamma,\mathcal{E}}([0,T],\mcb{M}_{ac})\\
+\infty,\;\text{otherwise}.
\end{cases}
\end{equation}
\begin{theorem}[Large Deviations Principle] \label{main_result_large_dev}
Assume that,  {for every $k \in \{1,\dots,l\}$, $B_k(\cdot)$ and $D_{-k}(\cdot)$ defined in \eqref{B_D}} are concave functions  {and that the initial density profile $\gamma\in C^{2+\beta}([0,1])$, for some $\beta>0$}. 
Fix $T > 0$ and a measure $\pi(dx) = 
\gamma(x) dx$ in $\mcb{M}_{ac}$. Consider  {the} sequence $\eta^N$ of configurations associated to $\gamma$. Then the measure $\mathbb{Q}_{\eta^N}$ satisfies a large deviation principle with speed $N$ and rate function $I_{[0,T]}(\cdot| \gamma)$. Namely, for each closed set
$\mathcal{C} \subset \mathcal{D}([0,T],\mcb M)$ and each open set $\mathcal{O} \subset \mathcal{D}([0,T],\mcb M)$,
\begin{align}
    \limsup_{N \to \infty} {\frac{1}{N} \log \mathbb{P}_{\eta^N}[\pi^N \in \mathcal{C}]} &\leq - \inf_{\pi \in \mathcal{C}} I_{[0,T]}(\pi| \gamma),\\
    \liminf_{N \to \infty} {\frac{1}{N} \log \mathbb{P}_{\eta^N}[\pi^N \in \mathcal{O}]} &\geq - \inf_{\pi \in \mathcal{O}} I_{[0,T]}(\pi| \gamma).
\end{align}
\end{theorem}

 {\begin{remark}
We observe that the assumption that, for every $k \in \{1,\dots,l\}$, $B_k(\cdot)$ and $D_{-k}(\cdot)$ defined in \eqref{B_D} are concave functions is fulfilled not only for the case of the choice boundary rates $R^\pm$ of \cite{Nahum2020} but also for the ones in Section \ref{model1_5}. We will prove this last statement there.
\end{remark}}

\subsection{Stationary Measure: Bulk Vs Boundary Dynamics.} \label{bulk_vs_boundary} 

Since this model is described by an irreducible continuous-time Markov chain with finite state space, then it has a unique stationary measure.

It is known that if $l = 1$ and  {the density of the boundary reservoirs are the same and equal to $\rho \in (0,1)$}, then the unique stationary measure is reversible and given by the Bernoulli product measure with density $\rho$,  {denoted by $\mu_\rho$:}
\begin{equation} \label{Bernoulli_product}
    \mu_\rho(\eta) = \prod_{x=\epsilon_N}^{\tau_N} \rho^{\eta(x)} (1-\rho)^{1-\eta(x)}.
\end{equation} See for example \cite{Nahum2020}.
However, as soon as  {those  {densities}  {are} different}, it is well known that this difference creates a current in the system, driving it out of equilibrium. The system is no longer reversible and, in particular, exhibits long-range correlations \cite{spohn1983long}. Moreover, when $l\geq 2$, particles are created and annihilated at more than one site in the bulk of $\Omega_N$, and
according to different rates which depend on the environment. For this reason, the boundary dynamics can not be represented by a one-site dynamics. The state space is here $\{0,1\}^{ {\Sigma_l^\pm}}$ instead of $\{0,1\}$ as in the exclusion dynamics of \cite{landim2021steady}, \cite{Goncalves19Hydro},\cite{FGLN23} and \cite{bertini2009dynamical}.

Consider a neighborhood  {$\{(N-k)\epsilon_N, \dots, \tau_N\}$} of the right boundary  {reservoir}. The stationary state of the bulk dynamics restricted to this set (we forbid the exchange of particles between $ {(N-k-1)\epsilon_N}$ and $ {(N-k)\epsilon_N})$ is the uniform measure over all configurations with a fixed number of particles. For $k$ large, by the equivalence of ensembles, locally this measure is close to a Bernoulli product measure with some fixed density.  {Only} in very special cases, the stationary state on  {$\Sigma^\pm_l$} induced by the generator  {$\mathscr{L}_N^{\pm}$} introduced above is a Bernoulli product measure. When this does not happen, there is a conflict between the bulk dynamics, which drives the system towards a Bernoulli product measure, and the  {boundary} dynamics, which propels the system to another stationary state. As the bulk dynamics is accelerated by $N^2$, while the boundary dynamics is sped up by $N$, the bulk dynamics  {prevails} and the state of the system at the boundary is close to a Bernoulli product measure. In particular,  {the Bernoulli product measures can be used as the reference measure in order to prove the hydrodynamic limit}. However, for models whose boundary and bulk dynamics  {are accelerated} on the same time scale, and such that the boundary dynamics is not reversible with respect to a product measure,  {one cannot use the same approach with the same reference measures,} as mentioned in \cite{erignoux2018hydrodynamic}.\\

\section{Hydrodynamic Limit.} \label{section_hydrodynamics}
Here we follow the usual strategy: prove tightness of the sequence of probability measures $(\mathbb{Q}_{\nu_N})_{N \in \mathbb{N}}$ which guarantees the weak convergence up to a subsequence to a limit point and characterize uniquely the latter.

\subsection{Tightness.} \label{tightness_section} We denote by $\Delta_N$ the discrete Laplace operator which is defined, for every $G: [0,1] \to \mathbb{R}$ and $x \in \Lambda_N$, as
 {\begin{equation*}
    \Delta_N G(x) := N^2[G(x+\epsilon_N) + G(x-\epsilon_N)- 2 G(x)].
\end{equation*}}
We also denote by $\nabla^+_N$ the discrete gradient operator which is given, for every $G: [0,1] \to \mathbb{R}$ and $x \in \Lambda_N$, by
 {\begin{equation*}
    \nabla^+_N G(x) := N[G(x+\epsilon_N) - G(x)].
\end{equation*}}

After some simple but long computations, we have that,  {for every $s \in [0,T]$,}
\begin{equation} \label{generator_on_eta}
\mathscr{L}_N\langle\pi_s^N,H\rangle=\langle h_{-}(\eta_s), H\rangle+\langle h_{+}(\eta_s), H\rangle + \langle \pi_s^N, \Delta_N H\rangle+\eta_s(\epsilon_N)\nabla^{+}_N H(0)-\eta_s(\tau_N)\nabla^{+}_N H(\tau_N)
\end{equation}
and
\begin{align*}
&\mathscr{L}_N\langle\pi_s^N,H\rangle^2-2\langle\pi_s^N,H\rangle\mathscr{L}_N\langle\pi_s^N,H\rangle \\
=& {\sum_{x\in\Lambda_N^\circ}[\eta_s(x+\epsilon_N)-\eta_s(x)]^2\left(\nabla^+_N H(x)\right)^2}+\frac{1}{N}\sum_{\xi\in {\{0,1\}^{\Sigma^-_l}}} R^{-}(\Pi_{1,l}\eta,\xi) A( {\Sigma^-_l},H,\xi,\eta)\\
&+\frac{1}{N}\sum_{\xi\in {\{0,1\}^{\Sigma^+_l}}} R^{+}(\Pi_{N-l,N-1}\eta,\xi)A( {\Sigma_l^+},H,\xi,\eta),
\end{align*}
where we define
\begin{align} \label{def_langle_h_minus} 
\langle h_{-}(\eta_s), H\rangle := \sum_{j \in  {\Sigma^-_l}}  {(1-2\eta_s(j ))H(j)}\sum_{\substack{\xi\in {\{0,1\}^{\Sigma^-_l}} \\
    \xi(j )\neq\eta(j)}} R^{-}(\Pi_{1,l}\eta_s,\xi),
\end{align}
\begin{align} \label{def_langle_h_plus}
    \langle h_{+}(\eta_s), H\rangle := \sum_{j \in  {\Sigma^+_l}} {(1-2\eta_s(j ))H(j)} \sum_{\substack{\xi\in {\{0,1\}^{\Sigma^+_l}} \\
    \xi(j )\neq\eta(j)}} R^{+}(\Pi_{N-l,N-1}\eta_s,\xi).
\end{align} and, for every $S \subset \Lambda_N$,
\begin{align}
     A(S,H,\xi,\eta) &:= \left[\sum_{y\in S}(\xi(y)-\eta(y))H(y)\right]^2.
\end{align}
Tightness for  {$(\mathbb{Q}_{\nu_N})_{N \in \mathbb{N}}$}
follows from the fact that 
    \begin{equation} \label{dynkin1}
    \mathcal{M}^H_N(t) = \langle \pi^N_{t}, H  \rangle - \langle \pi^N_{0}, H  \rangle - \int_0^t \mathscr{L}_N\langle \pi^N_s, H\rangle ds
\end{equation} and 
    \begin{equation*}
    \mathcal{N}^H_N(t) = (\mathcal{M}^H_N(t))^2 - \int_0^{t} \{\mathscr{L}_N[\langle \pi^N_s, H\rangle^2] - 2 \langle \pi^N_s, H\rangle\mathscr{L}_N\langle \pi^N_s, H\rangle \}ds
\end{equation*}
are zero mean martingales, also that, for every $H \in C^2([0,1])$, $|\mathscr{L}_N\langle \pi^N_t, H\rangle|$ is uniformly bounded in $N$ and $t$ and finally that  {$\displaystyle\lim_{N \to \infty} \mathbb{E}_{\nu_N} [\mathcal{M}^H_N(t)^2 ] = 0$} uniformly in $t$. For more details, see for example \cite{FGS23}.

\subsection{Characterization of the limit points.} 
From Prohorov's theorem, we know that tightness implies that there exists a subsequence  {$(\mathbb{Q}_{\nu_{N_k}})_{ {k \geq 1}}$ of $(\mathbb{Q}_{\nu_N})_{ {N \geq 1}}$} and a probability measure $\mathbb{Q}$ over $ {\mathcal{D}}([0,T],\mcb{M})$ such that, as $k \to \infty$,  {$\mathbb{Q}_{\nu_{N_k}}$} converge weakly to $\mathbb{Q}$. Let us now characterize $\mathbb{Q}$  {to} show that it is independent of the choice of subsequence,  {from which we conclude the weak convergence of the whole sequence $(\mathbb{Q}_{\nu_N})_{N \geq 1}$.}

\begin{proposition} \label{proposition_characterization_limit_point_hidro}
The limit point $\mathbb{Q}$ is the unique probability measure that has support over the absolutely continuous measure $\pi_\cdot = \rho_\cdot  dx$ where $\rho:  {\Omega_T} \to \mathbb{R}$ is the unique weak solution of \eqref{PDEproblem}.
\end{proposition}

First, by Riesz's Representation Theorem  {and the fact that the occupation variables are bounded}, $\mathbb{Q}$ has support over absolutely continuous measures  {- for details, see a similar proof in page 60 of Chapter 4 of \cite{KL99}.} Now we show that $\mathbb{Q}$ gives measure one to the set of all absolutely continuous measures with density $ {\rho_\cdot}$ that  {are} weak solutions of \eqref{PDEproblem}. 

 {We start by showing the regularity of the density profiles $\rho$. First we have that $\rho \in L^2([0,T];L^2([0,1]))$ since we can show that $0 \leq \rho(t,x) \leq 1$ for all $(t,x) \in [0,T] \times [0,1]$ as a consequence of the boundedness of the occupation variables - see Lemma \ref{rho_in_l2_l2}.} Then, by the energy estimate given in Lemma \ref{energy_estimate_lemma_statement}, we conclude that $\rho \in L^2([0,T];H^1([0,1]))$ - see Proposition \ref{proposition_rho_has_weak_deriv}.

\begin{lemma} \label{rho_in_l2_l2}
    All limit points $\mathbb{Q}$ of the sequence $( \mathbb{Q}_{\nu_N})_{ {N \geq 1}}$ are concentrated on one path of absolutely continuous measures whose density, which we denote by $\rho_t$, is such that
    \begin{equation*}
        \int_0^T \int_0^1 [\rho_t(u)]^2 du ds < \infty \quad \textrm{ $\mathbb{Q}$-a.s. }
    \end{equation*}
\end{lemma}

\begin{proposition} \label{proposition_rho_has_weak_deriv}
All limit points $\mathbb{Q}$ of the sequence $(\mathbb{Q}_{\nu_N})_{N \geq 1}$ are such that, $\mathbb{Q}$-a.s. there exists  {a} $L^2([0,T],H^1[0,1])$ function, called the weak derivative of $\rho$ and that we denote by $\nabla \rho_s(u)$, such that,  {for every $G$ smooth function, it holds that
\begin{align*}
    \int_0^T \left[\int_0^1 \rho_s(x) \nabla G_s(x) dx - F^+(\rho_t(1)) G(1) - F^-(\rho_t(0)) G(0) \right] ds  = - \int_0^T \int_0^1 \nabla \rho_s(x) G_s(x) dx ds
\end{align*}} and
\begin{align*}
    \int_0^T \int_0^1 \frac{|\nabla \rho_s(u)|^2}{\rho_s(u)[1-\rho_s(u)]} du ds < \infty.
\end{align*}
\end{proposition}

The proof of  {the previous proposition} relies on the following energy estimate:
\begin{lemma}[Energy estimate] \label{energy_estimate_lemma_statement}
    Let $\sigma$ represent the mobility, i.e. $\sigma(\alpha):=\alpha(1-\alpha)$. Then
    \begin{align} \label{energy_estimate_equation}
        \mathbb{E}_{\mathbb{Q}} \left[ \sup_{H \in C([0,T];C^1([0,1]))} \left\{ \int_0^T \int_0^1 \nabla H_s(x) \rho_s(x) dxds - 2 \int_0^T \int_0^1 [H_s(x)]^2 \sigma(\rho_s(x))dxds\right\}\right]  { < \infty.}
    \end{align}
\end{lemma}
The proof of Lemma \ref{energy_estimate_lemma_statement} relies on the Entropy Inequality and the bounds on the Dirichlet forms presented in Appendix \ref{appA}.
Proposition \ref{proposition_rho_has_weak_deriv} follows from this statement.

 {Since we know that all limit points $\mathbb{Q}$ give total measure for a trajectory of absolutely continuous measures whose density is in $L^2([0,T], H^1([0,1]))$, from here on, we can assume that all the density functions we will use next are in the Sobolev space $L^2([0,T];H^1([0,1]))$. For each absolutely continuous measure $\pi_t(dx) = \rho_t(x)dx$ (and so $\rho \in L^2([0,T], H^1([0,1]))$),} define, for every  {$H \in C^{1,2}( {\Omega_T})$},
 {\begin{align*}
     M^H(t) &:= \langle \rho_t, H_t \rangle - \langle \rho_0, H_0 \rangle - \int_0^t \langle \rho_s, \partial_t H_s + \Delta H_s \rangle ds \\
     &\int_0^t[\rho_s(0) \partial_u H_s(0) - \rho_s(1) \partial_u H_s(1) - F_+[ \rho_s(1)] H_s(1) -  F_-[ \rho_s(0)] H_s(0)] ds,
\end{align*}}where $\Delta$ denotes the Laplace operator. Our goal is to show that, for every $\delta > 0$ and for all  {$H \in C^{1,2}( {\Omega_T})$}, $\mathbb{Q}\left(A^\delta\right) = 0$, 
where  {\begin{equation} \label{set_A_epsilon}
    A^\delta := \Big\{\pi \in  {\mathcal{D}}([0,T], \mcb{M}) \ | \ \pi_t(du) = \rho_t(u) du,  {\ \rho \in L^2([0,1],H^1([0,1]))} \ \textrm{ and } \ \sup_{0\leq t\leq T}| M^H(t)| > \delta\Big\}.
\end{equation}}To do that, we would like to work with the measures  {$\mathbb{Q}_{\nu_{N_k}}$} and use Portmanteau's theorem to prove such a result. This can not be done directly since the set $A^\delta$ is not necessarily an open set in the Skorohod topology. To fix this problem, as in Section 5.2 of \cite{FGNslow}, let us define,  {for every $t \geq 0$ and $\epsilon \in (0,1)$,
\begin{equation*}
    \pi_t \ast \iota_\epsilon(y) := \begin{cases}\frac{1}{\epsilon}\int_0^\epsilon\pi_t(du), \textrm{ if } y \in [0,\epsilon_N], \\
    \frac{1}{\epsilon}\int_{1-\epsilon}^1\pi_t(du), \textrm{ if } y \in [\tau_N,1].
    \end{cases}
\end{equation*}}
Then  {setting}, for every $\epsilon > 0$ and every  {$H \in C^{1,2}( {\Omega_T})$,
\begin{align*}
     M_\epsilon^H(t) &:= \langle \pi_t, H_t \rangle - \langle \pi_0, H_0 \rangle - \int_0^t [\langle \pi_s, \partial_t H_s + \Delta H_s \rangle ]ds \\
    &- \int_0^t \left\{ F_+[ \pi_s \ast \iota_\epsilon(1)] H_s(1) + F_-[ \pi_s \ast \iota_\epsilon(0)] H_s(0) \right\} ds \\
    &+ \int_0^t [\pi_s \ast \iota_\epsilon(0) \partial_u H_s(0) - \pi_s \ast \iota_\epsilon(1)\partial_u H_s(1) ]ds,
\end{align*}} we have that
 {$\pi_t(du) \mapsto \sup_{0 \leq t \leq T} |M_\epsilon^H(t)|$} is continuous and so, for every $\delta > 0$, the set  {$$A_\epsilon^\delta :=\left\{ \pi \in  {\mathcal{D}}([0,T], \mcb{M}) \ \Big|  {\ \pi_t(du) = \rho_t(u) du,  {\ \rho \in L^2([0,1],H^1([0,1]))} \ \textrm{ and }} \ \displaystyle\sup_{0 \leq t \leq T} |M_\epsilon^H(t)| > \delta\right\}$$} is open in $ {\mathcal{D}}([0,T], \mcb{M})$. Then, for every $\delta > 0$,
 {\begin{align} \label{RHS}
    \mathbb{Q}\left(A^\delta\right) = \lim_{\epsilon \to 0} \mathbb{Q}\left(A_\epsilon^\delta\right)
    &\leq \lim_{\epsilon \to 0} \liminf_{k \to \infty} \mathbb{Q}_{\nu_{N_k}} \left(  {A_\epsilon^\delta}\right),
\end{align}}where the last inequality holds by Portmanteau's theorem.  {Define} $w_{\epsilon_{N_k}} := 0$, $w_{\tau_{N_k}} := 1$ and,  {for every $H \in C^{1,2}( {\Omega_T})$, $t \geq 0$ and $N \in \mathbb{N}$,
\begin{align*}
    M_\epsilon^{H, N}(t) &:= \langle \pi^N_t, H_t \rangle - \langle \pi^N_0, H_0 \rangle - \int_0^t [\langle \pi^N_s, \partial_t H_s + \Delta_N H_s \rangle ]ds \\
    &- \int_0^t \left\{ F_+[ \pi^N_s \ast \iota_\epsilon(\tau_N)] H_s(1) + F_-[ \pi^N_s \ast \iota_\epsilon(\epsilon_N)] H_s(0) \right\} ds \\
    &+ \int_0^t [\pi^N_s \ast \iota_\epsilon(\epsilon_N) \partial_u H_s(0) - \pi^N_s \ast \iota_\epsilon(\tau_N)\partial_u H_s(1) ]ds.
\end{align*}  {S}ince, for fixed $\epsilon>0$, taking $k$ sufficiently large such that $\epsilon_{N_k} < \epsilon$ and $\tau_{N_k} > 1 - \epsilon$, we have by definition of $\pi_t \ast \iota_\epsilon$ that $\pi_t \ast \iota_\epsilon(\tau_{N_k}) = \pi_t \ast \iota_\epsilon(1)$ and $\pi_t \ast \iota_\epsilon(\epsilon_{N_k}) = \pi_t \ast \iota_\epsilon(0)$, we can bound from above the RHS of the inequality in \eqref{RHS} by
\begin{align} \label{not_vanish}
    &\lim_{\epsilon \to 0} \liminf_{k \to \infty} \mathbb{P}_{\nu_{N_k}}\left(\sup_{0\leq t\leq T}| M_{\epsilon}^{H,N_k}(t)| > \delta\right).
\end{align}
Combining \eqref{generator_on_eta} with the definition of $\mathcal{M}^H_N(t)$ given in \eqref{dynkin1}, we have that
\begin{align*}
     \mathcal{M}^H_N(t) = \langle \pi^N_t, H_t \rangle &- \langle \pi^N_0, H_0 \rangle - \int_0^t [\langle \pi^N_s, \partial_t H_s + \Delta_N H_s \rangle ]ds \\
    &- \int_0^t [\langle h_{-}(\eta_s), H_s\rangle+ \langle h_{+}(\eta_s), H_s\rangle + \eta_s(\epsilon_N)\nabla^{+}_N H_s(0) - \eta_s(\tau_N)\nabla^{+}_N H_s(\tau_N) ]ds.
\end{align*} and so} $\mathbb{P}_{\nu_{N_k}}\left( \mathcal{M}^H_{N_k}(t) = 0, \ \forall t \in [0,T] \textrm{ and } \forall H \in C^{1,2}( {\Omega_T})\right) = 1$. Thus, for all $H \in C^{1,2}( {\Omega_T})$ and $\delta>0$,
\begin{equation} \label{to_use}
    \mathbb{P}_{\nu_{N_k}}\left(\sup_{0 \leq t \leq T} |\mathcal{M}^H_{N_k}(t)| > \delta \right) = 0.
\end{equation} So, we would like to relate $$\mathbb{P}_{\nu_{N_k}}\left(\sup_{0 \leq t \leq T} |M_\epsilon^{H,N_k}(t)| > \delta \right) \quad \textrm{ and } \quad \mathbb{P}_{\nu_{N_k}}\left(\sup_{0 \leq t \leq T} |\mathcal{M}^H_{N_k}(t)| > \delta \right)$$ to conclude our proof. This is the idea behind the following results, the so-called Replacement Lemmas.

\begin{lemma}[Basic Replacement Lemma]\label{repl_lemma_ave}
For every $H \in C([0,T])$, for any $t \in [0,T]$ and $z=\epsilon_N$, it holds
\begin{equation} \label{eq_1_RLbasic}
\lim_{\epsilon\to 0}\lim_{N\to+\infty}\mathbb{E}_{\nu_N} \left[\Big | \int_0^t [\eta_{sN^2}(z) - \overrightarrow{\eta}^{\lfloor \epsilon N\rfloor}_{sN^2}(z)] H(s)ds\Big| \right] =0,
\end{equation}
where, for every $L\in\mathbb N$ and $\eta \in \Omega_N$,
\begin{equation} \label{right_average}
\overrightarrow{\eta}^L(z) :=\frac{1}{L\epsilon_N}\sum_{y= z+\epsilon_N}^{z+L\epsilon_N}\eta(y).
\end{equation}
The same result holds for $z = \tau_N$ by replacing the average to the right, i.e. $\overrightarrow{\eta}^{\lfloor \epsilon N\rfloor}_{sN^2}$, by an average to the left, i.e. $\overleftarrow{\eta}^{\lfloor \epsilon N\rfloor}_{sN^2}$, which is defined, for every $L\in\mathbb N$ and $\eta \in \Omega_N$, as
\begin{equation}\label{left_average}
 \overleftarrow{\eta}^L(z) :=\frac{1}{L\epsilon_N}\sum_{y=z-L\epsilon_N}^{z-\epsilon_N}\eta(y).
\end{equation}
\end{lemma}

Remark that $\pi^N_t \ast \iota_\epsilon(\epsilon_N) = \overrightarrow{\eta}^{\lfloor \epsilon N \rfloor}(\epsilon_N)$ and $\pi^N_t \ast \iota_\epsilon(\tau_N) = \overleftarrow{\eta}^{\lfloor \epsilon N \rfloor}(\tau_N)$.

\begin{lemma}[Boundary Replacement Lemma]\label{repl_lemma_2}
For any function $H \in C^{1,2}( {\Omega_T})$ and any $t \in [0,T]$, it holds, for any $N \in \mathbb{N}$ and any $\epsilon > 0$,
 {
\begin{align*}
\mathbb{E}_{\nu_N} &\left[\Big | \int_0^t [\langle h_{-}(\eta_s), H_s\rangle - F_-(\overrightarrow{\eta}_{sN^2}^{\lfloor \epsilon N \rfloor}(\epsilon_N)) H_s(0)] ds \Big | \right] \\
&\hspace{4.5cm} \lesssim  {\frac{l}{\epsilon N}} +\epsilon ( l B + 1 + l) + \sum_{j = 2}^{l} \sum_{p=1}^{j-1} \binom{j-1}{p} \frac{2 l}{\epsilon^{p-1} N ^p} {+ \frac{1}{B}},
\end{align*}} for all $B > 0$ independent of $N$ and $\epsilon$. The result also holds by replacing $h_-, F_-$ and $\overrightarrow{\eta}$ by $h_+$, $F^+$ and $\overleftarrow{\eta}$, respectively.
Above $\overrightarrow{\eta}^{\lfloor \epsilon N\rfloor}_{sN^2}$ and $\overleftarrow{\eta}^{\lfloor \epsilon N\rfloor}_{sN^2}$ are the same as in \eqref{right_average} and \eqref{left_average} {, respectively, and $F_{\pm}(\cdot)$, $\langle h_-(\eta_s), H_s\rangle$ and $\langle h_+(\eta_s), H_s\rangle$ are as in \eqref{def_F_pm}, \eqref{def_langle_h_minus} and \eqref{def_langle_h_plus}, respectively.}
\end{lemma} The proof of Lemma \ref{repl_lemma_ave} can be adapted from the ideas of Lemma 7 of \cite{Goncalves19Hydro} and for that reason we omit it here. For the proof of Lemma \ref{repl_lemma_2} see Section \ref{replacement_lemmas_proof_section}. From these two lemmas, we can easily conclude the following results.

\begin{corollary} \label{corollary_repl_lemma_ave}
For any function $H \in C^{1,2}( {\Omega_T})$and for every $\delta > 0$, it holds that
\begin{equation} \label{eq_1_corollary_RL1}
\lim_{\epsilon\to 0}\lim_{N\to+\infty}\mathbb{P}_{\nu_N} \left( \sup_{0 \leq t \leq T} \Big | \int_0^t \eta_{sN^2}(\epsilon_N)\nabla^+_N H_s(0) - \pi^N_s \ast \iota_\epsilon(\epsilon_N)\partial_u H_s(0) ds\Big| > \delta \right) =0,
\end{equation} and
\begin{equation} \label{eq_2_corollary_RL1}
\lim_{\epsilon\to 0}\lim_{N\to+\infty}\mathbb{P}_{\nu_N} \left( \sup_{0 \leq t \leq T} \Big | \int_0^t \eta_{sN^2}(\tau_N)\nabla^+_N H_s(\tau_N) - \pi^N_s \ast \iota_\epsilon(\tau_N)\partial_u H_s(1) ds\Big| > \delta \right) =0.
\end{equation}
\end{corollary}

\begin{corollary} \label{corollary_repl_lemma_2}
For any function $H \in C^{1,2}( {\Omega_T})$ and for every $\delta >0$, it holds that
\begin{equation} \label{eq_1_corollary_RL2}
\lim_{\epsilon \to 0} \lim_{N\to+\infty}  \mathbb{P}_{\nu_N} \left( \sup_{0 \leq t \leq T}\Big | \int_0^t [\langle h_{-}(\eta_s), H_s\rangle - F_-(\pi^N_s \ast \iota_\epsilon(\epsilon_N)) H_s(0)] ds \Big | > \delta \right) =0,
\end{equation}
and
\begin{equation} \label{eq_2_corollary_RL2}
\lim_{\epsilon \to 0} \lim_{N\to+\infty}  \mathbb{P}_{\nu_N} \left( \sup_{0 \leq t \leq T} \Big | \int_0^t [\langle h_{+}(\eta_s), H_s\rangle - F_+(\pi^N_s \ast \iota_\epsilon(\tau_N)) H_s(1)] ds \Big | > \delta \right) =0.
\end{equation}
\end{corollary} The proof of Corollary \ref{corollary_repl_lemma_ave} and \ref{corollary_repl_lemma_2} can be found in Section \ref{replacement_lemmas_proof_section} after the proof of Lemma \ref{repl_lemma_2}.

By the previous results, we can now conclude that, for any $H \in C^{1,2}( {\Omega_T})$ and any $\delta >0$,
\begin{align} \nonumber
    &\lim_{\epsilon \to 0} \limsup_{k \to \infty} \mathbb{P}_{\nu_{N_k}}\left(\sup_{0\leq t\leq T}| M_{\epsilon}^{H,N_k}(t)| > \delta \right)\\ \label{to_bound1}
    &\leq \limsup_{k \to \infty} \mathbb{P}_{\nu_{N_k}} \left(\sup_{0 \leq t \leq T} \Big| \mathcal{M}^H_{N_k}(t)\Big| > \frac{\delta}{5} \right) \\ \label{to_bound2}
    &+ \lim_{\epsilon\to 0} \limsup_{k \to \infty} \mathbb{P}_{\nu_{N_k}} \left(\sup_{0 \leq t \leq T} \Big | \int_0^t [\langle h_{-}(\eta_s), H_s\rangle - F_-(\pi^{N_k}_s \ast \iota_\epsilon(\epsilon_{N_k})) H_s(0)] ds \Big | > \frac{\delta}{5} \right) \\ \label{to_bound3}
    &+ \lim_{\epsilon\to 0} \limsup_{k \to \infty} \mathbb{P}_{\nu_{N_k}} \left(\sup_{0 \leq t \leq T} 
    \Big|\int_0^t [\langle h_{+}(\eta_s), H_s\rangle - F_+(\pi^{N_k}_s \ast \iota_\epsilon(\tau_{N_k})) H_s(1)] ds\Big| > \frac{\delta}{5}\right) \\ \label{to_bound4}
    &+ \lim_{\epsilon\to 0} \limsup_{k \to \infty} \mathbb{P}_{\nu_{N_k}} \left( \sup_{0 \leq t \leq T} \Big | \int_0^t \eta_{s N^2_k}(\epsilon_{N_k})\nabla^+_{N_k} H_s(0) - \pi^{N_k}_s \ast \iota_\epsilon(\epsilon_{N_k})\partial_u H_s(0) ds\Big| > \frac{\delta}{5} \right)\\ \label{to_bound5}
    &+ \lim_{\epsilon\to 0} \limsup_{k \to \infty} \mathbb{P}_{\nu_{N_k}} \left(\sup_{0 \leq t \leq T} 
    \Big | \int_0^t \eta_{s N^2_k}(\tau_{N_k})\nabla^+_{N_k} H_s(\tau_{N_k}) - \pi^{N_k}_s \ast \iota_\epsilon(\tau_{N_k}) \partial_u H_s(1) ds\Big| > \frac{\delta}{5} \right).
\end{align}
By \eqref{to_use}, the RHS of \eqref{to_bound1} is zero; by Corollary \ref{corollary_repl_lemma_2},  \eqref{to_bound2} and \eqref{to_bound3} are zero; and, by Corollary \ref{corollary_repl_lemma_ave}, \eqref{to_bound4} and \eqref{to_bound5} are zero. This finishes our proof.\par

\subsection{Uniqueness} We start by remarking that $F_+, F_-:[0,1] \to \mathbb{R}$ are  {$C^1([0,1])$} functions (not necessarily linear) and that $ {\gamma}$  {is} a measurable function. Taking  {$\alpha=1$} and $\beta = 0$ in Theorem \ref{th_uniqueness}, we assure the uniqueness of a weak solution $\rho$ of the initial value problem \eqref{PDEproblem} that is in $L^2([0,T];H^1([0,1]))$. This proves that the limiting measure $\mathbb{Q}$ is completely characterized, and so, is unique. From here, the desired weak convergence of the sequence of probability measures $(\mathbb{Q}_{\nu_N})_{N \in \mathbb{N}}$ follows and so the Hydrodynamic Limit, i.e. Theorem \ref{th_hydrodynamics}.

\section{Large Deviations Principle} \label{large_dev_proof}
In this section we prove the dynamical large deviations principle for the empirical measure of the exclusion process considered here. We first present some classical properties of the rate function. Then, we construct the decomposition of the rate function, a technique introduced in \cite{FGLN23}, that allow us to take care of the boundary terms in the proof of the lower bound. After that we prove the upper and lower bounds in Theorem \ref{main_result_large_dev}, finishing with the proof of  {Theorem \ref{thm_I_density}}, the most technical step in the proof of the lower bound.

\subsection{The Rate Function $I_{[0,T]}(\cdot)$} \label{I-density}

In this subsection, we state several properties of the rate functional $I_{[0,T]}(\cdot|\gamma)$. We omit their proofs since they are very close to the ones presented in \cite{FGLN23}. Let us start with two elementary bounds.

Let $\tau_r \rho:\mathbb{R}_+\times [0,1]\to\mathbb{R}$ be the function given by $ {(\tau_r\rho)}_t(x)=\rho(t+r,x)$. The first bound asserts that, for any $0< S< T$, the cost of a trajectory in the time interval $[0,T]$ is bounded by its cost on $[0,S]$ plus the cost on the interval $[S,T]$ while the second one states that the cost of a trajectory in any subinterval of $[0,T]$ is bounded by the cost on the whole interval.
\begin{proposition} For all $\pi(t,dx)=\rho_t(x)dx$ in $ {\mathcal{D}}([0,T],\mathscr{M}_{ac})$ and $0<S<T$:
  \begin{equation}\label{cost_subinterval}
    I_{[0,T]}(\rho)\leq I_{[0,S]}(\rho)+I_{[0,T-S]}(\tau_S\rho)\quad\text{and}\quad I_{[0,S]}(\rho)\leq I_{[0,T]}(\rho)
\end{equation}  
\end{proposition}  {We remark that one can in fact show that equality holds in both cases of \eqref{cost_subinterval} but for our purposes this result is enough.}

The next result shows that a path with finite $I_{[0,T]}(\cdot|\gamma)$ rate function has a density whose initial condition is $\gamma$ and is continuous. The result is similar to Lemma 3.1 in \cite{FGLN23}. We present here for completeness.
\begin{lemma}
    Fix $T>0$ and $\gamma\in\mathscr{M}_{ac}$. Let $\rho\in  {\mathcal{D}}([0,T], \mathscr{M}_{ac})$ be such that $I_{[0,T]}(\rho|\gamma)<\infty.$ Then, $\rho(0,x)=\gamma(x).$ Moreover, for each $M>0$, $G\in C^2([0,1])$ and $\varepsilon>0$, there exists a $\delta>0$ such that 
    \begin{equation*}
        \sup_{\rho:I_T(\rho|\gamma)\leq M}\sup_{|t-s|\leq\delta}|\langle \rho_t, G\rangle-\langle \rho_s, G\rangle|\leq \varepsilon.
    \end{equation*}
    In particular, $\rho$ belongs to $C([0,T],\mathscr{M}_{ac})$.
\end{lemma}
The next proposition plays an important role in the proof of Theorem \ref{thm_I_density}.  {See \cite{FGLN23} for the proof strategy.}

\begin{proposition}
    There exists a constant $C_0>0$ such that 
    \begin{equation}
        \int_0^T\int_0^1 \frac{|\nabla\rho_t(x)|^2}{\sigma(\rho_t(x))}\;dx\;dt\leq C_0\{I_{[0,T]}(\rho)+1\}
    \end{equation}
    for any path $u\in  {\mathcal{D}}_{\mathcal{E}}([0,T],\mcb{M}_{ac})$.
\end{proposition}
\begin{corollary}\label{rate_function_0}
    The density $\rho$ of a path $\pi\in  {\mathcal{D}([0,T],\mcb{M}_{ac})}$ is the weak solution of the initial-boundary value problem \eqref{PDEproblem} if, and only if, $I_{[0,T]}(\rho|\gamma)=0.$
\end{corollary}
\begin{theorem}\label{rate_function_lsc_compac_ls}
    Fix $T>0$ and a measurable function $\gamma:[0,1]\to [0,1]$. Assume that,  {for every $k \in \{1,\dots,l\}$, $B_k(\cdot)$ and $D_{-k}(\cdot)$ defined in \eqref{B_D} are concave functions}. Then, the function $I_{[0,T]}(\cdot|\gamma): {\mathcal{D}}([0,T],\mcb{M})\to [0,+\infty]$ is  {convex}, lower-semicontinuous and has compact level sets.
\end{theorem}
\subsection{Deconstructing the Rate Function}\label{deconstruction_sec}
In this section we prove that the rate function $I_{[0,T]}(\cdot)$ can be decomposed as the sum of two functionals $I_{[0,T]}^{(1)}(\cdot)+I_{[0,T]}^{(2)}(\cdot)$, where the first one measures the cost of a trajectory in the bulk and the other on the boundary. This is the main tool in the proof of Theorem \ref{thm_I_density}. Before presenting the decomposition, we need to introduce the weighted Sobolev spaces.
\subsubsection{Weighted Sobolev Spaces.} For a non-negative function  {$f:\Omega_T\to\mathbb{R}_{+}$}, denote by $L^2(f)$ the Hilbert space induced by the smooth functions in $C^\infty(\Omega_T)$ endowed with the scalar product defined by
\begin{equation*}
    \langle\langle G,H\rangle\rangle_f=\int_0^T\int_0^1 f_t\; G_t\; H_t\;dx\;dt.
\end{equation*}
Denote by $C_c^\infty(\Omega_T)$ the space of smooth real functions $H:\Omega_T\to\mathbb{R}$ with support in $(0,T)\times(0,1)$. Let $H_0^1(f)$ be the Hilbert spaces induced by the set $C^\infty_c(\Omega_T)$ endowed with the scalar product $\langle\langle G,H \rangle\rangle_{1,f}$, defined by 
\begin{equation*}
    \langle\langle G,H \rangle\rangle_{1,f}=\langle\langle \nabla G,\nabla H \rangle\rangle_{f}.
\end{equation*}
Let $||\cdot||_{1,f}$ denote the associated norm. We define $H^{-1}(f)$ to be the dual space of $H_0^1(f)$. We say a linear functional $L:H_0^1(f)\to\mathbb{R}$ belongs to $H^{-1}(f)$ if 
\begin{equation*}
    ||L||_{-1,f}:=\sup_{G\in C_c^\infty(\Omega_T)}\{2L(G)-||G||^2_{1,f}\}<+\infty.
\end{equation*}
The space $H^{-1}(f)$ can be identified with the space $\{\nabla P\;:\; P\in L^2(1/f)\}$. This is the content of the next lemma. It  {corresponds to} Lemma 4.8 in \cite{bertini2009dynamical}, and therefore the proof is omitted here. We will use now and throughout the text the notation $\langle f\rangle:=\int_0^1 f(x)\;dx.$
\begin{lemma}\label{structureH-1}
    A linear functional $L:H_0^1(f)\to\mathbb{R}$ belongs to $H^{-1}(f)$ if, and only if, there exists $P\in L^2(1/f)$ such that 
    \begin{equation*}
        L(H)=\int_0^T\int_0^1 P_t\;\nabla H_t\;dx\;dt,
    \end{equation*}
    for every $H\in C_{c}^{\infty}(\Omega_T)$. In this case, 
    \begin{equation*}
        ||L||^2_{-1,f}=\int_0^T\{\langle P_t, P_t\rangle_{f(t)^{-1}}-c_t\}\;dt,
    \end{equation*}
    where $c_t=\{\langle P_t/f_t\rangle^2\;/\;\langle 1/f(t)\rangle\}\mathds{1}\{\langle 1/f(t)\rangle<\infty\}$.
\end{lemma}
\subsubsection{Decomposition} In this section we present a decomposition of the rate functional which will be the main step in the proof of Theorem \ref{thm_I_density}. Until the end of this section, we will assume that $\pi(t,dx)=\rho_t(x)dx$ is a path in $ {\mathcal{D}}_{\mathcal{E}}([0,T],\mathscr{M}_{ac})$, whose density $\rho$ is continuous in $\Omega_T$, smooth in time, that there exists an $\varepsilon>0$ such that $\varepsilon\leq \rho_t(x)\leq 1-\varepsilon$ and that $I_{[0,T]}(\rho)<\infty.$ In order to do the decomposition, let us separate the linear and non-linear parts of the rate functional. Let $\mathfrak{M}:C^{0,1}(\Omega_T)\to\mathbb{R}$ be the functional given by 
\begin{equation*}
    \mathfrak{M}(H):=\int_0^T\int_0^1\sigma(\rho_t)|\nabla H_t|^2\;dx\;dt+\int_0^T\Psi(t,H_t(0),H_t(1))\;dt,
\end{equation*}
where $\Psi(t,H_t(0),H_t(1))$ denotes the boundary term
\begin{equation*}
    \Psi(t,H_t(0),H_t(1))=\mathfrak{b}^-(\rho_t(0),H_t(0))+\mathfrak{b}^+(\rho_t(1),H_t(1)).
\end{equation*}
Defining 
\begin{equation*}
\mathfrak{L}(H)=\int_0^T\langle\partial_t\rho_t , H_t\rangle \;dt-\int_0^T\langle\nabla \rho_t,\nabla H_t\rangle\;dt,
\end{equation*}
we can write
\begin{equation*}
    I_{[0,T]}(\rho)=\sup_{H\in C^{1,2}(\Omega_T)}J_{T,H}(\rho)=\sup_{H\in C^{1,2}(\Omega_T)}\{\mathfrak{L}(H)-\mathfrak{M}(H)\}.
\end{equation*}
Now we decompose this supremum into two: one over functions that vanish at the boundary and the other on general $C^{1,2}$ functions that match the right values at the boundary. This is natural having in mind the case of Dirichlet boundary conditions, where the interaction with the boundary is so strong that we do not see any contribution of the boundary on the dynamical rate functional. 
Let $H:\Omega_T\to\mathbb{R}$ be any function in $C^{1,2}(\Omega_T)$. Suppose we have decomposed it into 
\begin{equation*}
    H=H^{(0)}+H^{(1)}
\end{equation*}
such that $H^{(0)}(t,0)=H^{(0)}(t,1)=0$, for all $0\leq t\leq T$, and in a way that $H^{(1)}$ depends on $H$ only through its values on the boundary. Then, since $\mathfrak{L}$ is linear, our rate  {functional} becomes 
\begin{equation*}
 \sup_{H\in C^{1,2}(\Omega_T)}\{\mathfrak{L}(H^{(0)})+\mathfrak{L}(H^{(1)})-\mathfrak{M}(H^{(0)}+H^{(1)})\}   
\end{equation*}
\begin{align*}
    =\sup_{H\in C^{1,2}(\Omega_T)}&\Big\{\mathfrak{L}(H^{(0)})-\int_0^T\int_0^1\sigma(\rho_t)|\nabla H^{(0)}_t|^2\;dx\;dt-2\int_0^T\int_0^1\sigma(\rho_t)\nabla H^{(0)}_t \nabla H^{(1)}_t\;dx\;dt\\
    & {+} \mathfrak{L}(H^{(1)})-\int_0^T\int_0^1\sigma(\rho_t)|\nabla H^{(1)}_t|^2\;dx\;dt-\int_0^T \Psi(t,H_t(0),H_t(1))\;dt\Big\}.
\end{align*}
The first line depends only on $H^{(0)}$ and a crossed term, while the second line depends only on the values of $H$ at the boundaries, namely, $H_t(0)$ and $H_t(1).$ Therefore, if we define in the decomposition $H^{(1)}$ so that we can kill the crossed term, we will have two independent variational problems.
\par
Let $\Xi:\Omega_T\to\mathbb{R}$ be the function given by
\begin{equation}\label{def_Xi}
    \Xi(t,x)=\frac{1}{\int_0^11/\sigma(\rho(t,y))\;dy}\int_0^x \frac{1}{\sigma(\rho(t,y) {)}}\;dy.
\end{equation}
Note that $\Xi\in C^{\infty,1}(\Omega_T)$ and that for all $0\leq t\leq T$, $\Xi(t,0)=0,\;\Xi(t,1)=1$. Given a function $H:\Omega_T\to\mathbb{R}$, we decompose it as $H=H^{(0)}+H^{(1)}$, 
where
\begin{equation*}
    H^{(1)}(t,x)= H(t,0) + [H(t,1)-H(t,0)]\;\Xi(t,x).
\end{equation*}
Note that $H^{(1)}(t,0)=H(t,0)$ and $H^{(1)}(t,1)=H(t,1)$. For this choice of the decomposition, one can easily check that $\nabla H^{(0)}$ and $\nabla H^{(1)}$ are orthogonal in $L^2(\sigma(\rho)):$
\begin{equation*}
    \int_0^T\int_0^1\sigma(\rho_t)\nabla H_t^{(0)}\nabla H_t^{(1)}\;dx\;dt=0.
\end{equation*}
Let $L:C^{0,1}(\Omega_T)\to\mathbb{R}$ denote any linear functional. We denote by $L_0$ its restriction to the subset $C_0^{0,1}(\Omega_T)$:
\begin{equation*}
    L_0(H)=L(H),\quad H\in C_0^{0,1}(\Omega_T),
\end{equation*}
where $C_0^{0,1}(\Omega_T)=\{H\in C^{0,1}(\Omega_T)\;:\;H_t(0)=H_t(1)=0,\;0\leq t\leq T\}$.
By the reasoning of the previous paragraphs, we have the following decomposition result that holds for any linear functional $L:C^{0,1}(\Omega_T)\to\mathbb{R}$.
\begin{lemma}
    Let $L:C^{0,1}(\Omega_T)\to\mathbb{R}$ be any linear functional and recall 
    \begin{equation*}
    \mathfrak{M}(H):=\int_0^T\int_0^1\sigma(\rho_t)|\nabla H_t|^2\;dx\;dt+\int_0^T\Psi(t,H_t(0),H_t(1))\;dt.
\end{equation*}
Then 
\begin{equation*}
    \sup_{C^{0,1}(\Omega_T)}\{L_0(H)-\mathfrak{M}(H)\}=S_1+S_2,
\end{equation*}
where
\begin{equation*}
    S_1=\sup_{C^{0,1}_0(\Omega_T)}\{L_0(G)-\int_0^T\int_0^1\sigma(\rho_t)|\nabla G_t|^2\;dx\;dt\}
\end{equation*}
and 
\begin{equation*}
    S_2=\sup_{h,g\in C([0,T])}\{L(g_t[1-\Xi_t])+L(h_t\Xi_t)-\int_0^T\zeta_t[h_t-g_t]^2\;dt-\int_0^T\Psi(t,g_t,h_t)\;dt\}.
\end{equation*}
In the above formula, $\zeta_t=1/\langle 1/\sigma(\rho_t) \rangle$.
\end{lemma}

This lemma was first proved in \cite{FGLN23}. We observe that the first variational problem is related to the interior of $\Omega_T$ while the second one is related to the boundary of the cylinder $\Omega_T$. Moreover, the only important fact about the boundary variational problem is that it depends on the functions $H$ only through its values on the boundary. Therefore, the proof is carried out in the same way as in \cite{FGLN23} and, for that reason, we omit the details here. We apply it now to the linear functionals appearing in the definition of the rate function. 

Let $L^{(\partial_t)}, L^{(\nabla)}:C^{0,1}(\Omega_T)\to\mathbb{R}$ denote the linear functionals 
\begin{equation*}
    L^{(\partial_t)}(H)=\int_0^T\langle \partial_t u_t, H_t\rangle\;dt,\quad L^{(\nabla)}(H)=\int_0^T\langle \nabla u_t,\nabla H_t\rangle\;dt,
\end{equation*}
so that $\mathfrak{L}=L^{(\partial_t)}+L^{(\nabla)}$. Define 
\begin{equation}\label{def_a_and_b}
    a(t)=\langle \partial_t\rho_t,[1-\Xi_t]\rangle - \langle \nabla\rho_t,\nabla\Xi_t\rangle\quad\text{and}\quad b(t)=\langle \partial_t\rho_t,\Xi_t\rangle + \langle \nabla\rho_t,\nabla\Xi_t\rangle,
\end{equation}
so that 
\begin{align*}
\mathfrak{L}(H)&=\mathfrak{L}(H^{(0)}+H^{(1)})=\mathfrak{L}_0(H^{(0)})+\mathfrak{L}(H_t(1)\Xi_t)+\mathfrak{L}(H_t(0)[1-\Xi_t])\\
&=\mathfrak{L}_0(H^{(0)})+\int_0^T a(t)H_t(1)\;dt+\int_0^T b(t)H_t(0)\;dt.
\end{align*}

Denote by $\Upsilon_t:\mathbb{R}^2\to\mathbb{R},\;0\leq t\leq T$ the map 
\begin{equation*}
    \Upsilon_t(x,y)=\zeta_t[x-y]^2+\mathfrak{b}^{-}(\rho_t(0),x)+\mathfrak{b}^{+}(\rho_t(1),y).
\end{equation*}

Let $\Phi_t:\mathbb{R}^2\to\mathbb{R},t\geq0$, denote its Legendre Transform
\begin{equation}\label{legendre_transf.bd}
    \Phi_t(a,b):=\sup_{x,y\in\mathbb{R}}\{ax+by-\Upsilon_t(x,y)\}.
\end{equation}

With this notation, the previous lemma translates into the following.
\begin{lemma}\label{decomposition_lemma}
    Fix a path $\pi(t,dx)=\rho_t(x)\;dx$ in $ {\mathcal{D}}([0,T],\mathscr{M}_{ac})$. Assume that $ {\rho}$ is continuous in $\Omega_T$, smooth in time and that there exists an $\epsilon>0$ such that $\epsilon\leq  {\rho}_t(x)\leq 1-\epsilon,$  {for all} $(t,x)\in\Omega_T$ and that $I_{[0,T]}(\rho)<\infty$. We have that $I_{[0,T]}(\rho)=I_{[0,T]}^{(1)}(\rho)+I_{[0,T]}^{(2)}(\rho)$ where
    \begin{equation*}
        I_{[0,T]}^{(1)}(\rho)=\frac{1}{4}||\mathfrak{L}_0||^2_{-1,\sigma(\rho)}\quad\quad\text{and}\quad\quad I_{[0,T]}^{(2)}(\rho)=\int_0^T\Phi_t(a_t,b_t)\;dt,
    \end{equation*}
     {where $a_t$ and $b_t$ were introduced in \eqref{def_a_and_b}.}
\end{lemma}

Now, exploring the structure of the space $H^{-1}(\sigma(\rho))$ given by Lemma \ref{structureH-1}, the first term can be simplified. We present here only the statement since the proof can be carried out by the same steps as in \cite{FGLN23}.

\begin{lemma}\label{operators_with_M}
     Fix a path $\pi(t,dx)=\rho_t(x)\;dx$ in $ {\mathcal{D}}([0,T],\mathscr{M}_{ac})$. Assume $ {\rho}$  {to satisfy the same hypothesis as in Lemma \ref{decomposition_lemma}.} Then, there exists an $M\in L^2(\sigma(\rho)^{-1})$ satisfying 
     \begin{equation}
         L_0^{(\partial_t)}(H)=\int_0^T\langle M_s,\nabla H_s\rangle\;ds\quad,\quad\int_0^1 \frac{M_t}{\sigma(\rho_t)}\;dx=0
     \end{equation}
     for almost every $0\leq t\leq T$ and for all $H\in C_0^{0,1}(\Omega_T)$ such that 
     \begin{equation}
         I^{(1)}_{[0,T]}(\rho)=\frac{1}{4}\int_0^T\Big\{||M_t+\nabla \rho_t||^2_{\sigma(\rho_t)^{-1}}-R_t\Big\}\;dt {,}
     \end{equation}
     where $R_t=\Big\langle \frac{\nabla \rho_t}{\sigma(\rho_t)} \Big\rangle^2\frac{1}{\langle \sigma(\rho_t)^{-1}\rangle}$.
\end{lemma}
The next theorem summarizes the results of this section.
\begin{theorem}\label{decomposition_summary}
    Fix a path $\pi(t,dx)=\rho_t(x)\;dx$ in $ {\mathcal{D}}([0,T],\mathscr{M}_{ac})$. Assume that $ {\rho}$ is continuous in $\Omega_T$, smooth in time and that there exists an $\epsilon>0$ such that $\epsilon\leq \rho_t(x)\leq 1-\epsilon,$  { for all } $(t,x)\in\Omega_T$ and that $I_{[0,T]}(\rho)<\infty$. Then,
    \begin{equation}
        I_{[0,T]}(\rho)=I^{(1)}_{[0,T]}(\rho)+I^{(2)}_{[0,T]}(\rho),
    \end{equation}
    for
    \begin{align}
        I^{(1)}_{[0,T]}(\rho)=\frac{1}{4}&\int_0^T\Big\{||M_t+\nabla u_t||^2_{\sigma(u_t)^{-1}}-R_t\Big\}\;dt,\\
        &I^{(2)}_{[0,T]}(\rho)=\int_0^T \Phi_t(a_t,b_t)\;dt,
    \end{align}
     {where $a_t$ and $b_t$ were introduced in \eqref{def_a_and_b} and $R_t$ was introduced at the end of Lemma \ref{operators_with_M}.}
\end{theorem}
\subsection{Upper Bound} \label{upperbound}

In this subsection, we prove the upper bound stated in Theorem \ref{main_result_large_dev}.
\subsubsection{Super-Exponential Estimates.}First we will establish the large deviations upper bound for compact sets and then for closed sets. Fix $G\in C^{1,2}( {\Omega_T})$. Consider the exponential martingale $\mathbb{M}_t^G$ defined by 
\begin{equation}\label{exponential_martingale}
    \mathbb{M}_t^G=\exp\left\{N\langle \pi^N_t, G_t\rangle-N\langle \pi^N_0, G_0\rangle -\int_0^t e^{-N\langle \pi^N_s, G_s\rangle}(\partial_s+\mathscr{L}_N) e^{N\langle \pi^N_s, G_s\rangle}\;ds\right\}.
\end{equation}
Elementary computations show that
\begin{align*}
  \mathbb{M}_t^G &= \exp\{N\langle \pi_t^N,G_t\rangle-N\langle \pi_0^N,G_0\rangle- N\int_0^t \langle \pi_s^N, {\Delta_N}G_s\rangle\;ds-N\int_0^t\langle\pi_s^N,\partial_s G_s\rangle\;ds\\
  &-\frac{1}{2}\int_0^t \sum_{x\in\Lambda_N}[( {\nabla^+_N} G_s)(x)]^2[\eta_s(x+\epsilon_N)-\eta_s(x)]^2\;ds\\
  &- N\int_0^t\sum_{\xi\in {\{0,1\}^{\Sigma^+_l}}}R^{+}(\Pi_{N-l-1,N-1}\eta_s,\xi)[e^{\sum_{x\in {\Sigma^+_l}}G_s(x)(\xi(x)-\eta_s(x))}-1]\;ds\\
  &- N\int_0^t\sum_{\xi\in {\{0,1\}^{\Sigma^-_l}} }R^{-}(\Pi_{1,l}\eta_s,\xi)[e^{\sum_{x\in {\Sigma^-_l}}G_s(x)(\xi(x)-\eta_s(x))}-1]\;ds\}.
\end{align*}

Some terms in the above expression cannot be expressed as a function of the empirical measure. Therefore, one  {needs the so-called} Super-Exponential estimates, to guarantee we can replace them  {by functions of the empirical measure}. Let $\Psi(\eta)=[\eta(x+\epsilon_N)-\eta(x)]^2$ and denote by $\Tilde{\Psi}(\alpha)$ the expectation of $\Psi$ with respect to the Bernoulli product measure with density $\alpha$:
\begin{equation*}
    \Tilde{\Psi}(\alpha):=\mathbb{E}_{\nu_\alpha}[\Psi(\eta)]=2\alpha(1-\alpha).
\end{equation*}
Define 
\begin{equation}\label{zetaminus}
     {\varphi_{\pm}(\eta, G_s)=\sum_{\xi\in\{0,1\}^{\Sigma^\pm_l}} R^{\pm}(\Pi\eta,\xi)\left[e^{\sum_{x\in\Sigma^\pm_l} G_s(x)(\xi(x)-\eta(x))}-1\right].}
\end{equation}
 {Recall the definition of $\mathfrak{b}^+$ given in \eqref{bfrakplus} and $\mathfrak{b}^-$ given in \eqref{bfrakminus}.} For $\epsilon>0$, let 
\begin{equation}
    V^{G,\Psi}_{N,\epsilon}(t,\eta_t):=\frac{1}{N}\sum_{x\in\Lambda_N}[(\nabla^+_N G_s)(x)]^2\{\Psi(\eta)-\Tilde{\Psi}(\overrightarrow{\eta}^{\lfloor \epsilon N\rfloor}_{t})\},
\end{equation}
\begin{equation}
     {W^{G,-}_{N,\epsilon}(t,\eta_t):=\Big\{\varphi_{-}(\eta_t,G_t)-\mathfrak{b}^{-}\Big(\overrightarrow{\eta}^{\lfloor \epsilon N\rfloor}_{t}(\epsilon_N),G_t(0)\Big)\Big\},}
\end{equation} and
\begin{equation}
     {W^{G,+}_{N,\epsilon}(t,\eta_t):=\Big\{\varphi_{+}(\eta_t, G_t)-\mathfrak{b}^{+}\Big(\overleftarrow{\eta}^{\lfloor \epsilon N\rfloor}_{t}(\tau_N),G_t(1)\Big)\Big\}.}
\end{equation}
\begin{theorem}{(Super-Exponential Estimate)}\label{super_exp_estimate}
    Fix any  {$G\in C^{0,1}(\Omega_T)$} and  {$(\eta^N)_{N\geq 1}$ a sequence of configurations}. Then, for any $\delta>0$  {and $T>0$,}
    \begin{equation}\label{superexpbulk}
    \lim_{\epsilon\to0}\limsup_{N\to\infty}\frac{1}{N}\log \mathbb{P}_{\eta^N}\Big[\Big|\int_0^T V^{G,\Psi}_{N,\epsilon}(t,\eta_t)\;dt \Big|>\delta \Big]=-\infty
    \end{equation}
    \begin{equation}\label{superexpbd}
    \lim_{\epsilon\to0}\limsup_{N\to\infty}\frac{1}{N}\log \mathbb{P}_{\eta^N}\Big[\Big|\int_0^T W^{G,\pm}_{N,\epsilon}(t,\eta_t)\;dt \Big|>\delta \Big]=-\infty
    \end{equation}   
\end{theorem}
\begin{proof}
The proof of \eqref{superexpbulk} in analogous to the one presented in \cite{bertini2009dynamical}. In order to prove the second one, we proceed as usual.  {We will prove only the result for the left boundary since for the right the argument is completely analogous.} First in view of the inequality
\begin{equation*}
        \limsup_{N\to\infty}\frac{1}{a_N}\log(b_N+c_N)\leq\max\{\limsup_{N\to\infty}\frac{1}{a_N}\log b_N, \limsup_{N\to\infty}\frac{1}{a_N}\log c_N\}
\end{equation*}
we can ignore the absolute value. For $a>0$, by an exponential Chebyshev inequality 
\begin{equation*}
        \mathbb{P}_{\eta^N}\Big[\int_0^T W^{G,-}_{N,\epsilon}(t,\eta_t)\;dt>\delta\Big]\leq e^{-a\delta N}\mathbb{E}_{\eta^N}\Big[\exp\Big\{aN\int_0^T  W^{G,-}_{N,\epsilon}(t,\eta_t)\;dt\Big\}\Big].
\end{equation*}
Thus, \eqref{superexpbd} is bounded above by 
\begin{equation} \label{to_bound_before_a_infty}
        -a\delta+\limsup_{N\to\infty}\frac{1}{N}\log\mathbb{E}_{\eta^N}\Big[\exp\Big\{aN\int_0^T  W^{G,-}_{N,\epsilon}(t,\eta_t)\;dt\Big\}\Big].
\end{equation}
Therefore, if we can prove that the second term in \eqref{to_bound_before_a_infty} is bounded by some $C_0 >0$ that does not depend on $N$, $\epsilon$ or $a$, then the proof finishes by sending $a\to\infty$.

By Feynmann-Kac's formula, see \cite[Appendix 1, Section 7]{KL99},
\begin{equation*}
\limsup_{N\to\infty}\frac{1}{N}\log\mathbb{E}_{\eta^N}\Big[\exp\Big\{aN\int_0^T  W^{G,-}_{N,\epsilon}(t,\eta_t)\;dt\Big\}\Big] \leq \limsup_{N\to\infty}\frac{ {T}}{N}\sup_{f}\Big\{aN\langle W^{G,-}_{N,\epsilon},f\rangle_{\nu^N_\rho}-\mathfrak{D}(\sqrt{f}) \Big\},
\end{equation*}where the supremum runs over all densities $f$  {and where $\mathfrak{D}$ denotes the Dirichlet Form introduced in Appendix \ref{app_dirichlet form}}.
The next step is to find a relationship between the inner product and the Dirichlet form.  {Since for each $\xi\in {\{0,1\}^{\Sigma^-_l}}$, $\varphi_{-}(\eta, G_s)$ is a bounded function that depends only on $G_s, \eta(1), \dots,\eta(l)$, then, by Lemma \ref{lemma_basis_L^2}, one can write it as a sum of linear combinations of products of the variables $\eta(1), \dots, \eta(l)$ as  
\begin{equation*}
  \varphi_{-}(\eta, G_s)=\sum_{\xi\in {\{0,1\}^{\Sigma^-_l}}}\sum_{A(\xi)\subset {\Sigma^-_l}}c(s,A(\xi))\prod_{x\in A(\xi)}\eta(x),
\end{equation*}
for some coeficients $c(s,A(\xi))$ and where $A(\xi)=\{x\in {\Sigma^-_l} : \xi(x)=1\}$. Hence, we can write
\begin{equation*}
    W^{G,-}_{N,\epsilon}(t,\eta_t)=\sum_{\xi\in {\{0,1\}^{\Sigma^-_l}}}\sum_{A(\xi)\subset {\Sigma^-_l}}c(t,A(\xi))\Big\{\prod_{x\in A(\xi)}\eta(x)-\left[\overrightarrow{\eta}^{\lfloor \epsilon N \rfloor}(\epsilon_N) \right]^{|A(\xi)|}\Big\}.
\end{equation*}
}From the bounds in Lemmas \ref{lemma_estimate_dirichlet_form}  {and the proof of Lemma \ref{repl_lemma_2}}, we get that
    \begin{align*}
        &\frac{t}{N}\sup_{f}\Big\{aN\langle W^{G,-}_{N,\epsilon},f\rangle_{\nu^N_\rho}-\mathfrak{D}(\sqrt{f}) \Big\}\\
        &\leq t\sup_{f}\Big\{ {\frac{1}{4N}\mathcal{D}(\sqrt{f})+O(\epsilon)+O(N^{-1})-\frac{1}{4N}\mathcal{D}(\sqrt{f})+C^{bulk}+\frac{C^{bd}}{N}}\Big\}.
    \end{align*}
    Therefore,  {after taking $N\to\infty$ and $\epsilon\to0$, we get that the second term in equation \eqref{to_bound_before_a_infty} is bounded by a constant independent of $a$}, and we are done.
\end{proof}
\subsubsection{Energy Estimates} To properly state the next result we will introduce a family of smooth approximations. Consider the function $\phi:\mathbb{R}\to[0,\infty)$ defined by
\begin{equation*}
    \phi(r)=\frac{1}{Z}\exp{\Big\{-\frac{1}{(1-r^2)}\Big\}}\boldsymbol{1}\{|r|<1\},
\end{equation*}
where $Z$ is a renormalization constant. For each $\delta>0$, let \begin{equation}\label{approx_unity}
    \phi^\delta(r)=\frac{1}{\delta}\phi\Big(\frac{r}{\delta}\Big).
\end{equation}
Fix a decreasing sequence $\{U_\varepsilon\}_{\varepsilon>0}$ such that $U_\varepsilon\downarrow1$ and let $\phi^\delta$ denote a approximation of unity \eqref{approx_unity}. For $\pi\in\mcb{M}$, let 
\begin{equation}
    \Theta_\varepsilon(\pi)(dx)=\frac{1}{U_\varepsilon}\int_0^1\phi^\delta(y-x)\pi(dy)\;dx.
\end{equation}
We set $\pi^{N,\varepsilon}=\Theta_\varepsilon(\pi^N)\in\mcb{M}_{ac}$, and denote its density by $\rho^{N,\varepsilon}.$
\begin{theorem}\label{energy_estimate_LDP}
    Fix a sequence $\{H_j:j\geq 1\}$ of functions in $C^{0,1}(\Omega_T)$ with compact support in $[0,T]\times(0,1)$ and a sequence $\{\eta^N:N\geq 1\}$ of configurations. There exists a finite constant $C=C(R^{\pm},l)$ such that 
    \begin{equation}
    \limsup_{\varepsilon\to0}\limsup_{N\to\infty}\frac{1}{N}\log\mathbb{P}_{\eta^N}^N\Big(\max_{1\leq j\leq k}Q_{H_j}(\rho^{N,\varepsilon})\geq m\Big)\leq -m +C(T+1),
    \end{equation}
    for all $k,m\geq1.$
\end{theorem}
This theorem corresponds to the Lemma 3.3 in \cite{bertini2009dynamical}. The details are left to the reader.
The proof of the large deviations upper bound now follows from Theorems \ref{super_exp_estimate} and \ref{energy_estimate_LDP}. We refer to \cite[Chapter 10]{KL99} for the precise steps.

\subsection{Lower Bound} \label{lowerbound}

In this subsection, we prove the lower bound in Theorem \ref{main_result_large_dev}. The strategy is by now classic but we outline here its main steps for convenience to the reader. 

\subsubsection{Perturbed Process} Fix a function $G\in C^{1,2}(\Omega_T)$. Consider the exponential martingale introduced in \eqref{exponential_martingale} and define a new probability measure as 
\begin{equation*}
    \mathbb{P}^G_N(\cdot):=\mathbb{E}(\mathds{1}\{\cdot\}\mathbb{M}^G_t).
\end{equation*}
It is well known that this perturbation creates a weak asymmetry on the  {initial process}, generating a new Markov process which we call the perturbed process. We will denote the generator of this new process by $L^G_N$ and the induced dynamics will be called the tilted dynamics. The first step in the proof of the large deviations lower bound is to obtain a law of large numbers for the tilted dynamics. First, we need to introduce some notation.
\par For $0<\alpha<1$ and $M\in\mathbb{R}$, let 
 {\begin{equation*}
    \mathfrak{p}^{\pm}(\alpha,M)=\mathbb{E}_{\nu_\alpha}\left[\sum_{\xi\in {\{0,1\}^{\Sigma^\pm_l}}}R^{\pm}(\Pi\eta,\xi)M\sum_{y\in {\Sigma^\pm_l}}(\xi(y)-\eta(y))\;e^{M\sum_{x\in {\Sigma^\pm_l}}[\xi(x)-\eta(x)]}\right].
\end{equation*}}

\begin{definition}[Weak solution of \eqref{LLN_perturbed_process}] Fix $G \in C^{0,1}(\Omega_T)$. A function  {$\rho \in L^2([0,T],H^1([0,1]))$} such that $0\leq \rho \leq 1$ a.e. in $\Omega_T$ is called a weak solution of 
\begin{align} \label{LLN_perturbed_process}
    \begin{cases}
    \partial_t \rho = \Delta \rho-2\nabla\{\sigma(\rho)\nabla G\}\\
    \nabla \rho_t(0) -2\sigma(\rho_t(0))\nabla G_t(0) = -\mathfrak{p}^-(\rho_t(0),G_t(0))\\
    \nabla \rho_t(1) -2\sigma(\rho_t(1))\nabla G_t(1)= \mathfrak{p}^+(\rho_t(1),G_t(1))\\
    \rho(0,x) = \gamma(x)
    \end{cases} 
\end{align} if 
\begin{align*}
    \int_0^1 \left[\rho_t H_t -\gamma H_0 -\int_0^t \rho_s\partial_s H_s\;ds \right] dx &=\int_0^t\int_0^1\{-\nabla\rho_s\nabla H_s+ 2\sigma(\rho_s)\nabla G_s\nabla H_s\}\;dx\;ds\\
    &+ \int_0^t \{\mathfrak{p}^-(\rho_s(0),G_s(0))H_s(0)+\mathfrak{p}^+(\rho_s(1),G_s(1))H_s(1)\} ds,
\end{align*}
for all $t\in [0,T]$ and all test functions $H\in C^{1,2}( {\Omega_T})$.
\end{definition}

 {The initial value problem \eqref{LLN_perturbed_process} has a weak solution, which is a consequence of the tightness of the sequence of probability measures associated with the perturbed process.  {Moreover,} this solution is unique - see Theorem \ref{th_uniqueness} taking $\alpha=1$ and $\beta = 0$.}

\begin{theorem} (Hydrodynamic Limit for the perturbed process)\label{hydrolimit_pertubed_process}
    Fix a  {measurable }density profile $\gamma:[0,1]\to[0,1]$ and let $\{\mu_N\}_N$ be a sequence of probability measures associated to it. Then, for each $t>0,\;\delta>0$ and $f\in C([0,1])$ we have
    \begin{equation*}
    \lim_{N\to\infty}\mathbb{P}_{\mu_N}^{G}\Big(|\langle\pi_t^N,f\rangle-\langle\rho^G_t,f\rangle|>\delta\Big)=0,
    \end{equation*}
    where $\rho^G\in L^2(0,T;H^1[0,1])$ is the unique weak solution of 
    \eqref{LLN_perturbed_process}.
\end{theorem}
The proof of Theorem \ref{hydrolimit_pertubed_process} is a consequence of the usual strategy, i.e. tightness and characterization of limit point, as we did in the first part of the article to prove the Hydrodynamic Limit of the model, jointly with the uniqueness result given in Appendix \ref{appB} {, Theorem \ref{th_uniqueness} taking $\alpha=1$ and $\beta = 0$.} We leave the details to the reader.

By Theorem \ref{hydrolimit_pertubed_process} and elementary computations, we have the following classical result which can be proven following the steps of Lemma 5.4 of \cite[Appendix 1]{KL99}.
\begin{lemma}\label{rate_function_rel_entropy}
For each continuous function $\gamma:[0,1]\to[0,1]$ and each smooth function $G\in C^{1,2}( {\Omega_T})$
\begin{equation}
    \lim_{N\to\infty}\frac{1}{N} H( \mathbb{P}^{G}_{\eta^N}\;|\;\mathbb{P}_{\eta^N})=I_{[0,T]}(\rho^G|\gamma),
\end{equation}
where $ {\rho^G}$ is the unique weak solution of \eqref{hydrolimit_pertubed_process}.
\end{lemma}
\begin{definition}
    Given $\gamma\in\mcb M_{ac}$, let $\Pi_\gamma$ be the set of paths $\pi(t,dx)=\rho_t(x)dx$ in $\mathcal{D}([0,T],\mcb M_{ac})$ such that
    \begin{enumerate}
        \item There exists a $t_0>0$ such that $\rho$ follows the hydrodynamic equation \eqref{PDEproblem} in the interval $[0,t_0]$.
        \item For every $0<\delta\leq T,\;\exists\;\varepsilon>0$ such that $\varepsilon\leq\rho_t(x)\leq1-\varepsilon$ for all $(t,x)\in  {[\delta,T]}\times[0,1].$
        \item $\rho$ is smooth in $(0,T]\times[0,1].$
    \end{enumerate}
\end{definition}
Fix an open set $\mathcal{O}\in  {\mathcal{D}}([0,T],\mcb{M})$ and a path $\pi(t,dx)=\rho_t(x)dx\in\Pi_\gamma$. By Lemma \ref{rate_function_rel_entropy} and the proof presented in \cite[Chapter 10]{KL99}, we get
\begin{equation*}
    \liminf_{N\to\infty}\frac{1}{N}\log \mathbb{P}_{\eta^N}(\pi^N\in\mathcal{O})\geq -\inf_{\mathcal{O}\cap \Pi_\gamma}I_{[0,T]}(\rho|\gamma),
\end{equation*}
and the result is proved for all paths $\pi\in\Pi_\gamma$. Therefore, to conclude the proof of Theorem \ref{main_result_large_dev}, we need to show that every trajectory $\pi$ with finite rate function, can be approximated by trajectories $\pi_n\in\Pi_\gamma$ in such a way that $I_{[0,T]}(\pi_n)\to I_{[0,T]}(\pi)$. This is the content of the next section,  {which contains} the so-called I-density Theorem.
\subsubsection{$I_{[0,T]}(\cdot)$-density} In this section, we prove that any trajectory $\pi\in  {\mathcal{D}}([0,T],\mcb{M})$ with finite rate function can be approximated by a sequence of smooth trajectories $\{\pi^n: n\geq1\}$ such that
\begin{equation*}
    \pi^n\rightarrow\pi\quad\text{and}\quad I_{[0,T]}(\pi^n|\gamma)\rightarrow I_{[0,T]}(\pi|\gamma).
\end{equation*}
In order to state the result more precisely, we need to introduce some terminology. 

\begin{definition}
    A subset $\mathcal{A}\subset \mathcal{D}([0,T],\mcb M)$ is said to be $I_{[0,T]}(\cdot|\gamma)$-dense if for any $\pi\in\mathcal{D}([0,T],\mcb M)$ such that $I_{[0,T]}(\pi|\gamma)<+\infty$, there exists $\{\pi^n:n\geq1\}$ in $\mathcal{A}$ with $\pi^n\to\pi$ in $\mathcal{D}([0,T],\mcb M)$ and $I_{[0,T]}(\pi^n|\gamma)\to I_{[0,T]}(\pi|\gamma)$.
\end{definition}
The main result of this section reads as follows.
\begin{theorem}\label{thm_I_density}
    For all $\gamma:[0,1]\to[0,1]$, the set $\Pi_\gamma$ is $I_{[0,T]}(\cdot|\gamma)$-dense.
\end{theorem}

The proof of Theorem \ref{thm_I_density} is divided in several  {lemmas.} It consists in approximating in several steps a general trajectory $\pi$ by a sequence
of profiles, smoother at each  {step.} Throughout this section, we will denote by $\rho^{(\gamma)}: {\Omega_T}\to[0,1]$ the unique weak solution of equation \eqref{PDEproblem}.
\begin{definition}
    Let $\Pi_1$ denote the set of all trajectories $\pi_t(dx)=\rho_t(x)dx$ in $ {\mathcal{D}}_{\mathcal{E}}([0,T],\mcb{M}_{ac})$, for which $\exists\;\delta>0$, such that $\rho(t,\cdot)=\rho^{(\gamma)}(t,\cdot)$, for all $0\leq t\leq\delta.$
\end{definition}
\begin{lemma}
    The set $\Pi_1$ is $I_{[0,T]}(\cdot|\gamma)$-dense.
\end{lemma}
\begin{proof}
    Fix a path $\pi_t(dx)=\rho_t(x)\;dx$ in $ {\mathcal{D}}([0,T],\mcb{M})$ with finite rate function. For each $\delta>0,$ let $\pi^\delta(t,dx)=\rho^\delta(t,x)\;dx$ be the path defined by 
    \begin{align*}
    \rho^\delta(t,x)=
    \begin{cases}
       \rho^{(\gamma)}(t,x),\quad\quad\quad\text{ if } t\in[0,\delta]\\
       \rho^{(\gamma)}(2\delta-t,x),\quad\text{ if } t\in[\delta,2\delta]\\
       \rho(t-2\delta,x),\quad\text{ if } t\in[2\delta,T].
    \end{cases}
    \end{align*}
    We have that $\pi^\delta\rightarrow\pi$ in $ {\mathcal{D}}([0,T],\mcb{M})$ as $\delta\downarrow0$ and $\pi^\delta\in\Pi_1.$ Therefore, it remains to show only that, as $\delta\downarrow0$, $I_{[0,T]}(\pi^\delta|\gamma)\rightarrow I_{[0,T]}(\pi|\gamma)$. By Theorem \ref{rate_function_lsc_compac_ls}, the rate function is lower semicontinuous. Hence, $I_{[0,T]}(\pi|\gamma)\leq\liminf_{\delta\to0}I_{[0,T]}(\pi^\delta|\gamma)$. In order to finish the proof, we need therefore to show 
    \begin{equation*}
       \limsup_{\delta\to0}I_{[0,T]}(\pi^\delta|\gamma)\leq I_{[0,T]}(\pi|\gamma).
    \end{equation*}
    To do so, we decompose the rate function $I_{[0,T]}(\pi^\delta|\gamma)$ into the sum of the contributions of each time interval $[0,\delta], [\delta,2\delta] \text{ and } [2\delta,T]$. By \eqref{cost_subinterval},
     {\begin{equation*}
        I_{[0,T]}(\pi^\delta)\leq I_{[0,\delta]}(\pi^\delta)+I_{[0,\delta]}(\tau_\delta\pi^\delta) + I_{[0,T]}(\pi),
    \end{equation*}}where $\tau_\delta\pi^\delta_t(dx)=\pi^\delta_{t+\delta}(dx).$ By Corollary \ref{rate_function_0}, the contribution of the first term is zero. Let $v^\delta=\tau_\delta\rho^\delta$. As $v^\delta(t,x)=\tau_\delta\rho^\delta(t,x)=  {\rho^{(\gamma)}(2\delta-t)}$ in the time interval $[0,\delta]$, the density $v^\delta$ solves the backward heat equation with the opposite boundary conditions  {
    \begin{equation}\label{backward_heat}
        \begin{cases}
            \partial_t v^\delta=-\Delta v^\delta\\
            \nabla v^\delta (t,0)=F_{-}(v^\delta (t,0))\\
            \nabla v^\delta (t,1)=-F_{+}(v^\delta (t,1))\\
            v^\delta(0,x)=u(\delta,x)
        \end{cases}
    \end{equation}}
    Recall  {the definition of $J_{T,H}$ in \eqref{def_J}.} By \eqref{def_weak_sol_eq_hydrodynamcis}, for each $H\in C^{1,2}(\Omega_T)$, $J_{\delta, H}(v^\delta)$ is equal to
     {\begin{align*}
        \int_0^\delta \{2\langle\nabla \rho^{(\gamma)}_t,\nabla H_t\rangle-\langle\sigma(\rho^{(\gamma)}_t),(\nabla H_t)^2\rangle\}\;dt -\int_0^\delta\{\hat{\mathfrak{b}}^{-}(\rho^{(\gamma)}_t(0),H_t(0))+\hat{\mathfrak{b}}^{+}(\rho^{(\gamma)}_t(1),H_t(1))\}dt,
    \end{align*}}
where 
 {\begin{equation*}
    \hat{\mathfrak{b}}^{\pm}(\alpha,M):=\mathbb{E}_{\nu_\alpha}\left[\sum_{\xi\in {\{0,1\}^{\Sigma^\pm_l}}}R^{\pm}(\Pi\eta,\xi)\Big(e^{M\sum_{x\in {\Sigma^\pm_l}}(\xi(x)-\eta(x))}-1+M\sum_{x\in {\Sigma^\pm_l}}(\xi(x)-\eta(x))\Big)\right].
\end{equation*}}
By Cauchy-Schwarz's inequality and Lemma \ref{our_B5}, we have that 
    \begin{equation*}
        J_{\delta,H}(v^\delta)\leq C(\delta),
    \end{equation*}
    where $C(\delta)$ is a constant independent of $H$ and that vanishes as $\delta\to 0$. Therefore, we conclude that 
    \begin{equation*}
         {\lim_{\delta\to0}I_{[0,\delta]}(\tau_\delta \rho^\delta)=0.}
    \end{equation*}
\end{proof}
\begin{definition}
    Let $\Pi_2$ be the set of all trajectories $\pi_t(dx)=\rho_t(x)\;dx$ in $\Pi_1$ such that $\forall \delta>0,\;\exists\;\epsilon>0$ such that $\epsilon\leq\rho_t(x)\leq1-\epsilon$, for all $(t,x)\in[\delta,T]\times[0,1]$.
\end{definition}
\begin{lemma}
    The set $\Pi_2$ is $I_{[0,T]}(\cdot|\gamma)$-dense.
\end{lemma}
\begin{proof}
    Consider a path $\pi_t(dx)=\rho_t(x)\;dx$ in $\Pi_1$ with finite rate function. Define for $\varepsilon\in(0,1)$ the path $\pi^\varepsilon_t(dx)=\rho^\varepsilon(t,x)\;dx$ where $\rho^\varepsilon=(1-\varepsilon)\rho+\varepsilon \rho^{(\gamma)}.$ Let us show that this sequence belongs to $\Pi_2$. First, since $\pi\in\Pi_1$, by definition there exists a $\delta>0$ such that $\pi_t^\varepsilon=\pi_t$ for $0\leq t\leq\delta.$ Hence, $\pi^\varepsilon$ follows the hydrodynamic equation on that time interval. Moreover, since the energy $Q_{[0,T]}$ is convex, we have that $Q_{[0,T]}(\rho^\varepsilon)$ is finite, because both $Q_{[0,T]}(\rho)$ and  {$Q_{[0,T]}(\rho^{(\gamma)})$} are. Therefore, $\pi^\varepsilon\in\Pi_1$. Now, by Lemma \eqref{bounded_away_0_1}, for all $\delta>0$, there exists an $\alpha>0$ such that $\alpha\leq\rho_t(x)\leq1-\alpha,$ for all $\delta\leq t\leq T.$ Since $0\leq\rho\leq1$, this property still holds for $\rho^\varepsilon$ (maybe with a different $\alpha$). Thus, $\pi^\varepsilon\in\Pi_2$.\par
    It is clear that $\pi^\varepsilon\to\pi$, as $\varepsilon\to0$. Hence, in view again of the lower-semicontinuity of the rate function, it suffices to prove that 
    \begin{equation*}
        \limsup_{\varepsilon\to\infty}I_{[0,T]}(\pi^\varepsilon|\gamma)\leq I_{[0,T]}(\pi|\gamma).
    \end{equation*}
 {By Theorem \ref{rate_function_lsc_compac_ls}, the rate function is convex, and we are done.}
    \end{proof}
\begin{definition}
    Let $\Pi_3$ be the set of all trajectories $\pi_t(dx)=\rho_t(x)dx$ in $\Pi_2$ whose density $\rho$ is continuous in $(0,T]\times[0,1]$ and smooth in time, meaning that for all $x\in[0,1],\;\rho(\cdot,x)\in C^{\infty}((0,T])$.
\end{definition}
\begin{lemma}
     The set $\Pi_3$ is $I_{[0,T]}(\cdot|\gamma)$-dense.
\end{lemma}
The proof that we can approximate any path with finite rate function by paths in $\Pi_3$ is completely analogous to the one presented in \cite{FGLN23}. It relies on Lemma \ref{our_B5} and Lemma \ref{solution_is_smooth}. We left the details for the reader.
\begin{definition}
    Let $\Pi_4$ be the set of all trajectories $\pi_t(dx)=\rho_t(x)dx$ in $\Pi_3$ whose density $\rho(t,\cdot)$ belongs to $C^\infty([0,1])$ for any $t\in(0,T].$
\end{definition}
Note that $\Pi_4=\Pi_\gamma$. For the next lemma we will need to introduce an auxiliary function. Consider a $\delta>0$ and let $\chi:[0,T]\to[0,1]$ be a smooth, non-decreasing function such that 
\begin{equation}\label{def_chi}
    \begin{dcases*}
        \chi(t)=0,\quad \text{if}\quad t\in[0,\delta],\\
        0<\chi(t)<1, \quad \text{if}\quad t\in(\delta,2\delta),\\
        \chi(t)=1, \quad \text{if}\quad t\in[2\delta,T].\\
    \end{dcases*}
\end{equation}
Set $\chi_n(t)=\frac{\chi(t)}{n}$ for $n\geq 1.$ Note that $\chi_n(t)=1/n$, for $t\geq 2\delta.$\par
Let $(P_t^{(R)}:t\geq0)$ denote the semigroup associated to the Laplacian with Robin boundary condition and $(P_t^{(D)}:t\geq0)$ the same operator with Dirichlet boundary conditions. We refer to \cite[Appendix A]{FGLN23}, for a precise definition and several properties  {of} those semigroups.\par 
Fix $\pi_t(dx)=\rho_t(x)dx$ in $\Pi_3$ with finite rate function $I_{[0,T]}(\rho|\gamma).$ In particular, $\pi$ belongs to $\Pi_1$, hence its density $\rho$ solves the hydrodynamic equation \eqref{PDEproblem} in some time interval $[0,3\delta]$, for some $\delta>0.$ Consider the sequence $\pi^n_t(dx)=\rho^n(t,x)\;dx$ given by 
\begin{equation}\label{approximation_i-density_4}
    \rho_t^n=w_t+P^{(D)}_{\chi_n(t)}[\rho_t-w_t],
\end{equation}
where $w_t=\rho_t(0)+[\rho_t(1)-\rho_t(0)]x$ and $\chi_n(t)$ was defined in \eqref{def_chi}.
\begin{lemma}
    For any $n\geq 1$, $\pi^n_t(dx)=\rho^n(t,x)\;dx$ belongs to $\Pi_4$ and the trajectory $\rho^n$ has finite energy.
\end{lemma}

The proof of this result is similar to the one presented in \cite{FGLN23}. Applying Lemma \ref{our_B5} and Lemma \ref{bounded_away_0_1}, we get the same result.\par
Now we use the decomposition of the rate function introduced in Section \ref{deconstruction_sec} to complete the proof of Theorem \ref{thm_I_density}. Note that $\Pi_4=\Pi_\gamma.$ Recall that the main goal of this decomposition is to write the rate function as the sum of two independent variational problems,  {with one of them running only over functions that vanish on the boundary of the interval.} This fact will be exploited several times in the proof of the next lemma.  {In a nutshell, to treat the contribution due to the evolution in the bulk,} this deconstruction allows us to perform several integrations by parts arguments and to use Poincaré's Inequality to obtain sharp bounds.  {Moreover, using} the structure of the function $\Phi_t$,  {introduced in \eqref{legendre_transf.bd},} together with bounds on the functions $a_t$ and $b_t$, introduced in \eqref{def_a_and_b},  {we can pass to the limit.}

\begin{lemma}
    The set $\Pi_4$ is $I_{[0,T]}(\cdot|\gamma)$-dense.
\end{lemma}
\begin{proof}
    Let $\pi(t,dx)=\rho_t(x)dx$ be a path in $\Pi_3$ such that $I_{[0,T]}(\rho|\gamma)<\infty.$ By the previous Lemma, the approximation $\rho^n$ defined in \eqref{approximation_i-density_4} belongs to $\Pi_4.$ By definition, as $n\to\infty,$ $\pi^n$ converges to $\pi$ in $ {\mathcal{D}}([0,T],\mathscr{M})$. Hence, again in view of the lower-semicontinuity of the rate function, it remains only to show that $\limsup_{n\to\infty}I_{[0,T]}(\pi^n|\gamma)\leq I_{[0,T]}(\pi|\gamma).$\\
    By Proposition \ref{cost_subinterval}, the cost of the trajectory $\pi^n$ on the time interval $[0,T]$ is bounded by the sum of its costs on the intervals $[0,\delta]$ and $[\delta,T]$. Since for $0\leq t\leq\delta,\;\rho^n=\rho^{(\gamma)}$, we have that 
    \begin{equation}
        I_{[0,\delta]}(\rho^n)=0.
    \end{equation}
Therefore, we turn the analysis to the time interval $[\delta,T].$ In view of Proposition \ref{cost_subinterval}, the cost on this interval is equal to $I_{[0,T-\delta]}I(\tau_\delta\rho^n)$. Letting $v=\tau_\delta\rho,\;v^n=\tau_\delta\rho^n,\;\hat{w}=\tau_\delta w,\;\hat{\chi}_n(t)=\chi_n(t-\delta)$ and $T_\delta=T-\delta$ we can write $v_t^n=\hat{w}_t+P^{(D)}_{\hat{\chi}_n(t)}[v_t-\hat{w}_t]$ and with this notation, 
\begin{equation*}
    { I_{[0,T-\delta]}(\tau_\delta\rho^n)=I_{[0,T_\delta]}(v^n).}
\end{equation*}
    Hence, to conclude the proof of the Lemma, it suffices to show that 
    \begin{equation}
        \limsup_{n\to\infty}I_{[0,T_\delta]}(v^n)\leq I_{[0,T_\delta]}(v).
    \end{equation}
    \par
    By Lemma \ref{decomposition_lemma}, $I_{[0,T_\delta]}(v^n)=I^{(1)}_{[0,T_\delta]}(v^n)+I_{[0,T_\delta]}^{(2)}(v^n)$, which gives us that 
    \begin{equation}\label{decomposition_lemma4}
        I_{[0,T_\delta]}(v^n)=\frac{1}{4}||\mathfrak{L}_0||^2_{-1,\sigma(v^n)}+\int_0^{T_\delta}\Phi_t^{v_n}(a_t^n,b_t^n)\;dt,
    \end{equation}
    where $a^n_t$ and $b^n_t$ are defined as $a_t$ and $b_t$  {in \eqref{def_a_and_b}} with $v^n$ instead of $\rho.$ 
    Our goal now is to show that
    \begin{equation}
        \limsup_{n\to\infty}I^{(j)}_{[0,T_\delta]}(v^n)\leq I^{(j)}_{[0,T_\delta]}(v)
    \end{equation}
    for $j=1,2.$ Since the analysis for $j=1$ is completely analogous to the one presented in \cite{FGLN23},  {here we only take care of the case $j=2$.}
    
    We turn to the second term in \eqref{decomposition_lemma4}. Our goal is to provide good bounds on $a_t^n,\;b_t^n$ and $\Phi$ in order to perform a dominated convergence argument and pass to the limit. Recall that 
    \begin{equation*}
    a^n_t=\langle \partial_t v^n_t,[1-\Xi_t^n]\rangle - \langle \nabla v^n_t,\nabla\Xi_t^n\rangle\quad\text{and}\quad b^n_t=\langle \partial_t v_t^n,\Xi_t^n\rangle + \langle \nabla v^n_t,\nabla\Xi_t^n\rangle,
\end{equation*}
where $\Xi_t^n$ is given by \eqref{def_Xi} with $\rho$ replaced by $v^n_t.$ If we denote by $(P^{(N)}_t)_{t\geq 0}$ the semigroup associated with the Laplacian on $[0,1]$ with Neumann boundary conditions, then an easy computation shows that 
    \begin{equation} \label{relation_P_N_andP_D}
        \nabla P^{(D)}_t f=P^{(N)}_t \nabla f,
    \end{equation}
    for any $f\in C^1([0,1]).$ By the definition of $v_t^n$, we have that 
\begin{equation*}
\begin{split}
    a^n_t &=\langle(I-P^{(D)}_{\hat{\chi}_n(t)}) \partial_t \hat{w}_t,[1-\Xi_t^n]\rangle+\langle P^{(D)}_{\hat{\chi}_n(t)}\partial_t v_t,[1-\Xi_t^n]\rangle\\
    &-\langle (I-P^{(N)}_{\hat{\chi}_n(t)})\nabla\hat{w}_t,\nabla\Xi_t^n\rangle-\langle P^{(N)}_{\hat{\chi}_n(t)}\nabla v_t,\nabla\Xi_t^n\rangle\\
    &+\hat{\chi}'_n(t)\langle\Delta P^{(D)}_{\hat{\chi}_n(t)}[v_t-\hat{w}_t],[1-\Xi_t^n]\rangle.
    \end{split}
\end{equation*}
Since the semigroup $P^{(D)}_t$ is symmetric  {and commutes with the Laplacian}, the last line of the above display can be written as 
\begin{equation*}
    \hat{\chi}'_n(t)\langle v_t-\hat{w}_t,\Delta P^{(D)}_{\hat{\chi}_n(t)}[1-\Xi_t^n]\rangle.
\end{equation*}
Moreover, since $v_t-\hat{w}_t$ vanishes for $x=0$ and $x=1$,  {by integration by parts, the relation \eqref{relation_P_N_andP_D} and symmetry} of the semigroup $P^{(N)}_t$, we get  {that}
\begin{equation*}
\begin{split}
    \hat{\chi}'_n(t)\langle v_t-\hat{w}_t,\Delta P^{(D)}_{\hat{\chi}_n(t)}[1-\Xi_t^n]\rangle&=-\hat{\chi}'_n(t)\langle \nabla[v_t-\hat{w}_t],\nabla P^{(D)}_{\hat{\chi}_n(t)}[1-\Xi_t^n]\rangle\\
    &=\hat{\chi}'_n(t)\langle P^{(N)}_{\hat{\chi}_n(t)}\nabla[v_t-\hat{w}_t],\nabla\Xi_t^n\rangle.
    \end{split}
\end{equation*}
 {Since $v_t^n$ is bounded away from 0 and 1, by the definition of $\Xi_t^n$, we know that }there exists a constant $C_0>0,$ such that 
\begin{equation}\label{bound_on_Xi^n}
    |\Xi_t^n|\leq C_0\quad \text{and} \quad |\nabla \Xi_t^n|\leq C_0,
\end{equation}
for all $t\in[0,T_\delta]$ and for all $n \in \mathbb{N}$.  {Therefore, because $\hat{\chi}'_n(t) \to 0$ as $n\to\infty$, we prove that }
\begin{equation*}
    \lim_{n\to\infty}\hat{\chi}'_n(t)\langle P^{(N)}_{\hat{\chi}_n(t)}\nabla[v_t-\hat{w}_t],\nabla\Xi_t^n\rangle=0.
\end{equation*}
Moreover, since $P^{(N)}_{\hat{\chi}_n(t)},P^{(D)}_{\hat{\chi}_n(t)}\to I$ as $n\to\infty$, by \eqref{bound_on_Xi^n}, 
\begin{equation*}
    \begin{split}
        &\lim_{n\to\infty}\langle(I-P^{(D)}_{\hat{\chi}_n(t)}) \partial_t \hat{w}_t,[1-\Xi_t^n]\rangle=0\\
        &\lim_{n\to\infty}\langle (I-P^{(N)}_{\hat{\chi}_n(t)})\nabla\hat{w}_t,\nabla\Xi_t^n\rangle=0.
    \end{split}
\end{equation*}
Hence, we can rewrite $a_t^n$ as 
\begin{equation*}
    a_t^n=\langle P^{(D)}_{\hat{\chi}_n(t)}\partial_t v_t,[1-\Xi_t^n]\rangle-\langle P^{(N)}_{\hat{\chi}_n(t)}\nabla v_t,\nabla\Xi_t^n\rangle+\Theta_n,
\end{equation*}
where $\lim_{n\to\infty}\Theta_n=0$.  {Using the Cauchy-Schwarz inequality combined with the fact that $P^{(D)}_t, P^{(N)}_t$ are contractions in $L^2([0,1])$ and that it holds \eqref{bound_on_Xi^n},} we get that 
\begin{equation}\label{bound_on_a_t^n}
    |a_t^n|^2\leq C_0[1+\langle(\partial_t v_t)^2\rangle+\langle(\nabla v_t)^2 \rangle],
\end{equation}
for every $n\geq1$ and $t\in[0,T_\delta].$ Since $\Xi_t^n\to\Xi_t$ and $\nabla\Xi_t^n\to\nabla\Xi_t^n$ in $L^2([0,1])$, as $n\to\infty$, we conclude that 
\begin{equation}\label{limit-a_t}
    \lim_{n\to\infty}a_t^n=\langle \partial_t v_t,[1-\Xi_t]\rangle-\langle \nabla v_t,\nabla\Xi_t\rangle=:a_t.
\end{equation}
The same reasoning gives us that 
\begin{equation}\label{limit-b_t}
    \lim_{n\to\infty}b_t^n=\langle \partial_t v_t,\Xi_t\rangle+\langle \nabla v_t,\nabla\Xi_t\rangle=:b_t.
\end{equation}
\par
Now, since the function $\Phi_t$ is a Legendre transform, it is convex and is also continuous. Taking $x=y=0$, we get that $\Phi_t(a,b)\geq0.$ Moreover, letting $\Phi^{-}_t(a)$ and $\Phi_t^{+}(b)$ denote the Legendre transforms of $\mathfrak{b}^{-}(\rho_t(0),\cdot)$ and $\mathfrak{b}^{+}(\rho_t(1),\cdot)$, respectively, then we have 
\begin{equation*}
    0\leq \Phi_t(a,b)\leq \Phi^{-}_t(a)+\Phi^{+}_t(b).
\end{equation*}
     {Above,} $\Phi^{-}_t(a)$, for example, is given by 
    \begin{equation*}
    \Phi^{-}_t(a)=\sup_{x\in\mathbb{R}}\{ax-\mathfrak{b}^{-}(\rho_t(0), {x})\}.
    \end{equation*}
    In particular, 
    \begin{equation} \label{bound_on_phi_minus}
        0\leq\Phi^{-}_t(a)\leq C_0\Big\{1+\frac{|a|}{l}+\frac{|a|}{l}\log\Big(\frac{|a|}{l}\Big)\Big\}.
    \end{equation}
    Indeed, by definition of $\Phi^{-}_t(a)$ and $\mathfrak{b}^{-}(\rho_t(0),\cdot)$, we have that
\begin{align*}
    \Phi^{-}_t(a)
    =\sup_{x\in\mathbb{R}}\left\{ax-\mathbb{E}_{\nu_{\rho_t(0)}}\left[\sum_{\xi\in {\{0,1\}^{\Sigma^-_l}}}R^{-}(\eta,\xi)\Big(e^{x\sum_{i\in {\Sigma^-_l}}(\xi(i)-\eta(i))}-1\Big)\right]\right\}.
\end{align*}
Since $\sum_{i\in {\Sigma^-_l}}(\xi(i)-\eta(i))\in\{-l,...,l\},$ this expression is  {bounded by} 
\begin{align*}
    \sup_{x\in\mathbb{R}}\left\{ax-(e^{-|x|l}-1)\mathbb{E}_{\nu_{\rho_t(0)}}\left[\sum_{\xi\in {\{0,1\}^{\Sigma^-_l}}}R^{-}(\eta,\xi)\right]\right\}.
\end{align*}
 {Denoting by} $C:=\mathbb{E}_{\nu_{\rho_t(0)}}\Big[\sum_{\xi\in {\{0,1\}^{\Sigma^-_l}}}R^{-}(\eta,\xi)\Big]<\infty$, the above expression is bounded by
\begin{equation*}
    \sup_{x\in\mathbb{R}} \{ax+C(e^{-|x|l}-1)\}.
\end{equation*}
 {This} supremum is attained at $x=(-1)^{\mathbbm{1}_{a> 0}}\frac{1}{l}\log\Big(\frac{|a|}{Cl}\Big)$,  {and so the bound on $\Phi_t^{-}$ given in \eqref{bound_on_phi_minus} must hold. Analogously, we get the same estimate for $\Phi_t^{+}(b)$. This then allows us to conclude that}
    \begin{equation}\label{bound_Phi_DCT}
        \Phi_t(a,b)\leq C_0\left[1+\frac{|a|}{l}+\frac{|a|}{l}\log\Big(\frac{|a|}{l}\Big)+\frac{|b|}{l}+\frac{|b|}{l}\log\Big(\frac{|b|}{l}\Big)\right],
    \end{equation}
    for some constant $C_0>0.$
    
    Hence, since $\Phi^v$ is continuous, by \eqref{bound_Phi_DCT}, \eqref{limit-a_t},\eqref{limit-b_t} and the Dominated Convergence Theorem, 
    \begin{equation*}
        \lim_{n\to\infty}\int_0^{T_\delta}\Phi^{v^n}(a_t^n,b_t^n)\;dt=\int_0^{T_\delta} \Phi(a_t,b_t)\;dt,
    \end{equation*} from which we conclude that
    \begin{equation*}
        \lim_{n\to\infty} I^{(2)}_{[0,T_\delta]}(v^n)=I_{[0,T_\delta]}^{(2)}(v),
    \end{equation*}  {which completes our proof.}
\end{proof}

\appendix
\section{Entropy estimate, Dirichlet Form and Replacement Lemma} \label{appA}

 {In this section, we will provide the proof of the auxiliary results that are crucial to prove the hydrodynamic limit result given in Theorem \ref{th_hydrodynamics}.}

\subsection{Entropy estimate}
 {Let $\rho: \Lambda_N \to (0,1)$ be such that there exists constants $\alpha,\beta \in (0,1)$ for which, for every $x \in \Lambda_N$, it holds that $0<\alpha\leq\rho(x)\leq\beta<1$. Let $\nu_{\rho(\cdot)}^N$ be the Bernoulli product measure on $\Omega_N$ given by
\begin{equation*}
\nu_{\rho(\cdot)}^N(\eta)=\prod_{x\in\Lambda_N}\rho(x)^{\eta(x)}(1-\rho(x))^{1-\eta(x)}.
\end{equation*} The next lemma shows that the relative entropy of a sequence of probability measures with respect to $\nu_{\rho(\cdot)}^N$ is of order $N$.}
\begin{lemma}
 {There exist $C>0$, such that
\begin{equation} \label{entropy_estimate}
    H(\nu_N | \nu^N_{\rho(\cdot)}) \leq C N,
\end{equation} for every sequence of probability measures $\nu_N$ in $\Omega_N$ and for any $N \in \mathbb{N}$, where $H(\mu|\nu)$ represents the relative entropy between $\mu$ and $\nu$.}
\end{lemma}  {For a short proof of the previous lemma see, for example, the begining of Section 4 of \cite{FGS23}.}

\subsection{Estimate of the Dirichlet form}\label{app_dirichlet form} {Fix a density profile $\rho:[0,1]\to[0,1]$ such that $\rho$ is Lipschitz and bounded away from 0 and 1.
We define the Dirchlet Form $\mathfrak{D}(\cdot)$  by 
\begin{equation}
\mathfrak{D}(\sqrt{f}) = -\langle\mathscr{L}_N\sqrt{f},\sqrt{f}\rangle_{\nu_{\rho(\cdot)}^N},
\end{equation}where $\nu_{\rho(\cdot)}^N$ is the Bernoulli product measure with marginals of parameter $\rho$ and $f$ is any density with respect to $\nu_{\rho(\cdot)}^N$.}

The next lemma provides an upper bound for  {$\mathfrak{D}$}.

\begin{lemma} \label{lemma_estimate_dirichlet_form}
Let $\rho:[0,1]\to[0,1]$ be a Lipschitz continuous function that is bounded away from zero and one, i.e. there exists $0 < \alpha < \beta < 1$ for which $\alpha \leq ||\rho||_\infty \leq \beta$. Given a density $f$ with respect to $\nu_{\rho(\cdot)}^N$, denote
\begin{equation} \label{notation_E_eta}
    E_\eta(\sqrt{f},\zeta) := [\sqrt{f}(\zeta)-\sqrt{f}(\eta)]^2,
\end{equation} and let $\mathcal{D}$ be the quadratic form defined by 
\begin{align*}
\mathcal{D}(\sqrt{f}) &:= \mathcal{D}_1(\sqrt{f}) + \mathcal{D}_2(\sqrt{f}) + \mathcal{D}_3(\sqrt{f}),
\end{align*} where 
\begin{align} \nonumber
\mathcal{D}_1(\sqrt{f}) &:= N\;\mathbb{E}_{\nu_{\rho(\cdot)}^N}\left[\sum_{\xi\in {\{0,1\}^{\Sigma^-_l}}}R^-\left(\Pi_{1,l} \eta, \xi \right)E_\eta(\sqrt{f},\xi||\Pi_{l+1,N-1} \eta) \right];\\ \label{Dirichlet_bulk_contribution}
\mathcal{D}_2(\sqrt{f}) &:= N^2\;\mathbb{E}_{\nu_{\rho(\cdot)}^N}\left[\sum_{x\in\Lambda_N^{\circ}}E_\eta(\sqrt{f},\sigma^{x,x+\epsilon_N}\eta)\right];\\
\nonumber
\mathcal{D}_3(\sqrt{f}) &:= N\;\mathbb{E}_{\nu_{\rho(\cdot)}^N}\left[\sum_{\xi\in {\{0,1\}^{\Sigma^+_l}}} R^+\left(\Pi_{N-l,N-1} \eta, \xi \right)E_\eta(\sqrt{f},\Pi_{1,N-l-1} \eta || \xi)]\right].
\end{align}
Then,  {there exists positive constants $C^{bulk}$ and $C^{bd}$ such that
\begin{equation} \label{bound_dirichlet_total}
 {-}\mathfrak{D}(\sqrt{f})\leq - \frac{1}{4}\mathcal{D}(\sqrt{f})+  {(C^{bulk}+ C^{bd})N},
\end{equation} for every density $f$ with respect to $\nu^N_{\rho(\cdot)}$.}
\end{lemma}

\begin{proof}
By definition,
\begin{equation*}
    -\mathfrak{D}(\sqrt{f})=\int_{\Omega_N} \left([\mathscr{L}^{-}_N\sqrt{f}(\eta) + \mathscr{L}_N^{bulk}\sqrt{f}(\eta) + \mathscr{L}^{+}_N\sqrt{f}(\eta)]\sqrt{f}(\eta) \right)\;d\nu_{\rho(\cdot)}^N(\eta).
\end{equation*}

For the left boundary, by completing the square, we get that
\begin{align} \nonumber
&\int_{\Omega_N}(\mathscr{L}^{-}_N\sqrt{f})(\eta)\sqrt{f}(\eta)\;d\nu_{\rho(\cdot)}^N(\eta) \\ 
\label{first_term_bound}
&=-\frac{N}{2}\int_{\Omega_N}\sum_{\xi\in {\{0,1\}^{\Sigma^-_l}}}R^-\left(\Pi_{1,l} \eta, \xi \right)E_\eta(\sqrt{f},\xi||\Pi_{l+1,N-1} \eta)\;d\nu_{\rho(\cdot)}^N(\eta)
\\ \label{second_term_bound}
&+\frac{N}{2} \int_{\Omega_N}\sum_{\xi\in {\{0,1\}^{\Sigma^-_l}}}R^-\left(\Pi_{1,l} \eta, \xi \right)[f(\xi||\Pi_{l+1,N-1} \eta)-f(\eta)]\;d\nu_{\rho(\cdot)}^N(\eta).
\end{align}
 {The term in \eqref{first_term_bound} is equal to $-\frac{1}{2} \mathcal{D}_1(\sqrt{f})$}, so we only need to estimate \eqref{second_term_bound}. To do that, we will perform a change of variables $\eta \mapsto \zeta = \xi||\Pi_{l+1,N-1} \eta$. We start by computing
\begin{equation} \label{Radon_Nikodym_to_compute}
    \frac{d\nu_{\rho(\cdot)}^N(\zeta)}{d\nu_{\rho(\cdot)}^N(\eta)}
\end{equation} using the definition of $\nu_{\rho(\cdot)}^N(\zeta)$, which we recall is given by
\begin{equation*}
\nu_{\rho(\cdot)}^N(\zeta)=\prod_{x\in\Lambda_N}[\rho(x)]^{\zeta(x)}[1-\rho(x)]^{1-\zeta(x)}.
\end{equation*} Let $\mathscr{S}=\{x\in\Lambda_N\;:\;\zeta(x)\neq\eta(x)\} \subset  {\Sigma_l^-}$. Then,
\begin{align*}
\nu_{\rho(\cdot)}^N(\zeta)&=\prod_{x\in\Lambda_N/\mathscr{S}}[\rho(x)]^{\eta(x)}[1-\rho(x)]^{1-\eta(x)}\prod_{x\in\mathscr{S}}[\rho(x)]^{\zeta(x)}[1-\rho(x)]^{1-\zeta(x)}\\
&=\nu_{\rho(\cdot)}^N(\eta)\prod_{x\in\mathscr{S}}[\rho(x)]^{\zeta(x)}[1-\rho(x)]^{1-\zeta(x)}\prod_{x\in\mathscr{S}}[\rho(x)]^{-\eta(x)}[1-\rho(x)]^{\eta(x)-1},
\end{align*}
and, since $\zeta(x)-\eta(x)=1-2\eta(x)$ for every $x \in  {\mathscr{S}}$, we can rewrite the expression above to get 
\begin{equation*}
\nu_{\rho(\cdot)}^N(\zeta)=\nu_{\rho(\cdot)}^N(\eta)\prod_{x\in\mathscr{S}}\left(\frac{\rho(x)}{1-\rho(x)}\right)^{1-2\eta(x)} \Longleftrightarrow \nu_{\rho(\cdot)}^N(\eta)=\nu_{\rho(\cdot)}^N(\zeta)\prod_{x\in\mathscr{S}}\left(\frac{\rho(x)}{1-\rho(x)}\right)^{1-2\zeta(x)}.
\end{equation*}
Therefore, using the change of variables $\eta \mapsto \zeta = \xi||\Pi_{l+1,N-1}\eta$, we obtain
\begin{align*}
&\int_{\Omega_N}\sum_{\xi\in {\{0,1\}^{\Sigma^-_l}}}R^-\left(\Pi_{1,l} \eta, \xi \right) f(\xi||\Pi_{l+1,N-1} \eta)\;d\nu_{\rho(\cdot)}^N(\eta)\\
=&\int_{\Omega_N}\sum_{\mu\in {\{0,1\}^{\Sigma^-_l}}}R^-(\mu,\Pi_{1,l}\zeta)f(\zeta)\prod_{x\in\mathscr{S}}\left(\frac{\rho(x)}{1-\rho(x)}\right)^{1-2\zeta(x)}d\nu^N_{\rho(\cdot)}(\zeta).
\end{align*}
Hence, 
\begin{align*}
&N\int_{\Omega_N}\sum_{\xi\in {\{0,1\}^{\Sigma^-_l}}}R^-\left(\Pi_{1,l} \eta, \xi \right)[f(\xi||\Pi_{l+1,N-1} \eta)-f(\eta)]\;d\nu_{\rho(\cdot)}^N(\eta)\\
=&N\int_{\Omega_N}\sum_{\xi\in {\{0,1\}^{\Sigma^-_l}}}\left[R^-(\xi,\Pi_{1,l}\eta)\prod_{x\in\mathscr{S}}\left(\frac{\rho(x)}{1-\rho(x)}\right)^{1-2\eta(x)}-R^-(\Pi_{1,l}\eta,\xi)\right]f(\eta)\;d\nu_{\rho(\cdot)}^N(\eta).
\end{align*}

Since $\mathscr{S} \subset  {\Sigma^-_l}$, the bracket term is composed by finitely many rates of change between configurations in $ {\{0,1\}^{\Sigma^-_l}}$, $\rho(\cdot)$ is bounded away from zero and one and $f$ is a density, then, the last display is bounded by  {$N$ times} a constant depending only on $l$, $||\rho||_\infty$ and $R_-$,  {that we will denote by $C^{bd}_-$.}

Summarizing,
\begin{equation} \label{bound_dirichlet_lb}
    \int_{\Omega_N}(\mathscr{L}^{-}_N\sqrt{f})(\eta)\sqrt{f}(\eta)\;d\nu_{\rho(\cdot)}^N(\eta) \leq -\frac{1}{2}\mathcal{D}_1 +  C_-^{bd}. {N},
\end{equation}
and we emphasize that the constant $C_-^{bd}$ does not depend on $N$. 

Analogously,
\begin{equation} \label{bound_dirichlet_rb}
    \int_{\Omega_N}(\mathscr{L}^{+}_N\sqrt{f})(\eta)\sqrt{f}(\eta)\;d\nu_{\rho(\cdot)}^N(\eta) \leq -\frac{1}{2}\mathcal{D}_3 +  C_+^{bd}. {N},
\end{equation}
and,  { by treating the bulk term as usually handled, see for example \cite{FGS23} proof of 
Lemma 4.1, we get that}
 \begin{equation} \label{bound_dirichlet_bulk}
    \int_{\Omega_N}(\mathscr{L}_N^{bulk}\sqrt{f})(\eta)\sqrt{f}(\eta)\;d\nu_{\rho(\cdot)}^N(\eta) \leq -\frac{1}{4}\mathcal{D}_2+C^{bulk}.N,
\end{equation} with both constants, $C_+^{bd}$ and $C^{bulk}$ independent of $N$,  {and this completes the proof.}
\end{proof}

\subsection{Proof of the Replacement Lemma and corollaries} \label{replacement_lemmas_proof_section}
In this section, we prove the Boundary Replacement Lemma and corollaries needed in order to close \textrm{Dynkin's martingale} $ {\mathcal{M}_N^H(t)}$ in terms of the empirical measure and characterize the limit point.

\subsubsection{Proof of Corollary \ref{corollary_repl_lemma_ave}} 

Here we will only give the proof of \eqref{eq_1_corollary_RL1} and leave to the reader the proof of \eqref{eq_2_corollary_RL1} since it is completely analogous.

If $H(s,t) = 0$ for every $(s,t) \in  {\Omega_T}$, then the result is trivial. Let us take a function $H \in C^{1,2}( {\Omega_T}) \setminus\{0\}$ and let $\delta > 0$. We now follow the same strategy as in Section 6.1. of \cite{FGNslow}.  {Denote, for every $s \in [0,T]$ and $\eta \in \Omega_N$,}
\begin{align*}
    f_H^N(s,\eta,\epsilon):= \eta_{sN^2}(\epsilon_N) N [ H_s(\epsilon_N) - H_s(0)] - \pi^N_s \ast \iota_\epsilon(\epsilon_N)\partial_u H_s(0).
\end{align*} We want to prove that, for any $\delta > 0$, 
\begin{equation*}
    \lim_{\epsilon \to 0} \lim_{N \to \infty} \mathbb{P}_{\nu_N} \left( \sup_{0 \leq t \leq T} \Big | \int_0^t f_H^N(s,\eta,\epsilon)ds\Big| > \delta \right) = 0.
\end{equation*} Remark that we have a supremum inside a probability, and so to prove the result, we will start  {by reducing the problem to prove} that, for every $\delta > 0$, for every function $H \in C^{1,2}( {\Omega_T})$ and $t \in [0,T]$,
\begin{equation} \label{what_to_rpove_here}
\lim_{\epsilon\to 0}\lim_{N\to+\infty}\mathbb{P}_{\nu_N} \left( \Big | \int_0^t f_H^N(s,\eta,\epsilon) ds\Big| > \delta \right) =0.
\end{equation}

To do that, we take a partition $t_0 = 0 <
t_1 < \dots < t_j = T$ of $[0,T]$ with the size of the mesh bounded by
\begin{equation} \label{constant}
   \delta [ {  2 }(|| \partial^2_u H_s||_\infty + || \partial_u H_s||_\infty)]^{-1},
\end{equation} where $\displaystyle||G_s||_{\infty}:= \sup_{0 \leq s \leq T} \sup_{ x \in [0,1]} |G_s(x)|$, for $G \in C( {\Omega_T})$. Then
\begin{align*}
\sup_{0\leq t \leq T} \Big | \int_0^t f_H^N(s,\eta,\epsilon) ds\Big|
\leq \sup_{0 \leq t \leq T} \left(\Big | \int_{0}^{t_{K_t}} f_H^N(s,\eta,\epsilon) ds\Big| + \Big | \int_{t_{K_t}}^t f_H^N(s,\eta,\epsilon) ds\Big| \right),
\end{align*}
where, for every $t \in [0,T)$, $t_{K_t}$ is the unique point of the mesh for which $t_{K_t} \leq t < t_{K_t+1}$ and, for $t = T$, $t_{K_t} = t_{j-1}$. Then, for every $t \in [0,T]$,
\begin{align*}
    \Big | \int_{0}^{t_{K_t}} f_H^N(s,\eta,\epsilon) ds\Big| = \Big| \int_{0}^{t_{i}} f_H^N(s,\eta,\epsilon) ds\Big|,
\end{align*}
for some $i \in \{0,\dots,j\}$, and
\begin{align} \label{remaining}
    \Big | \int_{t_{K_t}}^{t} f_H^N(s,\eta,\epsilon) ds\Big| \lsim |t-t_{K_t}| (|| \partial_u H_s||_\infty + || \partial^2_u H_s||_\infty) \leq \delta/2,
\end{align}  { where the last estimate is a consequence of the choice of $\delta$ given in \eqref{constant}},  {the fact that the occupation variables are bounded and $H \in C^{1,2}(\Omega_T)$.} Thus,
\begin{align*}
    &\mathbb{P}_{\nu_N} \left( \sup_{0 \leq t \leq T} \Big | \int_0^t f_H^N(s,\eta,\epsilon)ds\Big| > \delta \right) \\
    &\leq \sum_{i=0}^j\mathbb{P}_{\nu_N} \left( \Big | \int_0^{t_i} f_H^N(s,\eta,\epsilon)ds\Big| > \delta/2 \right) + \underbrace{\mathbb{P}_{\nu_N} \left( \sup_{0 \leq t \leq T} \Big | \int_{t_{K_t}}^t f_H^N(s,\eta,\epsilon)ds\Big| > \delta/2 \right)}_{ = 0 \textrm{ because of } \eqref{remaining}}.
\end{align*} Since $j$ only depends on $\delta$ and  {$H$} and is finite, to prove \eqref{eq_1_corollary_RL1}, it is enough to show  {\eqref{what_to_rpove_here}, as we wanted.} Then, for every $\delta > 0$,
\begin{align} \label{RHSterm}
\mathbb{P}_{\nu_N} \left(\Big | \int_0^t f_H^N(s,\eta,\epsilon) ds\Big| > \delta \right) 
&\leq \mathbb{P}_{\nu_N} \left(\Big | \int_0^t [\eta_{sN^2}(\epsilon_N) - \pi^N_s \ast \iota_\epsilon(\epsilon_N)] \partial_u  H_s(0) ds \Big | > \frac{\delta}{2}\right) \\ \label{is_zero}
&+ \mathbb{P}_{\nu_N} \left(  {\int_0^t \frac{||\partial^2_u H_s||_\infty}{2N} ds} > \frac{\delta}{2} \right).
\end{align}
For every fixed $\delta >0$, we can find $N$ sufficiently large such that, $ \int_0^t\frac{||\partial^2_u  {H}_s||_\infty}{2N} ds < \frac{\delta}{2}$, therefore the term in \eqref{is_zero} will vanish as we take $N$ to infinity. Moreover, since $s \mapsto \partial_u   {H}_s(0)$ is a continuous function in $[0,T]$, using Markov's inequality and Lemma \ref{repl_lemma_ave},  {we obtain that the term on the RHS of the inequality in \eqref{RHSterm}} goes to zero when taking the limit as $N$ goes to infinity and then $\epsilon$ to zero. This finishes the proof of the result.

\subsubsection{Proof of Lemma \ref{repl_lemma_2} and its Corollary \ref{corollary_repl_lemma_2}} 

\begin{proof}[Proof of Lemma \ref{repl_lemma_2}]
Observe that
\begin{align} \label{boundbound}
    \langle h_{-}(\eta),  {H}_s\rangle &=  \sum_{j \in  {\Sigma^-_l}} \sum_{\xi\in {\{0,1\}^{\Sigma^-_l}}}R^{-}(\Pi_{1,l}\eta,\xi) (\xi(j)-\eta(j))[  {H}_s(j) -  {H}_s(0)] \\ \nonumber
    &+ \sum_{j \in  {\Sigma^-_l}} \sum_{\xi\in {\{0,1\}^{\Sigma^-_l}}}R^{-}(\Pi_{1,l}\eta,\xi) (\xi(j)-\eta(j))  {H}_s(0).
\end{align}
Since $ {H}_\cdot$ is  {$C^{1,2}( {\Omega_T})$}, in particular, $ {H}_s \in C^1[0,1]$ for every $s \in [0,T]$, then,  {the term on the RHS of \eqref{boundbound}} can be bounded by a constant $C = C(R^-,l)$ times $1/N$.  {So to prove our result, it is enough to estimate
\begin{align}
    \label{almost_what_to_prove_RL2}
    \mathbb{E}_{\nu_N} \left[\Big | \int_0^t V_-(\epsilon,\eta_s)  {H}_s(0) ds \Big | \right],
\end{align}
where
\begin{equation*}
    V_-(\epsilon,\eta_s) := h_-(\eta_s) - F_-(\overrightarrow{\eta_s}^{\lfloor \epsilon N \rfloor}(\epsilon_N)).
\end{equation*}}

Since $h_- \in L^2(\nu^N_{\rho(\cdot)})$, by the proof of Lemma \ref{lemma_basis_L^2}, meaning the identity in \eqref{equation_f_in_prod_eta}, we have, for every $\eta \in  {\{0,1\}^{\Sigma^-_l}}$, that $h_-$ can be written in terms of linear combinations of products of $\eta(1), \dots, \eta(l)$ as
\begin{equation*}
    h_-(\eta) = \sum_{\xi \in  {\{0,1\}^{\Sigma^-_l}}} \sum_{B \subset A^c_\xi}  (-1)^{|B|} h_-(\xi) D(\eta,B \cup A_\xi),
\end{equation*}
where  {$A_\xi = \{x \in \Sigma^-_l \ | \ \xi(x) = 1\}$ as it was defined in \eqref{def_A_xi_and_complement} and} $D(\gamma, C) = \prod_{x \in C} \gamma(x)$ for every $\gamma \in  {\{0,1\}^{\Sigma^-_l}}$ and $C \subset  {\Sigma^-_l}$. Thus,
\begin{align} \label{v_eta_def}
    V_-(\epsilon,\eta_s) 
    &=\sum_{\xi \in  {\{0,1\}^{\Sigma^-_l}}} \sum_{B \subset A^c_\xi}  (-1)^{|B|} h_-(\xi) \left[D(\eta,B \cup A_\xi) - \left[\overrightarrow{\eta}^{\lfloor \epsilon N \rfloor}(\epsilon_N) \right]^{|B \cup A_\xi|}\right]
\end{align} where to obtain the second equality we used the product form of $\nu^N_{\overrightarrow{\eta}^{\lfloor \epsilon N \rfloor}(\epsilon_N)}$.

By the usual technique, i.e. using the entropy inequality, the fact that $e^{|x|} \leq e^x + e^{-x}$, Jensen's inequality and Feynman-Kac's formula, we obtain that
\begin{align*}
    &\mathbb{E}_{\nu_N} \left[\Big | \int_0^t V_-(\epsilon,\eta_s)  {H}_s(0) ds \Big | \right] \leq \frac{ {H(\nu_N | \nu^N_{\rho(\cdot)})} + \log(2)}{B N} \\
    &\quad + \max\left\{\int_0^t | {H}_s(0)|ds, 1 \right\} \sup_{f \textrm{density}} \Big\{ \pm \langle V_-(\epsilon,\cdot), f \rangle_{\nu^N_{\rho(\cdot)}} + \frac{1}{BN} \langle \mathscr{L}_N \sqrt f, \sqrt f \rangle_{\nu^N_{\rho(\cdot)}} \Big\}.
\end{align*}
Since we already estimated in Lemma \ref{lemma_estimate_dirichlet_form} the term $\langle \mathscr{L}_N \sqrt f, \sqrt f \rangle_{\nu^N_{\rho(\cdot)}}$, for every density $f$, it is enough to show that $\langle \pm V_-(\epsilon,\cdot), f \rangle_{\nu^N_{\rho(\cdot)}}$ goes to zero as $N$ goes to infinity and then $\epsilon$ go to zero.  {In what follows, we only give the details for the term $\langle V_-(\epsilon,\cdot), f \rangle_{\nu^N_{\rho(\cdot)}}$ since the proof with the minus sign is analogous.}

 {We would like now to reduce the bound of the last display to the same type of arguments used in Lemma \ref{repl_lemma_ave}. So, remark that, for every $\xi \in  {\{0,1\}^{\Sigma^-_l}}$, denoting $\{x_1,\dots,x_{|B \cup A_\xi|}\} := B \cup A_\xi$, with $x_j \neq x_i$ if $i \neq j$ and $i,j \in \{1,\dots |B \cup A_\xi|\}$, by summing and subtracting some terms, we have that
\begin{align*}
    D(\eta,B \cup A_\xi) - \left[\overrightarrow{\eta}^{\lfloor \epsilon N \rfloor}(\epsilon_N) \right]^{|B \cup A_\xi|}= \sum_{j=1}^{|B \cup A_\xi|} \mathfrak{g}(\eta,B \cup A_\xi, j , \epsilon N)\left[ \eta(x_j) - \overrightarrow{\eta}^{\lfloor \epsilon N \rfloor}(\epsilon_N) \right],
\end{align*} where
\begin{align}
\mathfrak{g}(\eta,B \cup A_\xi, j , \epsilon N) :=  D(\eta,B \cup A_\xi\setminus \cup_{i=1}^j\{x_i\})\left[\overrightarrow{\eta}^{\lfloor \epsilon N \rfloor}(\epsilon_N) \right]^{j-1}.
\end{align}}Thus, using \eqref{v_eta_def} and the previous observation, we obtain that 
\begin{align} \label{V_f_expression_in_term_of_Is}
\langle V_-(\epsilon,\cdot), f \rangle_{\nu^N_{\rho(\cdot)}} &= \sum_{\xi \in  {\{0,1\}^{\Sigma^-_l}}} \sum_{B \subset A^c_\xi}  (-1)^{|B|} h_-(\xi) \sum_{j = 1}^{|B \cup A_\xi|} I^\epsilon_j,
\end{align}
where, fixing $\xi \in  {\{0,1\}^{\Sigma^-_l}}$ and $B \subset A^c_\xi$, for each $j \in \{1,\dots,|B \cup A_\xi| \}$,
\begin{align*}
    I^\epsilon_j &:= \int_{\Omega_N} \mathfrak{g}(\eta,B \cup A_\xi, j , \epsilon N) \left[ \eta(x_j) - \overrightarrow{\eta}^{\lfloor \epsilon N \rfloor}(\epsilon_N) \right] f(\eta) d\nu^N_{\rho(\cdot)}(\eta).
\end{align*} We are now reduced to estimate each of the integrals $I^\epsilon_j$,  {that now are only comparing one $\eta(x_j)$ with the mean number of particles on a box around $\epsilon_N$ given by $\overrightarrow{\eta}^{\lfloor \epsilon N \rfloor}(\epsilon_N)$.}

 {To simplify notation, let us set $z = \epsilon_N$.} Recall that
\begin{equation*}
    \overrightarrow{\eta}^{\lfloor \epsilon N\rfloor}(z) = \frac{1}{2\epsilon N+1}\sum_{|x-z|\leq\epsilon} \eta(x).
\end{equation*}  {Thus, we can decompose $I^\epsilon_j$ as $I^\epsilon_j = I^{\epsilon,+}_j + I^{\epsilon,-}_j$, where}
\begin{align}
\label{1_int_I_1}
    I^{\epsilon,+}_j := \sum_{w = 0}^{\epsilon N} \int_{\Omega_N} \frac{\mathfrak{g}(\eta,B \cup A_\xi, j , \epsilon N)}{2 \epsilon N +1} [\eta(x_j) - \eta(z + w/N)] \mathds{1}_{\{z + w/N \notin B \cup A_\xi \setminus \cup_{i=1}^j\{x_i\}\}} f(\eta) d\nu^N_{\rho(\cdot)}(\eta)
\end{align} and
\begin{align}\label{2_int_I_1}
    I^{\epsilon,-}_j := \sum_{w = 0}^{\epsilon N} \int_{\Omega_N} \frac{\mathfrak{g}(\eta,B \cup A_\xi, j , \epsilon N)}{2 \epsilon N +1} [\eta(x_j) - \eta(z + w/N)]\mathds{1}_{\{z + w/N \in B \cup A_\xi \setminus \cup_{i=1}^j\{x_i\}\}} f(\eta) d\nu^N_{\rho(\cdot)}(\eta).
\end{align} Since $|B \cup A_\xi \setminus\cup_{i=1}^j\{x_i\}| \leq l$, the sum in \eqref{2_int_I_1} has at most $l$ terms. Moreover, because the occupation variables are bounded, we get that $I^{\epsilon,-}_j$ is of order $\frac{l}{\epsilon N}$.  {To control $I^{\epsilon,+}_j$, observe that,} for each $w \in \{0, \dots, \epsilon N\}$ such that $z +w/N \notin B \cup A_\xi \setminus \cup_{i=1}^j\{x_i\}$,  {by creating telescopic sums, }we have that

\begin{align} \nonumber
    &\int_{\Omega_N} \mathfrak{g}(\eta,B \cup A_\xi, j , \epsilon N)[\eta(x_j) - \eta(z + w/N)] f(\eta) d\nu^N_{\rho(\cdot)}(\eta)\\
    \label{estimate_1_RL2}
    =& \sum_{k=0}^{w - N x_j}\int_{\Omega_N} \mathfrak{g}(\eta,B \cup A_\xi, j , \epsilon N) [\eta(w_k) - \eta(w_k + \epsilon_N)] f(\eta) \mathbbm{1}_{\{x_j < z + w/N\}} d\nu^N_{\rho(\cdot)}(\eta)\\ \label{estimate_2_RL2}
    -& \sum_{m=0}^{N x_j - w - 2} \int_{\Omega_N} \mathfrak{g}(\eta,B \cup A_\xi, j , \epsilon N) [\eta(v_m) - \eta(v_m + \epsilon_N)] f(\eta) \mathbbm{1}_{\{x_j > z + w/N\}} d\nu^N_{\rho(\cdot)}(\eta),
\end{align}
where, for each $k \in \{0,\dots,w-N x_j\}$ if $x_j < z + w/N$  {(resp. for each $m \in \{0,\dots,Nx_j-w\}$ if $x_j > z + w/N$)}, 
\begin{equation}
w_k := x_j+k/N \quad  { ( \textrm{resp.} \ v_m := z+ (w+m)/N )}.
\end{equation}
Remark that $w_{k+1} = w_k + \epsilon_N$ (resp. $v_{m+1} = v_m + \epsilon_N$). 

 {Now, we estimate \eqref{estimate_1_RL2} and remark that, to bound \eqref{estimate_2_RL2}, the proof is analogous. First, we remark that,} for every $w \in \{0,\dots, \epsilon N\}$ with $z + w/N \notin B \cup A_\xi \setminus (\cup_{i=1}^{j}\{x_i\})$, and every $k \in \{0,\dots,N x_j - w\}$ whenever $x_j < z + w/N$, we have that
\begin{align*}
    \overrightarrow{\sigma^{w_k, w_k+ \epsilon_N}\eta}^{\lfloor \epsilon N \rfloor}(\epsilon_N) &= \overrightarrow{\eta}^{\lfloor \epsilon N \rfloor}(\epsilon_N) + \frac{1}{2\epsilon N+1} [\eta(w_k+\epsilon_N) - \eta(w_k)] \mathbbm{1}_{|w_k - \epsilon_N| \leq \epsilon N} \mathbbm{1}_{|w_k| > \epsilon N}\\
    &+ \frac{1}{2\epsilon N+1} [\eta(w_k) - \eta(w_k+\epsilon_N)] \mathbbm{1}_{|w_k| \leq \epsilon N} \mathbbm{1}_{|w_k - \epsilon_N| > \epsilon N}.
\end{align*}

Denoting, for every $s \in \Lambda_N \setminus \{\tau_N\}$,
\begin{equation}
    \tau^{s}(\eta) := \eta(s) - \eta(s + \epsilon_N),
\end{equation} 
by changing variables $\eta \mapsto \xi = \sigma^{w_k,w_k+\epsilon_N}\eta$ in 
\begin{align*}
    \int_{\Omega_N} \mathfrak{g}(\eta,B \cup A_\xi, j , \epsilon N) \tau^{w_k}(\eta) f(\eta) \mathbbm{1}_{\{x_j < z + w/N\}} d\nu^N_{\rho(\cdot)}(\eta),
\end{align*} and recalling the binomial expansion, we get that
\begin{align*}
    \int_{\Omega_N} 2 \mathfrak{g}(\eta,B \cup A_\xi, j , \epsilon N) \tau^{w_k}(\eta) f(\eta) \mathbbm{1}_{\{x_j < z + w/N\}} d\nu^N_{\rho(\cdot)}(\eta),
\end{align*}
\begin{itemize}
    \item if $|[B \cup A_\xi \setminus (\cup_{i=1}^j \{x_i\}) ]\cap \{w_k,w_k+\epsilon_N\}| \in \{0,2\}$, is equal to
\begin{align*}
&\int_{\Omega_N} \mathfrak{g}(\eta,B \cup A_\xi, j , \epsilon N) \tau^{w_k}(\eta) \left\{ f(\eta) - f(\sigma^{w_k, w_k+ \epsilon_N}\eta) a_{w_k+\epsilon_N,w_k}(\eta) \right\} d\nu^N_{\rho(\cdot)}(\eta)\\
&+ \sum_{p=1}^{j-1} \binom{j-1}{p} \int_{\Omega_N} (2\epsilon N+1)\left[D(\eta,B \cup A_\xi\setminus \cup_{i=1}^j\{x_i\}) \left[\overrightarrow{\eta}^{\lfloor \epsilon N \rfloor}(z) \right]^{j-1-p} \mathbbm{1}_{\{x_j < z + w/N\}} \right. \\
&\hspace{2.5cm}\left[\frac{\tau^{w_k}(\eta)}{2\epsilon N+1}\right]^{p+1}\left( \left[-1\right]^{p+1}\mathbbm{1}_{|w_k - \epsilon_N| \leq \epsilon N} \mathbbm{1}_{|w_k| > \epsilon N} - \mathbbm{1}_{|w_k - \epsilon_N|> \epsilon N} \mathbbm{1}_{|w_k| \leq \epsilon N}  \right)\\
& \left.\hspace{7.3cm} f(\sigma^{w_k, w_k+ \epsilon_N}\eta) a_{w_k+\epsilon_N,w_k}(\eta) \right]d\nu^N_{\rho(\cdot)}(\eta);
\end{align*}

\item if $w_k \in B \cup A_\xi$ but $ w_k+\epsilon_N \notin B \cup A_\xi$, is equal to
\begin{align*}
&\displaystyle \int_{\Omega_N} \left[\overrightarrow{\eta}^{\lfloor \epsilon N \rfloor}(z)\right]^{j-1}D(\eta,B \cup A_\xi\setminus [\cup_{i=1}^j\{x_i\} \cup \{w_k\}]) \tau^{w_k}(\eta) \mathbbm{1}_{\{x_j < z + w/N\}} \\
&\hspace{3.5cm} \left\{\eta(w_k) f(\eta) - \eta(w_k+ \epsilon_N) f(\sigma^{w_k, w_k+ \epsilon_N}\eta) a_{w_k+\epsilon_N,w_k}(\eta) \right\}  d\nu^N_{\rho(\cdot)}(\eta)\\
&\displaystyle 
+ \sum_{p=1}^{j-1} \binom{j-1}{p} \int_{\Omega_N} (2\epsilon N+1)\left[D(\eta,B \cup A_\xi\setminus [\cup_{i=1}^j\{x_i\} \cup \{w_k\}]) \left[\overrightarrow{\eta}^{\lfloor \epsilon N \rfloor}(z) \right]^{j-1-p}\right. \\
&\hspace{0.8cm}\mathbbm{1}_{\{x_j < z + w/N\}}\left[\frac{\tau^{w_k}(\eta)}{2\epsilon N+1}\right]^{p+1}\left( \left[-1\right]^{p+1} \mathbbm{1}_{|w_k - \epsilon_N| \leq \epsilon N} \mathbbm{1}_{|w_k| > \epsilon N} 
- \mathbbm{1}_{|w_k - \epsilon_N| > \epsilon N} \mathbbm{1}_{|w_k| \leq \epsilon N} \right)\\
&\hspace{6cm} \left.\eta(w_k+ \epsilon_N) f(\sigma^{w_k, w_k+ \epsilon_N}\eta) a_{w_k+\epsilon_N,w_k}(\eta) \right]d\nu^N_{\rho(\cdot)}(\eta);
\end{align*}

\item analogously to the previous case, if $w_k + \epsilon_N \in B \cup A_\xi$ but $ w_k \notin B \cup A_\xi$, is equal to
\begin{align*}
&\displaystyle \int_{\Omega_N} \left[\overrightarrow{\eta}^{\lfloor \epsilon N \rfloor}(z)\right]^{j-1} D(\eta,B \cup A_\xi\setminus [\cup_{i=1}^j\{x_i\} \cup \{w_k+\epsilon_N\}]) \tau^{w_k}(\eta) \mathbbm{1}_{\{x_j < z + w/N\}} \\
&\hspace{3.5cm} \left\{\eta(w_k+\epsilon_N) f(\eta) - \eta(w_k) f(\sigma^{w_k, w_k+ \epsilon_N}\eta) a_{w_k+\epsilon_N,w_k}(\eta) \right\}  d\nu^N_{\rho(\cdot)}(\eta)\\
&\displaystyle 
+ \sum_{p=1}^{j-1} \binom{j-1}{p} \int_{\Omega_N} (2\epsilon N+1) \left[D(\eta,B \cup A_\xi\setminus [\cup_{i=1}^j\{x_i\} \cup \{w_k+\epsilon_N\}]) \left[\overrightarrow{\eta}^{\lfloor \epsilon N \rfloor}(z) \right]^{j-1-p} \right.\\
&\hspace{1cm}\mathbbm{1}_{\{x_j < z + w/N\}} \left[\frac{\tau^{w_k}(\eta)}{2\epsilon N+1}\right]^{p+1}\left( \left[-1\right]^{p+1}\mathbbm{1}_{|w_k - \epsilon_N| \leq \epsilon N} \mathbbm{1}_{|w_k| > \epsilon N} 
- \mathbbm{1}_{|w_k - \epsilon_N| > \epsilon N} \mathbbm{1}_{|w_k| \leq \epsilon N}\right) \\
& \hspace{7cm}\left.\eta(w_k) f(\sigma^{w_k, w_k+ \epsilon_N}\eta) a_{w_k+\epsilon_N,w_k}(\eta) \right] d\nu^N_{\rho(\cdot)}(\eta).
\end{align*}
\end{itemize}

Using the identity $x^2-y^2 = (\sqrt{x} - \sqrt{y})(\sqrt{x} + \sqrt{y})$ for $x,y \geq 0$ and Young's inequality:
\begin{itemize}
    \item if $|[B \cup A_\xi \setminus (\cup_{i=1}^j \{x_i\}) ]\cap \{w_k,w_k+\epsilon_N\}| \in \{0,2\}$, we get that the expression we wrote previously can be bounded from above by
    \begin{align*}
    &\frac{1}{2AN^2}\int_{\Omega_N} N^2\left\{ \sqrt{f}(\sigma^{w_k, w_k+ \epsilon_N}\eta) - \sqrt{f}(\eta)\right\}^2 d\nu^N_{\rho(\cdot)}(\eta) \\
    &+ A\int_{\Omega_N} \left[\mathfrak{g}(\eta,B \cup A_\xi, j , \epsilon N) \tau^{w_k}(\eta)\right]^2 \left\{ f(\eta) + f(\sigma^{w_k, w_k+ \epsilon_N}\eta)\right\} d\nu^N_{\rho(\cdot)}(\eta) \\ 
    &+\int_{\Omega_N} \mathfrak{g}(\eta,B \cup A_\xi, j , \epsilon N) \tau^{w_k}(\eta) f(\sigma^{w_k, w_k+ \epsilon_N}\eta)\left[1 - a_{w_k+\epsilon_N,w_k}(\eta)\right] d\nu^N_{\rho(\cdot)}(\eta)\\
    &+ \sum_{p=1}^{j-1} \binom{j-1}{p} \int_{\Omega_N} (2\epsilon N+1)\left[D(\eta,B \cup A_\xi\setminus \cup_{i=1}^j\{x_i\}) \left[\overrightarrow{\eta}^{\lfloor \epsilon N \rfloor}(z) \right]^{j-1-p} \mathbbm{1}_{\{x_j < z + w/N\}} \right. \\
    &\hspace{2.5cm}\left( \mathbbm{1}_{|w_k - \epsilon_N| \leq \epsilon N} \mathbbm{1}_{|w_k| > \epsilon N} + \mathbbm{1}_{|w_k - \epsilon_N|> \epsilon N} \mathbbm{1}_{|w_k| \leq \epsilon N}  \right) \left[\frac{\tau^{w_k}(\eta)}{2\epsilon N+1}\right]^{p+1}\\
    & \left.\hspace{6.5cm} f(\sigma^{w_k, w_k+ \epsilon_N}\eta) a_{w_k+\epsilon_N,w_k}(\eta) \right]d\nu^N_{\rho(\cdot)}(\eta);
    \end{align*}

    \item if $w_k \in B \cup A_\xi$ but $ w_k+\epsilon_N \notin B \cup A_\xi$, the expression we wrote previously can be bounded from above by
    \begin{align*}
    &\frac{1}{2AN^2}\int_{\Omega_N} N^2\left\{ \sqrt{f}(\sigma^{w_k, w_k+ \epsilon_N}\eta) - \sqrt{f}(\eta)\right\}^2 d\nu^N_{\rho(\cdot)}(\eta) \\
    &+ A\int_{\Omega_N} \left(\mathfrak{g}(\eta,B \cup A_\xi, j , \epsilon N) \tau^{w_k}(\eta)\right)^2 \mathbbm{1}_{\{x_j < z + w/N\}}\left\{ f(\eta) + f(\sigma^{w_k, w_k+ \epsilon_N}\eta) \right\}  d\nu^N_{\rho(\cdot)}(\eta)\\
    &+\int_{\Omega_N} \left[\overrightarrow{\eta}^{\lfloor \epsilon N \rfloor}(z)\right]^{j-1} D(\eta,B \cup A_\xi\setminus [\cup_{i=1}^j\{x_i,\} \cup \{w_k\}]) \tau^{w_k}(\eta) \eta(w_k+\epsilon_N) \mathbbm{1}_{\{x_j < z + w/N\}}\\
    &\hspace{3.5cm} f(\sigma^{w_k, w_k+ \epsilon_N}\eta)\left[1 - a_{w_k+\epsilon_N,w_k}(\eta)\right] d\nu^N_{\rho(\cdot)}(\eta)\\
    &+\displaystyle \int_{\Omega_N} \left[\overrightarrow{\eta}^{\lfloor \epsilon N \rfloor}(z)\right]^{j-1} D(\eta,B \cup A_\xi\setminus [\cup_{i=1}^j\{x_i\} \cup \{w_k\}]) [\tau^{w_k}(\eta)]^2 \mathbbm{1}_{\{x_j < z + w/N\}} \\
    &\hspace{3.5cm} f(\sigma^{w_k, w_k+ \epsilon_N}\eta) d\nu^N_{\rho(\cdot)}(\eta)\\
    &\displaystyle + \sum_{p=1}^{j-1} \binom{j-1}{p} \int_{\Omega_N} (2\epsilon N +1)\left[D(\eta,B \cup A_\xi\setminus [\cup_{i=1}^j\{x_i\} \cup \{w_k\}]) \left[\overrightarrow{\eta}^{\lfloor \epsilon N \rfloor}(z) \right]^{j-1-p} \mathbbm{1}_{\{x_j < z + w/N\}} \right. \\
    &\hspace{2.5cm}\left( \mathbbm{1}_{|w_k - \epsilon_N| \leq \epsilon N} \mathbbm{1}_{|w_k| > \epsilon N} 
    + \mathbbm{1}_{|w_k - \epsilon_N| > \epsilon N} \mathbbm{1}_{|w_k| \leq \epsilon N} \right) \left[\frac{\tau^{w_k}(\eta)}{2\epsilon N+1}\right]^{p+1} \eta(w_k+\epsilon_N)\\
    &\hspace{3.5cm} \left. f(\sigma^{w_k, w_k+ \epsilon_N}\eta) a_{w_k+\epsilon_N,w_k}(\eta) \right]d\nu^N_{\rho(\cdot)}(\eta);
\end{align*}

    \item analogously to the previous case, if $w_k + \epsilon_N \in B \cup A_\xi$ but $ w_k \notin B \cup A_\xi$, then the expression we wrote previously can be bounded from above by
    \begin{align*}
    &\frac{1}{2AN^2}\int_{\Omega_N} N^2\left\{ \sqrt{f}(\sigma^{w_k, w_k+ \epsilon_N}\eta) - \sqrt{f}(\eta)\right\}^2 d\nu^N_{\rho(\cdot)}(\eta) \\
    &+ A\int_{\Omega_N} \left(\mathfrak{g}(\eta,B \cup A_\xi, j , \epsilon N) \tau^{w_k}(\eta)\right)^2 \mathbbm{1}_{\{x_j < z + w/N\}}\left\{ f(\eta) + f(\sigma^{w_k, w_k+ \epsilon_N}\eta) \right\}  d\nu^N_{\rho(\cdot)}(\eta)\\
    &+\int_{\Omega_N} \left[\overrightarrow{\eta}^{\lfloor \epsilon N \rfloor}(z)\right]^{j-1} D(\eta,B \cup A_\xi\setminus [\cup_{i=1}^j\{x_i\} \cup \{w_k+\epsilon_N\}])\tau^{w_k}(\eta) \eta(w_k) \mathbbm{1}_{\{x_j < z + w/N\}}\\
    &\hspace{3.5cm} f(\sigma^{w_k, w_k+ \epsilon_N}\eta)\left[1 - a_{w_k+\epsilon_N,w_k}(\eta)\right] d\nu^N_{\rho(\cdot)}(\eta)\\
    &+\displaystyle \int_{\Omega_N} \left[\overrightarrow{\eta}^{\lfloor \epsilon N \rfloor}(z)\right]^{j-1} D(\eta,B \cup A_\xi\setminus [\cup_{i=1}^j\{x_i\} \cup \{w_k + \epsilon_N\}]) [\tau^{w_k}(\eta)]^2 \mathbbm{1}_{\{x_j < z + w/N\}} \\
    &\hspace{3.5cm} f(\sigma^{w_k, w_k+ \epsilon_N}\eta) d\nu^N_{\rho(\cdot)}(\eta)\\
    &\displaystyle + \sum_{p=1}^{j-1} \binom{j-1}{p} \int_{\Omega_N} (2\epsilon N +1)\left[D(\eta,B \cup A_\xi\setminus [\cup_{i=1}^j\{x_i\} \cup \{w_k+\epsilon_N\}]) \left[\overrightarrow{\eta}^{\lfloor \epsilon N \rfloor}(z) \right]^{j-1-p} \right. \\
    &\hspace{1cm}\mathbbm{1}_{\{x_j < z + w/N\}} \left( \mathbbm{1}_{|w_k - \epsilon_N| \leq \epsilon N} \mathbbm{1}_{|w_k| > \epsilon N} 
    + \mathbbm{1}_{|w_k - \epsilon_N| > \epsilon N} \mathbbm{1}_{|w_k| \leq \epsilon N} \right) \left[\frac{\tau^{w_k}(\eta)}{2\epsilon N+1}\right]^{p+1} \eta(w_k)\\
    &\hspace{3.5cm} \left. f(\sigma^{w_k, w_k+ \epsilon_N}\eta) a_{w_k+\epsilon_N,w_k}(\eta) \right]d\nu^N_{\rho(\cdot)}(\eta).
    \end{align*}
\end{itemize}
Since $1-a_{w_k+\epsilon_N,w_k}(\eta)$ is of order $1/N$, the occupation variables are bounded by $1$ and $f$ is a density, then, in all the three cases we considered before, we have that the last expressions are bounded from above by
\begin{align*}
    &\frac{1}{2AN^2} \mathcal{D}^{w_k,w_k+\epsilon_N}(\sqrt{f}) + A C_0 + \frac{C_1}{N} \\
    &+ C_2 \mathbbm{1}_{w_k \in B \cup A_\xi} \mathbbm{1}_{w_k+\epsilon_N \notin B \cup A_\xi} + C_3 \mathbbm{1}_{w_k \notin B \cup A_\xi} \mathbbm{1}_{w_k+\epsilon_N \in B \cup A_\xi} \\
    &+ \sum_{p=1}^{j-1} \binom{j-1}{p} \left( \mathbbm{1}_{|w_k - \epsilon_N| \leq \epsilon N} \mathbbm{1}_{|w_k| > \epsilon N} 
    + \mathbbm{1}_{|w_k - \epsilon_N| > \epsilon N} \mathbbm{1}_{|w_k| \leq \epsilon N} \right) \frac{C_4}{\epsilon^{p} N ^{p}},
\end{align*} for some positive constants $C_0, C_1,  C_2, C_3$ and $C_4$. Finally, putting together the previous estimates we get, for every $j \in \{1,\dots,|B \cup A_\xi|\}$, that
\begin{align*}
    I^\epsilon_j &\leq \frac{C \cdot l}{\epsilon N} + \sum_{w = 0}^{\epsilon N} \sum_{k=0}^{w - N x_j} \mathbbm{1}_{\{x_j < z + w/N\}} \left[ \frac{1}{2\epsilon A N^3} \mathcal{D}^{w_k,w_k+\epsilon_N}(\sqrt{f}) + \frac{A C_0}{\epsilon N} + \frac{C_1}{\epsilon N^2} \right.\\
    & + \frac{C_2}{\epsilon N} \mathbbm{1}_{w_k \in B \cup A_\xi} \mathbbm{1}_{w_k+\epsilon_N \notin B \cup A_\xi} + \frac{C_3}{\epsilon N} \mathbbm{1}_{w_k \notin B \cup A_\xi} \mathbbm{1}_{w_k+\epsilon_N \in B \cup A_\xi} \\
    & \left. + \sum_{p=1}^{j-1} \binom{j-1}{p} \left( \mathbbm{1}_{|w_k - \epsilon_N| \leq \epsilon N} \mathbbm{1}_{|w_k| > \epsilon N} 
    + \mathbbm{1}_{|w_k - \epsilon_N| > \epsilon N} \mathbbm{1}_{|w_k| \leq \epsilon N} \right) \frac{C_4}{\epsilon^{p+1} N ^{p+1}} \right]\\
    &-\sum_{w = 0}^{\epsilon N}  \sum_{m=0}^{N x_j - w - 2} \int_{\Omega_N} \mathbbm{1}_{\{x_j > z + w/N\}} \left[ \frac{1}{2\epsilon A N^3} \mathcal{D}^{v_m,v_m+\epsilon_N}(\sqrt{f}) + \frac{A C_0}{\epsilon N} + \frac{C_1}{\epsilon N^2} \right.\\
    &+ \frac{C_2}{\epsilon N} \mathbbm{1}_{v_m \in B \cup A_\xi} \mathbbm{1}_{v_m+\epsilon_N \notin B \cup A_\xi} + \frac{C_3}{\epsilon N} \mathbbm{1}_{v_m \notin B \cup A_\xi} \mathbbm{1}_{v_m+\epsilon_N \in B \cup A_\xi} \\
    &\left.+ \sum_{p=1}^{j-1} \binom{j-1}{p} \left( \mathbbm{1}_{|v_m - \epsilon_N| \leq \epsilon N} \mathbbm{1}_{|v_m| > \epsilon N} 
    + \mathbbm{1}_{|v_m - \epsilon_N| > \epsilon N} \mathbbm{1}_{|v_m| \leq \epsilon N} \right) \frac{C_4}{\epsilon^{p+1} N ^{p+1}} \right],
\end{align*} for some $C > 0$. But, for every $w \in \{0,\dots,\epsilon N\}$ where $x_{j} < z + w/N$ there exists at most one value of $k_1,k_2 \in \{0,\dots, w - N x_{j}\}$ for which $\mathbbm{1}_{w_{k_1} \leq \epsilon N}\mathbbm{1}_{ w_{k_1} + \epsilon_N > \epsilon N} = 1$ and $\mathbbm{1}_{w_{k_2} > \epsilon N}\mathbbm{1}_{ w_{k_2} + \epsilon_N \leq \epsilon N} = 1$. Moreover, since $|B \cup A_\xi| \leq l$ then, for each $w \in \{0,\dots,\epsilon N\}$ such that $z + w/N \notin B \cup A_\xi \setminus[ \cup_{i=1}^j\{x_i\}]$, there exists at most a number of order $O(\epsilon l)$ of values of $k$ for which $\mathbbm{1}_{w_k + \epsilon_N \in B \cup A_\xi} \mathbbm{1}_{w_k \notin B \cup A_\xi} = 1$ or $\mathbbm{1}_{w_k + \epsilon_N \in B \cup A_\xi} \mathbbm{1}_{w_k \notin B \cup A_\xi} = 1$. Thus,
\begin{align*}
I^\epsilon_j &\leq \frac{C \cdot l}{\epsilon N} + \frac{1}{2AN^2} \mathcal{D}(\sqrt{f}) + A C_0 \epsilon N + \frac{C_1 \epsilon^2 N^2}{\epsilon N^2} + C_2\epsilon l + C_3\epsilon l + \sum_{p=1}^{j-1} \binom{j-1}{p} 2 \frac{C_4 \epsilon^2 l N}{\epsilon^{p+1} N ^{p+1}}.
\end{align*} Choosing $A = 2 \tilde{C}(l, ||h_-||_{\infty})lB/N$, we obtain that
\begin{align*}
    \langle V_-(\epsilon,\cdot), f \rangle_{\nu^N_{\rho(\cdot)}} \leq C(l) \left[ \frac{C \cdot l}{\epsilon N} + \frac{1}{4 \tilde{C}(l, ||h_-||_{\infty})lB N} \mathcal{D}(\sqrt{f}) + 2 \epsilon \tilde{C}(l, ||h_-||_{\infty})l B C_0 \right.\\
    \left.+ C_1 \epsilon + (C_2+C_3)\epsilon l +  {\sum_{j = 2}^{l}} \sum_{p=1}^{j-1} \binom{j-1}{p} \frac{2 C_4 l}{\epsilon^{p-1} N ^p}\right],
\end{align*} where $\tilde{C}(l, ||h_-||_{\infty}) = \frac{C(l, ||h_-||_{\infty})}{C(l)}$.
To conclude, we have that
\begin{align*}
    &\mathbb{E}_{\nu_N} \left[\Big | \int_0^t V_-(\epsilon,\eta)  {H}_s(0) ds \Big | \right]
    \\
    &\lsim  {\frac{l}{\epsilon N}} + C_1(l, ||h_-||_{\infty})\left[2 \epsilon l B + \epsilon + \epsilon l +  {\sum_{j = 2}^{l}} \sum_{p=1}^{j-1} \binom{j-1}{p} \frac{2 l}{\epsilon^{p-1} N ^p}\right] +  {\frac{1}{B}},
\end{align*}  {for some $C_1(l, ||h_-||_{\infty} > 0$,} as we wanted.  {We remark that the last expression} goes to zero as we take $N \to \infty$, then $\epsilon \to 0$, and finally $B \to \infty$.

\end{proof}

\begin{proof}[Proof of Corollary \ref{corollary_repl_lemma_2}]
Here we will only give the outline of the proof of \eqref{eq_1_corollary_RL2} and leave the details to the read and also the proof of \eqref{eq_2_corollary_RL2} since it is completely analogous.

If we do the same reasoning as at the beginning of the proof of Corollary \ref{corollary_repl_lemma_ave}, we can reduce the problem to just prove that, for any function $ {H} \in C^{1,2}( {\Omega_T})$ and every $t \in [0,T]$, it holds, for every $\delta >0$, that
\begin{equation*}
\lim_{\epsilon \to 0} \lim_{N\to+\infty}  \mathbb{P}_{\nu_N} \left( \Big | \int_0^t [\langle h_{-}(\eta_s),  {H}_s\rangle - F_-(\pi^N_s \ast \iota_\epsilon(\epsilon_N))  {H}_s(0)] ds \Big | > \delta \right) =0.
\end{equation*}
Then, applying Markov's inequality and using Lemma \ref{repl_lemma_2}, the result follows.
\end{proof}

\subsection{Orthonormal basis of \texorpdfstring{$L^2(\nu^N_\alpha)$}{TEXT}}

\begin{lemma} \label{lemma_basis_L^2}
Let us denote by $\nu^N_\alpha$ the Bernoulli product measure in $ {\{0,1\}^{\Sigma^-_l}}$ of parameter $\alpha \in (0,1)$. Denote, for every $x \in  {\Sigma^-_l}$, the $L^2(\nu_\alpha^N)$ functions
\begin{equation*}
    v_x(\eta) := \frac{\eta(x)-\alpha}{\alpha(1-\alpha)},
\end{equation*}
where $\eta \in  {\{0,1\}^{\Sigma^-_l}}$. Then, $\{v_x \ | \ x \in  {\Sigma^-_l}\}$ is an orthonormal basis of $L^2(\nu_\alpha^N)$.
In particular,
 {$\{\eta^{k_1}(\epsilon_N)\dots\eta^{k_l}(l\epsilon_N): k_j \in \{0,1\}, j \in \{1,\dots,l\}, \eta \in  {\{0,1\}^{\Sigma^-_l}}\}$} generates $L^2(\nu_\alpha^N)$ {, where here $\eta^{k_j}(z)$ represents $1$ if $k_j=0$ and $\eta(z)$ if $k_j=1$, for every $z \in  {\Sigma^-_l}$.}
\end{lemma}

\begin{proof}
Let $f: {\{0,1\}^{\Sigma^-_l}} \to \mathbb{R}$ be a $L^2(\nu_\alpha^N)$ function. Then, for every $\eta \in  {\{0,1\}^{\Sigma^-_l}}$,
\begin{equation} \label{eq_f_eta}
    f(\eta) = \sum_{\xi \in  {\{0,1\}^{\Sigma^-_l}}} f(\xi) \mathds{1}_{\xi = \eta} = \sum_{\xi \in  {\{0,1\}^{\Sigma^-_l}}} f(\xi) \prod_{x \in  {\Sigma^-_l}} \mathds{1}_{\xi(x) = \eta(x)}.
\end{equation}
Since, for every $x \in  {\Sigma^-_l}$, $\eta(x) \in \{0,1\}$, then, defining, for each $\xi \in  {\{0,1\}^{\Sigma^-_l}}$, the level sets
\begin{align} \label{def_A_xi_and_complement}
    A_\xi := \{ x \in  {\Sigma^-_l} \ | \ \xi(x) = 1\} \textrm{ and } A^c_\xi := \{ x \in  {\Sigma^-_l} \ | \ \xi(x) = 0\},
\end{align} and remarking that
\begin{equation*}
    \prod_{x \in A^c_\xi} [1-\eta(x)] = \sum_{B \subset A^c_\xi} (-1)^{|B|} \prod_{x \in B} \eta(x),
\end{equation*} we get that
\begin{equation*}
    \prod_{x \in  {\Sigma^-_l}} \mathds{1}_{\xi(x) = \eta(x)} = \prod_{x \in A_\xi} \eta(x) \prod_{x \in A^c_\xi} [1-\eta(x)] = \sum_{B \subset A^c_\xi} (-1)^{|B|} \prod_{x \in B \cup A_\xi} \eta(x).
\end{equation*}
Using the previous identity in \eqref{eq_f_eta}, we obtain that
\begin{equation} \label{equation_f_in_prod_eta}
     f(\eta) = \sum_{\xi \in  {\{0,1\}^{\Sigma^-_l}}} \sum_{B \subset A^c_\xi}  (-1)^{|B|} f(\xi) \prod_{x \in B \cup A_\xi} \eta(x).
\end{equation}
Now, observe that, for every $x \in  {\Sigma^-_l}$ and $\eta \in  {\{0,1\}^{\Sigma^-_l}}$,
\begin{equation*}
    \eta(x) = \alpha(1-\alpha) v_x(\eta) + \alpha,
\end{equation*}
therefore
\begin{align*}
    f(\eta) &= \sum_{\xi \in  {\{0,1\}^{\Sigma^-_l}}} \sum_{B \subset A^c_\xi}  (-1)^{|B|} f(\xi) \prod_{x \in B \cup A_\xi} [\alpha(1-\alpha) v_x(\eta) + \alpha] \\
    &= \sum_{\xi \in  {\{0,1\}^{\Sigma^-_l}}} \sum_{B \subset A^c_\xi}  (-1)^{|B|} f(\xi) \alpha^{|B \cup A_\xi|}\sum_{C \subset B \cup A_\xi} \prod_{x \in C} [(1-\alpha)v_x(\eta)] \\
    &= \sum_{\xi \in  {\{0,1\}^{\Sigma^-_l}}} \sum_{B \subset A^c_\xi} \sum_{C \subset B \cup A_\xi} (-1)^{|B|}\alpha^{|B|+|A_\xi|}(1-\alpha)^{|C|}f(\xi)\prod_{x \in C} v_x(\eta).
\end{align*}
Since $\mathbb{E}_{\nu^N_\alpha}(v_x v_y) = \mathds{1}_{y=x}$ due to the fact that $\nu^N_\alpha$ is of product form and $v_x$ is centered (in the sense that it has mean zero), for every $x \in  {\Sigma^-_l}$, we have then proved that $\{v_x \ | \ x \in  {\Sigma^-_l}\}$ is an orthonormal basis of $L^2(\nu_\alpha^N)$.
\end{proof}

\section{A Model with more than one stationary profile}\label{model_1.5}

One of the goals for the construction of these exclusion processes with non-reversible boundary dynamics was to capture an example of a family of particle systems for which its hydrodynamic equation has a unique solution but there exist multiple solutions to the stationary case. Below, we provide an explicit example of choices of boundary rates for which the model shows this property.

\subsection{The Exclusion $l3$ model} \label{model1_5}

Take $l=3$. Consider $0 < a < b$ and let $a_1 = a > 0$. Choose $a_2 \geq a + 2b > 0$ and $a_0 = 2 a_2 + 4b - a \geq a + 8b > 0$, and set
\begin{equation}
    (\mathscr{L}^{-}_N f)(\eta) = N c^L(\eta) [f(\sigma^2 \eta) - f(\eta)],
\end{equation} where
\begin{equation} \label{rates_boundary_model_l3_left}
    c^L(\eta) = a_2 \mathds{1}_{\eta(\epsilon_N) \neq \eta(3\epsilon_N)} + a_1 \mathds{1}_{\eta(\epsilon_N) = \eta(3\epsilon_N) = \eta(2\epsilon_N)} + a_0 \mathds{1}_{\eta(\epsilon_N) = \eta(3\epsilon_N) \neq \eta(2\epsilon_N)}
\end{equation} and 
\begin{equation}
    \sigma^2 \eta (z) = \begin{cases}
        \eta(z), \textrm{ if } z \neq 2\epsilon_N,\\
        1 - \eta(2\epsilon_N), \textrm{ otherwise}.
    \end{cases}
\end{equation} Analogously, set \begin{equation}
    (\mathscr{L}^{+}_N f)(\eta) = N c^R(\eta) [f(\sigma^{N-2} \eta) - f(\eta)],
\end{equation} where
\begin{equation} \label{rates_boundary_model_l3_right}
    c^R(\eta) = a_2 \mathds{1}_{\eta(\tau_N) \neq \eta(\tau_N - 2\epsilon_N)} + a_1 \mathds{1}_{\eta(\tau_N) = \eta(\tau_N - 2\epsilon_N) = \eta(\tau_N - \epsilon_N)} + a_0 \mathds{1}_{\eta(\tau_N) = \eta(\tau_N - 2\epsilon_N) \neq \eta(\tau_N-\epsilon_N)}
\end{equation} and 
\begin{equation}
    \sigma^{N-2} \eta (z) = \begin{cases}
        \eta(z), \textrm{ if } z \neq \tau_N -\epsilon_N,\\
        1 - \eta(\tau_N -\epsilon_N), \textrm{ otherwise}.
    \end{cases}
\end{equation}
 {Observe that the choice of boundary rates given in \eqref{rates_boundary_model_l3_left} (resp. \eqref{rates_boundary_model_l3_right}) are such that, in order to add/remove a particle at the site $2\epsilon_N$ (resp. $(N-2)\epsilon_N$), the system has to look at the occupation of the sites $\epsilon_N$ and $3\epsilon_N$ (resp. $(N-3)\epsilon_N$ and $\tau_N$). Depending on the configuration of the system on the sites $\epsilon_N, 2\epsilon_N$ and $3\epsilon_N$ (resp. $(N-3)\epsilon_N, (N-2)\epsilon_N$ and $\tau_N$), the rate at which a change happens at site $2\epsilon_N$ (resp. $(N-2)\epsilon_N$) changes. }

A simple computation shows that  {for this model}
\begin{equation*}
    h^-(\eta) = (1- 2\eta(2 \epsilon_N)) c^L(\eta)
\quad \textrm{ and } \quad
    h^+(\eta) = (1- 2\eta(\tau_N- \epsilon_N)) c^R(\eta).
\end{equation*} Thus
\begin{align} \nonumber
    F_{\mp}(\rho) &= \pm \frac{1}{4} \left[ (a_0 - 3 a_1 - 2 a_2)(2\rho - 1) - (a_0 + a_1 - 2 a_2)(2\rho - 1)^3\right] \\ \label{Fminus_plus}
    &= \pm \left[(b-a)(2\rho - 1) - b(2 \rho - 1)^3 \right].
\end{align}

Using Theorem \ref{th_uniqueness}, we have that the initial value problem \eqref{PDEproblem} with $F_\mp$ as defined in \eqref{Fminus_plus}  { and initial profile at time $t=0$ given by the measurable function $\gamma:[0,1] \to [0,1]$} has a unique solution in $L^2([0,T], H^1([0,1]))$.  {Nevertheless,} there exists more than one stationary solution of \eqref{PDEproblem} for our choice of $F_\pm$.  {Indeed, here} $F_+ = - F_-$, so it is enough to show that $F_-$ (or equivalently $F_+$) has more than one zero. As observed in Remark \ref{stationary_rmk}, each of $F_-$ (equivalently $F_+$) zeros is a constant which is also a stationary solution of \eqref{PDEproblem}. Then
\begin{equation}
    F_-(\rho) = (2 \rho - 1) [(b-a)- b (2 \rho - 1)^2] = 0 \Longleftrightarrow \rho \in \left\{\frac{1}{2}, \frac{1 \pm \sqrt{(b-a)/b}}{2}\right\} \cap [0,1].
\end{equation} Because there exists choices of $0 < a < b$ for which $F_-(\rho) = 0$ has three solutions, we are done.


 {We remark that this choice of boundary rates satisfies the hypothesis of Theorem \ref{main_result_large_dev}, i.e. for every $j \in \{1,2,3\}$, $B_j(\cdot)$ and $D_j(\cdot)$ are concave functions. Indeed, we will prove our previous statement for the case of $B_j$, since for $D_j$ it is analogous. For our choice of boundary rates, 
\begin{align*}
\mathfrak{b}^-(\alpha,M) &= a_1[(e^M-1)(1-\alpha)^3 + (e^{-M}-1)\alpha^3]\\
&+ a_0[(e^M-1)\alpha^2(1-\alpha) + (e^{-M}-1)\alpha(1-\alpha)^2]\\
&+ 2 a_2 [(e^M-1)\alpha(1-\alpha)^2 + (e^{-M}-1)\alpha^2(1-\alpha)],
\end{align*} for every $M \in \mathbb{R}$ and every $\alpha \in (0,1)$. Thus, by our choice of $a_0, a_1$ and $a_2$, we have that
\begin{equation*}
    B_1(\alpha) = a(1-\alpha)^3 + (2a_2+4b-a)\alpha^2(1-\alpha) + 2 a_2 \alpha(1-\alpha)^2\textrm{ and } B_2 = B_3 = 0,
\end{equation*} where $b>a>0$ and $a_2\geq a+2b$. Now, a simple computation shows that
\begin{align*}
    (B_1(\cdot))''(\alpha) = -4a_2+4a+8b-24b < 0,
\end{align*} because $b >0$ and $a_2\geq a+2b$. Then, $B_1$ is concave, and clearly $B_2$ and $B_3$ are also concave.}

\section{Initial-Value problems with non-linear Robin boundary conditions }\label{appB}

We present in this section some results on the problems \eqref{PDEproblem} and \eqref{LLN_perturbed_process}.
Let us consider the initial value problem given by
\begin{align}
\begin{cases}\label{nonlinear_robin}
\partial_t \rho_t(x) = \Delta \rho_t(x) -\alpha \nabla(2\sigma(\rho_t(x))\nabla H_t(x)), \textrm{ for } (t,x) \in  {\Omega_T},\\
\nabla \rho(t,0) - \alpha 2\sigma(\rho_t(0))\nabla H_t(0)= -\beta F^-(\rho(t,0)) -\alpha \mathfrak{p}^-(\rho_t(0),H_t(0)),\textrm{ for } t \in [0,T],\\
\nabla \rho(t,1) - \alpha 2\sigma(\rho_t(1))\nabla H_t(1) = \beta F^+(\rho(t,1))+ \alpha\mathfrak{p}^+(\rho_t(1),H_t(1)), \textrm{ for } t \in [0,T],\\
\rho(0,x) = \gamma(x), \textrm{ for } x \in [0,1].
\end{cases}
\end{align} where $\alpha,\beta \geq 0$, $\gamma:[0,1] \to \mathbb{R}$ is a measurable function,  {$\sigma:[0,1] \to [0,1]$ is given by $\sigma(x)=x(1-x)$} and $F^\pm: [0,1] \to \mathbb{R}$ are defined as
\begin{align*}
F^\pm(\alpha) = \mathbb{E}_{\nu_\alpha} \left[ \sum_{\xi \in  {\{0,1\}^{\Sigma^-_l}}} R^\pm(\eta,\xi) \sum_{x \in  {\Sigma^-_l}} \left[\xi(x) - \eta(x) \right]	\right],
\end{align*} with $\nu_\alpha$ a Bernoulli product measure in $ {\{0,1\}^{\Sigma^-_l}}$ with parameter $\alpha$ and $ {R^\pm(\cdot,\cdot) \geq 0}$. Also, $\mathfrak{p}^\pm:[0,1] \times \mathbb{R} \to \mathbb{R}$ are defined, for $0<a<1$ and $M\in\mathbb{R}$, as
\begin{equation*}
    \mathfrak{p}^{\pm}(a,M)=\mathbb{E}_{\nu_a}\left[\sum_{\xi\in {\{0,1\}^{\Sigma^\pm_l}}}R^{\pm}(\eta,\xi)M\sum_{y\in {\Sigma^\pm_l}}[\xi(y)-\eta(y)]\;e^{M\sum_{x\in {\Sigma^\pm_l}}[\xi(x)-\eta(x)]}\right].
\end{equation*}

\subsection{Uniqueness of weak solution}

 {The proof method we use here finds inspiration in the proof of Theorem 7.4 of Chapter 5 of \cite{PDEbook}  {and \cite{ThesisAngele}.}}
\begin{theorem} \label{th_uniqueness}
For every $\alpha, \beta \geq 0$, the initial value problem \eqref{nonlinear_robin} has a unique weak solution $\rho$ in $L^2([0,T];H^1([0,1]))$.
\end{theorem}

\begin{proof} \label{uniqueness_PDE_hydro}
Observe that, if $a,b \in [0,1]$, since $F_\pm$,  {$\sigma \in C^1([0,1])$}, then
\begin{align*}
    F_\pm(a) - F_\pm(b) = (a-b) (F_\pm)'( {\xi_{ab}}) \quad \textrm{ and } \quad
    \sigma(a) - \sigma(b) = (a-b) \sigma'(\gamma),
\end{align*} for  {some $\xi_{ab}$ and $ \gamma_{ab}$} between $a$ and $b$. Let
\begin{align*}
     \Tilde{F}_\pm(a,b) := \int_0^1 (F_\pm)'(s a + (1-s) b) ds \quad \textrm{ and } \quad
     \Tilde{\sigma}(a,b) := \int_0^1 \sigma'(s a + (1-s) b) ds,
\end{align*} and observe that $(F_\pm)'( {\xi_{ab}}) = \Tilde{F}_\pm(a,b)$ and $\sigma'( {\gamma_{ab}}) = \Tilde{\sigma}(a,b)$. Analogously, for every $M \in \mathbb{R}$, $\mathfrak{p}^{\pm}(\cdot,M) \in C^1([0,1])$ and so, if $a,b \in [0,1]$, then 
\begin{align*}
    \mathfrak{p}^{\pm}(a,M) - \mathfrak{p}^{\pm}(b,M) = (a-b) (\mathfrak{p}^{\pm}(\cdot,M))'( {\delta_{ab}}),
\end{align*} for $ {\delta_{ab}}$ between $a$ and $b$. Let
\begin{equation*}
     \Tilde{\mathfrak{p}}^{\pm}(a,b,M):= \int_0^1 (\mathfrak{p}^{\pm}(\cdot,M))'(s a + (1-s) b) ds,
\end{equation*} and observe that $(\mathfrak{p}^{\pm}(\cdot,M))'( {\delta_{ab}}) = \Tilde{\mathfrak{p}}^{\pm}(a,b,M)$. Let $\omega$ denote the difference between two solutions $\rho$ and $\Tilde{\rho}$, i.e.
\begin{equation*}
    \omega = \rho - \Tilde{\rho}.
\end{equation*}Then $\omega$ is  {a weak solution of}
\begin{align} \label{PDE_problem_linearized}
    \begin{cases}
    \partial_t \omega_t(x) = \Delta \omega_t(x) - \alpha \nabla(2\omega_t(x)\Tilde{\sigma}(\rho_t(x),\Tilde{\rho}_t(x)) \nabla H_t(x)), \textrm{ for } (t,x) \in  {\Omega_T},\\
    \nabla \omega_t (0) = - \beta\omega_t(0)\Tilde{F}_-(\rho_t(0),\Tilde{\rho}_t(0)) -\alpha \omega_t(0)\Tilde{\mathfrak{p}}^-(\rho_t(0),\Tilde{\rho}_t(0),H_t(0)), \textrm{ for } t \in [0,T],\\
    \nabla \omega_t (1) = \beta\omega_t(1)\Tilde{F}_+(\rho_t(1),\Tilde{\rho}_t(1)) +\alpha \omega_t(1)\Tilde{\mathfrak{p}}^+(\rho_t(1),\Tilde{\rho}_t(1),H_t(1)), \textrm{ for } t \in [0,T],\\
    \omega_0(x) = 0, \textrm{ for } x \in [0,T].
    \end{cases}
\end{align}
We want now to analyze this ``linearized'' problem.  {Consider the following problem}
\begin{align} \label{PDE_problem_linear}
    \begin{cases}
    \partial_t v_t(x) = \Delta v_t(x) - \nabla(v_t(x)f_t(x)), \textrm{ for } (t,x) \in  {\Omega_T},\\
    \nabla v_t (0) = - v_t(0)h_t(0), \textrm{ for } t \in [0,T],\\
    \nabla v_t (1) = v_t(1)h_t(1), \textrm{ for } t \in [0,T],\\
    v_0(x) = 0, \textrm{ for } x \in [0,T],
    \end{cases}
\end{align} where $$h_t(0) := \beta\Tilde{F}_-(\rho_t(0),\Tilde{\rho}_t(0)) + \alpha \Tilde{\mathfrak{p}}^-(\rho_t(0),\Tilde{\rho}_t(0),H_t(0)),$$ $$h_t(1) := \beta\Tilde{F}_+(\rho_t(1),\Tilde{\rho}_t(1)) + \alpha \Tilde{\mathfrak{p}}^+(\rho_t(1),\Tilde{\rho}_t(1),H_t(1)),$$ and $$f_t(x):= 2\alpha \Tilde{\sigma}(\rho_t(x),\Tilde{\rho}_t(x)) \nabla H_t(x).$$ 
 {Remark that, since we fixed $\rho$ and $\tilde{\rho}$, $h_t(0)$ and $h_t(1)$ can be seen as functions of only $t$ and analogously $f_t(x)$ can be seen as a functions of $x$ and $t$ only. In this sense, \eqref{PDE_problem_linear} is a linear problem.}
 {
\begin{definition}[Weak solution of \eqref{PDE_problem_linear}] \label{def_weak_sol_eq_linear_prob}
We say that $v$ is a weak solution of \eqref{PDE_problem_linear} if the following two conditions are satisfied:
\begin{enumerate}
    \item $v \in L^2([0,T], H^1([0,1]))$;
    \item For every $G \in C^{1,2}( {\Omega_T})$ and every $t \in [0,T]$, 
    \begin{equation} \label{def_weak_sol_eq}
    F(v,G,t) = 0,
    \end{equation} where
    \begin{align*}
        F(v,G,t) &:= \int_0^1 v_t(x) G_t(x)dx - \int_0^t \int_0^1 v_s(x) \partial_s G_s(x) dx ds + \int_0^t \int_0^1 \nabla v_s(x) \nabla G_s(x) dx ds \\
        &- \int_0^t \int_0^1 v_s(x) f_s(x) \nabla G_s(x) dx ds - \sum_{y \in \{0,1\}} \int_0^t v_s(y)[h_s(y)-f_s(y)] G_s(y) ds.
    \end{align*} 
\end{enumerate}
\end{definition}}
We claim that, for every fixed $\rho$ and $\Tilde{\rho}$, the problem \eqref{PDE_problem_linear} has a unique  {weak} solution in $L^2([0,T];H^1([0,1])$ which is $v = 0$. Since $\omega$ is a solution of \eqref{PDE_problem_linearized}, the previous argument implies that $\omega$ is unique and identically zero, which shows the uniqueness of the weak solution of the initial value problem \eqref{nonlinear_robin}.

 {We turn now to the proof of the uniqueness of weak solution for the linearized problem \eqref{PDE_problem_linear}}. We start by assuming that we are looking for a more regular solution, i.e. $v \in C^{1,2}( {\Omega_T})$, and proceed by an energy estimate. So, for every $t \in [0,T]$,
\begin{align} \label{what_we_should_improve}
    \frac{1}{2}\frac{d}{dt} ||v_t||_2^2 = - ||\nabla v_t||_2^2 +\frac{1}{2}\int_0^1 \nabla([v_t(x)]^2)f_t(x)dx + \sum_{y \in \{0,1\}} (v_t(y))^2 [-f_t(y) + h_t(y)].
\end{align}
Now note that, for every $t \in [0,T]$, applying Young's inequality, we have that, for $A>0$ and $a \in (0,1)$,
\begin{align*}
(v_t(0))^2 &\leq \left( 1 + \frac{1}{A}\right) \left(\frac{1}{a}\int_0^a [v_t(0)-v_t(x)]dx \right)^2 + \frac{ 1 + A}{a^2} \left(\int_0^a v_t(x)dx \right)^2.
\end{align*}
Using the Fundamental Theorem of Calculus and Fubini's Theorem, we can rewrite the RHS of the last display as
\begin{align*}
\frac{1}{a^2}\left( 1 + \frac{1}{A}\right) \left(\int_0^a [\nabla v_t(y)](a-y) dy\right)^2 + \frac{1 + A}{a^2} \left(\int_0^a v_t(x)dx \right)^2.
\end{align*} Applying Cauchy-Schwartz's inequality on both terms, we arrive at
\begin{align} \label{poincare1}
    (v_t(0))^2 \leq \frac{a}{3}\left( 1 + \frac{1}{A}\right) \int_0^1 [\nabla v_t(y)]^2 dy + \frac{ 1 + A}{a} \int_0^1 [v_t(x)]^2 dx.
\end{align}
Analogously, for $B > 0$,
\begin{align} \label{poincare2}
    (v_t(1))^2 \leq \frac{a}{3}\left( 1 + \frac{1}{B}\right) \int_0^1 [\nabla v_t(y)]^2 dy + \frac{ 1 + B}{a} \int_0^1 [v_t(x)]^2 dx.
\end{align}
Using \eqref{poincare1}, \eqref{poincare2} and Young's inequality, we can rewrite the RHS of  \eqref{what_we_should_improve} to conclude that $\frac{1}{2}\frac{d}{dt} ||v_t||_2^2$ is  {bounded} from above by a constant times
\begin{align*}
    &\left( \frac{a}{3}\left( 1 + \frac{1}{B}\right) C(1) + \frac{a}{3}\left( 1 + \frac{1}{A}\right) C(0) - 1 + A'\right)||\nabla v_t||_2^2 \\
    &+ \left( \frac{ 1 + B}{a} C(1) + \frac{ 1 + A}{a} C(0) + \frac{\max_{x \in [0,1]} |f_t(x)|^2}{A'}\right)||v_t||_2^2,
\end{align*} where, for $x \in \{0,1\}$,
\begin{equation}
    C(x) := \max\left\{ \max_{t \in [0,T]} [h_t(x)-f_t(x)], 0 \right\}.
\end{equation}
Since $F^+$ and $F^-$ are $C^1[0,1]$ and, for every $M \in \mathbb{R}$, $\mathfrak{p}^{+}(\cdot,M)$ and $\mathfrak{p}^{-}(\cdot,M)$ are $C^1[0,1]$, then there exists $C > 0$ such that, for $x \in \{0,1\}$,  {$C(x) \leq C$. Moreover, since $\sigma \in C^1[0,1]$ and $\sup_{t \in [0,T]} |\nabla H_t(x)| < \infty$ for $x \in \{0,1\}$, then there exists $C > 0$ (possibly different than the previous one) such that $\sup_{t \in [0,T]} \max_{y \in [0,1]} |f_t(y)|^2 \leq C$.} This implies that, taking $a$ and $A'$ sufficiently small depending on the chosen values of $A$ and $B$, for every $t \in [0,T]$ there exists $C_1 > 0$ and $C_2 \geq 0$ such that
\begin{equation*}
    \frac{d}{dt} ||v_t||_2^2 \leq -C_1 ||\nabla v_t||_2^2 + C_2 ||v_t||_2^2.
\end{equation*} By Gronwall's inequality, we conclude that, for every $t \in [0,T]$,
\begin{equation}
    ||v_t||^2_2 \leq ||v_0||^2_2 e^{C_2 t}.
\end{equation} Because, $v_0 = 0$, we can immediately conclude that, for all $t \in [0,T]$, $v_t = 0$, as we wanted to show. 

\begin{remark}
This last argument not only shows the uniqueness of the solution in $C^{1,2}( {\Omega_T})$ but also the dependence of the solution on the initial states.
\end{remark}

Now we treat the case with lower regularity.  {We want to construct from the solution that we know to exist in $L^2([0,T], H^1([0,1])$ a more regular one for which we could apply similar arguments as in the case of a $C^{1,2}(\Omega_T)$ solution to guarantee uniqueness. To do that, we follow the ideas used in Section 5.3 of \cite{Nahum2020} and explore the reduction of the initial problem to a linear problem where it makes sense to define mild solutions. We start by defining what is a mild solution of \eqref{PDE_problem_linear}.}

\begin{definition}
    We call mild solution of \eqref{PDE_problem_linear} to any function $v: {\Omega_T}\to [0,1]$ satisfying $M(v,t) := v_t -  {(Sv)_t} = 0$ a.s. with
    \begin{equation}
         {(S v)_t} (\cdot):= \int_0^t \left[ P_{t-s}(\cdot,0)v_s(0)h_s(0) - P_{t-s}(\cdot,1)v_s(1)h_s(1) - \int_0^1 P_{t-s}(\cdot,x)\nabla(v_s(x)f_s(x)) dx \right]ds ,
    \end{equation} where $P_t(u,v) = \sum_{w \in  {\psi}^{-1}(v)} \Phi_t(u,w)$ is the density kernel generated by the Laplacian $\Delta$ on $[0,1]$ with reflecting Neumann boundary conditions and is related to the heat kernel
    \begin{equation}
        \Phi_t(u,w) = \frac{1}{(4\pi t)^{1/2}} e^{-\frac{(u-w)^2}{4t}},
    \end{equation} by the reflection map $\psi:\mathbb{R} \to [0,1]$ defined as 
    \begin{equation}
        \psi(u+k) = \begin{cases}
            u, \quad u \in [0,1], k \textrm{ even,}\\
            1-u, \quad u \in [0,1], k \textrm{ odd,}
        \end{cases}
    \end{equation} extended to $\mathbb{R}$ by symmetry as $\psi(v) = \psi(-v)$, for $v \in \mathbb{R}$. 
\end{definition}

The idea to introduce mild solutions here is to  {take advantage of} the  {regularity} of $S$, which is a consequence of the regularizing effect of $P_\cdot$, to conclude the uniqueness of the weak solution of the initial value problem \eqref{PDE_problem_linear}.  {We remark that $Sv$ is differentiable in time and $C^\infty$ in space since $P_t(\cdot,v)$ is a $C^\infty$ function.}

Our goal is to relate weak solutions of \eqref{PDE_problem_linear} with mild solutions of \eqref{PDE_problem_linear} with an analogous to Proposition 5.10 of \cite{Nahum2020} and from here conclude the uniqueness.  {The next proposition creates exactly such relation.}

\begin{proposition} \label{regularity_of_Sv}
For any mild solution $v: {\Omega_T} \to [0,1]$ of \eqref{PDE_problem_linear} that belongs to $L^2([0,T],H^1([0,1]))$, we have that $ {v} \in H^{1}([0,T], C^2([0,1]))$.  {Moreover, if $\rho: {\Omega_T} \to \mathbb{R}$ is a function satisfying $\langle M(\rho,t),G \rangle = 0$ for any $G \in C^{1,2}( {\Omega_T})$, then $F(S \rho,G,t) = 0$ for any $G \in C^{1,2}([0,T]\times[0,1])$ and every $t \geq 0$, and so $S \rho$ is a (strong) solution of \eqref{PDE_problem_linear}.}
\end{proposition}

\begin{proof}
Let us suppose we have a mild solution $v: {\Omega_T} \to [0,1]$ of \eqref{PDE_problem_linear} that belongs to $L^2([0,T],H^1([0,1]))$. We start by observing that the heat kernel is smooth in time and also smooth in both variables of space. Moreover, because $v \in L^2([0,T], H^1([0,1]))$, we also have by Lebesgue's differentiation theorem that $\partial_t S v_t \in L^2([0,T], C^{\infty}([0,1]))$  with 
\begin{align*}
    \partial_t S v_t &= v_t(0)h_t(0) - v_t(1)h_t(1) - \int_0^1 \nabla(v_t(x)f_t(x)) dx \\
    &+ \int_0^t \left[ \partial_t P_{t-s}(\cdot,0)v_s(0)h_s(0) - \partial_t P_{t-s}(\cdot,1)v_s(1)h_s(1) - \int_0^1 \partial_t P_{t-s}(\cdot,x)\nabla(v_s(x)f_s(x)) dx \right]ds,
\end{align*} in the almost sure sense. Thus, $S v_t \in H^1([0,T], C^{\infty}([0,1])) \subset H^1([0,T], C^2([0,1]))$, and this proves the  {first} statement of the proposition.  {The proof of the second statement follows as Proposition 5.10 of \cite{Nahum2020}.}
\end{proof}

In the next proposition, we observe that every weak solution is also a mild solution, in the same spirit as Proposition 5.10 of \cite{Nahum2020}.

\begin{proposition} \label{auxiliary_prop_for_uniqueness_without_reg}
All functions $v:  {\Omega_T} \to [0,1]$ that are weak solution of \eqref{PDE_problem_linear} in the sense of Definition \ref{def_weak_sol_eq_linear_prob} are also mild solutions of \eqref{PDE_problem_linear}. 
\end{proposition}
The proof of this result follows by adapting the proof of Proposition 5.10 of \cite{Nahum2020},  {and for that reason we omit it here}.

Now, we know that $\omega$ is a weak solution of \eqref{PDE_problem_linear}, then, by Proposition \ref{auxiliary_prop_for_uniqueness_without_reg}, $\omega$ is also a mild solution. Remark that $C^{1,2}( {\Omega_T})$ is a dense subset of $H^1([0,T], C^2([0,1]))$, thus, by standard arguments, we also have that $F(\omega,G,t) = 0$ for every $G \in H^1([0,T],C^2([0,1]))$ and $t \in [0,T]$. On the other hand, by definition of mild solution, we have that $\omega_t =  {(S \omega)}_t$ a.s.. Then, Proposition \ref{regularity_of_Sv} implies that in fact $\omega \in H^1([0,T],C^2([0,1]))$ and so $0 = F(\omega, \omega, t)$. We can write this equality has 
\begin{align*}
0 &= \int_0^1 [\omega_t(x)]^2dx - \int_0^t \int_0^1 \omega_s(x) \partial_s \omega_s(x) dx ds + \int_0^t \int_0^1 [\nabla \omega_s(x)]^2 dx ds \\
&- \int_0^t \int_0^1 \omega_s(x) f_s(x) \nabla \omega_s(x) dx ds - \sum_{y \in \{0,1\}} \int_0^t [\omega_s(y)]^2[h_s(y)-f_s(y)]ds,
\end{align*} for every $t \geq 0$. We observe that
\begin{align*}
\int_0^1 [\omega_t(x)]^2dx - \int_0^t \int_0^1 \omega_s(x) \partial_s \omega_s(x) dx ds 
&=  \int_0^t \int_0^1 \omega_s(x) \partial_t \omega_s(x) dx ds \\
&= \int_0^t \frac{1}{2}\lim_{h \to 0} \int_0^1 \frac{[\omega_{s+h}(x)]^2 - [\omega_{s}(x)]^2}{h} dx ds,
\end{align*} where the last identity follows by the Dominated Convergence Theorem and the fact that, since $\omega \in H^1([0,T],C^2([0,1]))$, by the Sobolev embeddings, we have that $\omega \in C([0,T],C^2([0,1]))$. So furthermore,
\begin{align*}
\int_0^1 [\omega_t(x)]^2dx - \int_0^t \int_0^1 \omega_s(x) \partial_s \omega_s(x) dx ds 
&= \int_0^t \frac{1}{2} \frac{d}{ds} \int_0^1 [\omega_{s}(x)]^2 dx ds,
\end{align*} and thus
\begin{align*}
0 &= \int_0^t \frac{1}{2} \frac{d}{ds} \int_0^1 [\omega_{s}(x)]^2 dx ds + \int_0^t \int_0^1 [\nabla \omega_s(x)]^2 dx ds \\
&- \int_0^t \int_0^1 \omega_s(x) f_s(x) \nabla \omega_s(x) dx ds - \sum_{y \in \{0,1\}} \int_0^t [\omega_s(y)]^2[h_s(y)-f_s(y)]ds,
\end{align*} for every $t \geq 0$. By Lebesgue's differentiation theorem, we obtain that 
\begin{align*}
\frac{1}{2} \frac{d}{dt} \int_0^1 [\omega_{t}(x)]^2 dx &= -\int_0^1 \left\{[\nabla \omega_t(x)]^2 -  \omega_t(x) f_t(x) \nabla \omega_t(x) \right\}dx + \sum_{y \in \{0,1\}} [\omega_t(y)]^2 [h_t(y)-f_t(y)],
\end{align*} a.s. in $t$. Then because $\omega_t$ is a $C^2$ function in space, we can repeat the arguments we did in the uniqueness in the space $C^{1,2}$ to bound the RHS of the last display and obtain that 
\begin{align*}
\frac{d}{dt} \int_0^1 [\omega_{t}(x)]^2 dx &\leq - C_1 \int_0^1 [\nabla \omega_{t}(x)]^2 dx + C_2 \int_0^1 [\omega_{t}(x)]^2 dx
\end{align*}a.s. in $t$. 
 {Since} $\omega$ is continuous in time, by the integral form of Gronwall's inequality it holds that
\begin{align*}
\int_0^1 [\omega_{t}(x)]^2 dx \leq \int_0^1 [\omega_{0}(x)]^2 dx \ e^{ C_2 t} = 0
\end{align*}a.s. in $t$, and so $\omega = 0$ a.s. in $t$ and $x$. With this, we showed that $\rho_t = \Tilde{\rho}_t$ a.e for all $t \geq 0$, which completes the proof of the uniqueness of the weak solution of \eqref{PDE_problem_linear} in the sense of definition \eqref{def_weak_sol_eq_linear_prob}. 

\begin{remark}
    Here is a summary with the main ideas behind the proof strategy of the uniqueness of the initial value problem \eqref{PDEproblem}:
    \begin{enumerate}
        \item Fix two possible solutions, consider the difference of the two and write the initial value problem associated to the difference;
        
        \item Consider the linearized problem associated to the one solved by the function of the difference and construct $\Phi: C \times C \to M$ as $\Phi(\rho,\Tilde{\rho}) = v$ where $C$ is the set of  {weak} solutions of \eqref{PDEproblem} and $M$ is the set of  {weak} solutions of the linearized problem. 
        
        \item Start by considering more regularity of your solutions, i.e. that are  {$C^{1,2}( {\Omega_T})$} or, in other words, that are in the same space as the space of test functions in the definition of weak solution and show uniqueness in  {$C^{1,2}( {\Omega_T})$} for the linearized problem. Extend the uniqueness to $L^2([0,T],H^1([0,1]))$ with the help of the definition of mild solutions.
        
        \item Conclude by observing that, for each fixed $\rho,\Tilde{\rho}$ because the linearized problem has a  {weak} solution which is unique and identically equal to zero, we guarantee that $\Phi$ is well defined and, in fact, it is the zero function in $L^2([0,T],H^1([0,1]))$. Because, for each fixed $\rho, \Tilde{\rho}$, we show that $\omega = \rho - \Tilde{\rho}$ is solution to the linearized problem, then $\omega = \Phi(\rho,\Tilde{\rho}) = 0$, proving the uniqueness of the non-linear problem.
    \end{enumerate}
\end{remark}
\end{proof}
\subsection{Properties of the solution of the hydrodynamic equation}
Let us consider the initial value problem given by \eqref{nonlinear_robin} with $\alpha=0$ and $\beta =1$.  { In this section, we fix $\gamma:[0,1]\to[0,1]$, which will be the initial condition at time zero, to be a $C^{2+\epsilon}([0,1])$ function with $0 < \epsilon < 1$ which satisfies the boundary conditions $\nabla \gamma(0) = - F_-(\gamma)$ and $\nabla \gamma(1) = F_+(\gamma(1))$.}
\begin{lemma} \label{thm_general_heateq} The problem \eqref{nonlinear_robin} has a weak solution in $L^2([0,T], H^1([0,1]))$. Moreover, if it has a solution in $C^{1,2}( {\Omega_T})$, it is unique. The uniqueness of solution extends to weak solutions in $L^2([0,T], H^1([0,1]))$.
\end{lemma}

\begin{lemma}\label{solution_is_smooth}
 {The problem \eqref{nonlinear_robin} has a smooth solution in $(0,T]\times[0,1]$.}
\end{lemma}
\begin{proof}
 {Observe that \eqref{nonlinear_robin} can be identified with the initial value problem given by (7.1)-(7.3) of Section 4 of Chapter 5 of \cite{PDEbook} taking
\begin{align*}
    a_{ij} = \mathbbm{1}_{i=j}, \quad b(x,t,u, u_x) = 0, \quad \psi_0 = \gamma, \quad \textrm{ and } \quad 
\psi(x,t,u) = \begin{cases} F^-(u_t(x)), \textrm{ if } x = 0,\\
- F^+(u_t(x)), \textrm{ if } x = 1.
\end{cases}
\end{align*}
To apply Theorem 7.4 of Chapter 5 of \cite{PDEbook} we would need to have $\gamma = 0$. Since this is not the case, consider $v = \rho - \gamma$. Now, we observe that we can write an initial value problem for $v$ also of the form of the one given by (7.1)-(7.3) of Section 4 of Chapter 5 of \cite{PDEbook} where now the initial condition at time $t$ is zero, as we wanted. 
Thus, applying \cite[Chapter 5, Theorem 7.4]{PDEbook} to $v$ and then recovering the connection between $v$ and $\rho$, we get that there exists a weak solution $\rho$ of \eqref{nonlinear_robin} in $H^{1+\beta/2, 2+\beta}( {\Omega_T})$ for some $\beta > 0$, that in our case is strictly bigger than $1$.} Since all functions in  {$H^{1+\beta/2, 2+\beta}( {\Omega_T})$ are continuous with all derivatives of the form $\partial^r_t \partial^s_x$ with $0 \leq r < (\beta+1)/2$ and $0 \leq s < 3/2 + \beta$ also continuous and in our case $\beta > 1$,} the solution $\rho$ is in $C^{1,2}( {\Omega_T})$. By a bootstrap argument, one can conclude that the solution is smooth in $(0,T]\times[0,1]$.
\end{proof}
\begin{lemma}\label{weak_max_principle}
The problem \eqref{nonlinear_robin} satisfies a weak maximum principle.
\end{lemma}
\begin{proof}
Let $\rho$ be the solution of the problem in $C^{1,2}( {\Omega_T})$.
We want to show that
\begin{align*}
\max_{(t,x)\in  {\Omega_T}} \rho_t(x) = \max_{(t,x) \in S^1\cup S^2 \cup S^3} \rho_t(x),
\end{align*} where 
\begin{align*}
S^1 := \{(t,0) \ | \ t \in [0,T]\}; \ S^2 := \{(t,1) \ | \ t \in [0,T]\}; \ S^3 := \{(0,x) \ | \ x \in (0,1)\}.
\end{align*}
To show this, it is enough to prove that $\rho$ can not attain its maximum in $(0,T]\times(0,1)$.  {Suppose this is the case, i.e. there exists $(x_0,t_0) \in (0,T]\times(0,1)$ such that $$\max_{(t,x)\in  {\Omega_T}} \rho_t(x) = \rho_{t_0}(x_0) > \max_{(t,x) \in S^1\cup S^2 \cup S^3} \rho_t(x).$$ Then, since $\rho \in [0,1]$, by Theorem $2.$ of \cite{MaximumPrinciples}, we have that $\rho_t(x) = \rho_{t_0}(x_0)$ for all $(t,x) \in  {\Omega_T}$  {and $(t,x)$ can} be connected by a horizontal and a vertical line segment included in $ {\Omega_T}$ to the point $(t_0,x_0)$. This is a contradiction because, for example, $(t_0,1)$ is a point in the conditions above and by assumption $$\rho_{t_0}(1) \leq \max_{(t,x) \in S^1\cup S^2 \cup S^3} \rho_t(x) < \rho_{t_0}(x_0).$$}

Repeating the same argument for $-\rho$ we get that
\begin{align*}
\min_{(t,x)\in  {\Omega_T}} \rho_t(x) = \min_{(t,x) \in S^1\cup S^2 \cup S^3} \rho_t(x),
\end{align*} and we are done.
\end{proof}

\begin{lemma}\label{bounded_away_0_1}
    For all $0<t_0\leq T$, there exists an $\varepsilon>0$ such that $\varepsilon\leq \rho_t(x)\leq1-\varepsilon$ for all $(t,x)\in[t_0,T]\times[0,1]$, where $\rho$ is unique solution of \eqref{nonlinear_robin}.
\end{lemma}
\begin{proof}
By Lemma \ref{weak_max_principle}, we have that, for every $\delta > 0$,
\begin{align*}
\max_{(t,x)\in [\delta,T]\times[0,1]} \rho_t(x) = \max_{(t,x) \in S_\delta^1\cup S_\delta^2} \rho_t(x),
\end{align*} where 
\begin{align*}
&S_\delta^1 := \{(t,0) \ | \ t \in [\delta,T]\}; \  S_\delta^2 := \{(t,1) \ | \ t \in [\delta,T]\},
\end{align*} { and the analogous replacing $\max$ by $\min$.} Then, if 
\begin{align*}
\max_{(t,x) \in S_\delta^1\cup S_\delta^2} \rho_t(x) = \rho_{t_0}(0)
\end{align*} for some $t_0 \in [\delta,T]$, then  {$\nabla \rho_{t_0}(0) \leq 0$ and so, using the boundary condition $\nabla \rho_{t_0}(0) = - F^-(\rho_{t_0}(0))$, $F^-(\rho_{t_0}(0)) \geq 0$}, i.e.
\begin{align*}
\mathbb{E}_{\nu_{\rho_{t_0}(0)}} \left[ \sum_{\xi \in  {\{0,1\}^{\Sigma^-_l}}} R^-(\eta,\xi) \sum_{x \in  {\Sigma^-_l}} \left[\xi(x) - \eta(x) \right]	\right]  { \geq } 0.
\end{align*}
 {On the other hand, if we had that $\rho_{t_0}(0) = 1$, then, denoting by $\mathbbm{1}$ the configuration in $ {\{0,1\}^{\Sigma^-_l}}$ with a particle in every point $x \in  {\Sigma^-_l}$, we have that
\begin{align*}
F^-(\rho_{t_0}(0)) = \mathbb{E}_{\nu_{1}} \left[ \sum_{\xi \in  {\{0,1\}^{\Sigma^-_l}}} R^-(\eta,\xi) \sum_{x \in  {\Sigma^-_l}} \left[\xi(x) - 1 \right]	\right] = \sum_{\xi \in  {\{0,1\}^{\Sigma^-_l}} \setminus \{\mathbbm{1}\}} R^-(\mathbbm{1},\xi) \sum_{x \in  {\Sigma^-_l}} \left[\xi(x) - 1 \right] < 0,
\end{align*} which is a contradiction. For the last inequality we used that $R^-(\mathbbm{1},\mathbbm{1}) = 0$ and that, for every $\xi \in  {\{0,1\}^{\Sigma^-_l}} \setminus \{\mathbbm{1}\}$, $R^-(\mathbbm{1},\xi) > 0$.
From here we conclude that if there exists $t_0 \in [\delta,T]$ such that
\begin{align*}
\max_{(t,x) \in S_\delta^1\cup S_\delta^2} \rho_t(x) = \rho_{t_0}(0),
\end{align*} then $\rho_{t_0}(0) < 1$.
Adapting the previous argument, we can easily conclude that, if there exists $t_1 \in [\delta,T]$ such that
\begin{align*}
\min_{(t,x) \in S_\delta^1\cup S_\delta^2} \rho_t(x) = \rho_{t_1}(1),
\end{align*} then $\rho_{t_1}(1) > 0$. Using Theorem 2 of \cite{MaximumPrinciples} we deduce that it can not exist $t \in [\delta, T]$ such that $\max_{(t,x) \in S_\delta^1\cup S_\delta^2} \rho_t(x) = \rho_t(1)$ with $\rho_t(1) = 1$ nor $\min_{(t,x) \in S_\delta^1\cup S_\delta^2} \rho_t(x) = \rho_t(0)$ with $\rho_t(0) = 0$. Putting together the previous results and recalling the maximum principle  - Theorem \ref{weak_max_principle} - we can conclude that} there exists $\alpha>0$ and $\beta < 1$ such that, for every $(t,x) \in [\delta,T]\times[0,1]$,
\begin{equation}
0 < \alpha \leq \rho_t(x) \leq \beta < 1,
\end{equation} which finishes our proof.
\end{proof}

\begin{lemma}\label{our_B5}
    There exists a finite constant $C_0$ that depends on $R^{\pm}$ and $l$ such that 
    \begin{equation}
    \begin{split}
    &\int_0^t\int_0^1\frac{(\nabla\rho_s)^2}{\sigma(\rho_s)}dx\;ds+\int_0^t\Big|F_{-}(\rho_s(0))\log\Big(\frac{\rho_s(0)}{1-\rho_s(0)}\Big)\Big|\;ds\\
    &+\int_0^t\Big|F_{+}(\rho_s(1))\log\Big(\frac{\rho_s(1)}{1-\rho_s(1)}\Big)\Big|\;ds\leq C_0 t+\int_0^1 V_0(\gamma)\;dx-\int_0^1 V_0(\rho_t)\;dx,
    \end{split}
    \end{equation}
    where $V_0(x)=x\log(x)+(1-x)\log(1-x).$
\end{lemma}
 {The last} result is analogous to \cite[Lemma B.5]{FGLN23}.  {It} relies on the fact that the function
\begin{equation*}
    f_{\pm}(\alpha):=F_{\pm}(\alpha)\log\Big(\frac{\alpha}{1-\alpha}\Big)
\end{equation*}
is bounded from below by a positive constant, which is a consequence of the fact that $F_{\pm}\in C^1([0,1])$ and by Lemma \ref{bounded_away_0_1}.

\section*{Acknowledgement}

C. L. has been partially supported by FAPERJ CNE E-26/201.117/2021, by CNPq Bolsa de Produtividade em Pesquisa PQ 305779/2022-2. B. S. thanks FCT/Portugal for the financial support through the PhD scholarship with reference 2022.13270.BD. B. S. thanks Thyago Santos for the fruitful conversations. J. M. thanks CNPq for the financial support through the PhD scholarship.

\section*{Declarations}

{}

The authors have no competing interests to declare that are relevant
to the content of this article.

Data sharing not applicable to this article as no datasets were
generated or analysed during the current study.

\end{document}